\documentclass[reqno,11pt]{amsart}
\usepackage{geometry}
\usepackage{mathtools,amssymb,amsthm,mathrsfs,color,lineno,paralist,graphicx,float}
\usepackage[colorlinks,
linkcolor=blue,
anchorcolor=green,
citecolor=blue, 
]{hyperref}
\usepackage{tikz}
\usepackage{graphicx}
\usepackage{longtable}
\usepackage{array}
\usepackage{booktabs}
\usepackage{comment}
\usetikzlibrary{shapes.geometric, arrows.meta, positioning, fit, backgrounds, shadows}
\usepackage{amsmath,amssymb}
\usetikzlibrary{arrows.meta}

\usepackage{enumitem}
\usepackage[T1]{fontenc}
\usepackage[utf8]{inputenc}
\usepackage{subfig} 
\usepackage[justification = centering, labelsep =period]{caption} 
\usepackage[mathcal]{eucal}

\usepackage{calc}
\definecolor{bleu1}{RGB}{0,57,128}
\def\bleu1{\color{bleu1}}

\usepackage{etoolbox}
\patchcmd{\section}{\normalfont}{\normalfont \bleu1}{}{}
\patchcmd{\subsection}{\normalfont}{\normalfont \bleu1}{}{}
\patchcmd{\subsubsection}{\normalfont}{\normalfont \bleu1}{}{}

\usepackage{tikz}
\usetikzlibrary{arrows.meta}

\newtheorem{proposition}{Proposition}[section]
\newtheorem{theorem}{Theorem}[section]
\newtheorem{definition}{Definition}[section] 
\newtheorem{lemma}{Lemma}[section]
\newtheorem{remark}{Remark}[section]
\newtheorem{corollary}{Corollary}[section]

\newtheorem{problem}{Problem}
\newtheorem{conjecture}{Conjecture}

\newcommand{\Z}{{\mathbb Z}}
\newcommand{\C}{{\mathbb C}}
\newcommand{\R}{{\mathbb R}}
\newcommand{\Q}{{\mathbb Q}}
\newcommand{\T}{{\mathbb T}}
\newcommand{\pos}[1]{(#1)_+}

\def\0{{\mathbf 0}}

\usepackage{tikz,tikz-3dplot}
\tdplotsetmaincoords{80}{45}
\tdplotsetrotatedcoords{-90}{180}{-90}
\tikzset{surface/.style={draw=blue!70!black, fill=blue!40!white, fill opacity=.6}}

\newcommand{\tangentlabel}[5]{
    \pgfmathsetmacro{\ang}{#3}
    \pgfmathsetmacro{\rot}{\ang+90}
    \pgfmathparse{ifthenelse(\rot>90,\rot-180,ifthenelse(\rot<-90,\rot+180,\rot))}
    \let\finalrot\pgfmathresult
    \node[
        #1,
        rotate=\finalrot,
        anchor=south
    ] at ({(#2+#5)*cos(\ang)},{(#2+#5)*sin(\ang)}) {#4};
}

\usepackage{pgfplots}
\usepgfplotslibrary{polar}
	\usepgfplotslibrary{polar}
\pgfplotsset{compat=1.17}
\makeatletter
\tikzset{reuse path/.code={\pgfsyssoftpath@setcurrentpath{#1}}}
\makeatother

\usepackage{soul}

\begin{document}
	\title[Acceleration and Gap Collapse]{Acceleration and Gap Collapse: Effective Dispersion in Quasiperiodic Schr\"odinger Operators} 

\author{Xianzhe Li}
\address{Department of Mathematics, University of California, Berkeley, CA 94720, USA} 
 \email{xianzhe@berkeley.edu}

 \author{Duxiao Wang}
\address{Chern Institute of Mathematics and LPMC, Nankai University, Tianjin 300071, China}
\email{2542721099@qq.com}

\author{Zhenfu Wang}
\address{
	Chern Institute of Mathematics and LPMC, Nankai University, Tianjin 300071, China
}
\email{zhenfuwang@mail.nankai.edu.cn}

\author{Qi Zhou}
\address{
	Chern Institute of Mathematics and LPMC, Nankai University, Tianjin 300071, China
}
\email{qizhou@nankai.edu.cn}

\begin{abstract}
We answer Avila's acceleration problem \cite{Avila15} for Diophantine frequencies.
Consider analytic perturbations of $2\lambda\cos(2\pi Nx)$, with
$N\geqslant2$ and $\lambda>1$. For each direction outside a finite
union of proper closed real linear hyperplanes, every sufficiently
small positive amplitude realizes all accelerations $1,\ldots,N$,
each on a spectral set of positive Lebesgue measure. The coexistence of
acceleration strata persists on open sets of potentials, disproving
the two whole-spectrum Type~I density conjectures of
Ge--Jitomirskaya--You \cite{GJY26}.
For the same
directions, on a set of amplitudes of density one at $0^+$, the
perturbed potential has fundamental period one and infinitely many
collapsed gap labels in every residue class modulo $N$.
In particular, this gives a negative answer to  Goldstein's question \cite{G16} on the finiteness of collapsed gaps. 

The infinite-range Aubry
dual  has quadratic leading dispersion near the split
outer edge and linear dispersion with nonzero slopes at the
collapsed energies we construct. Near the edge we obtain uniform
joint asymptotics for the integrated density of states and the
real and complexified Lyapunov exponents, including in spectral
gaps. 

The proof reduces the slow dynamics of the infinite-range 
dual to a $2N$-dimensional analytic cocycle. We use its Hermitian
symplectic structure, monotonicity in energy, and first variation
with respect to the potential to derive the two dispersion laws.
\end{abstract}

\maketitle

\section{Introduction}

How do distinct types of motion emerge from a degenerate Hamiltonian
system and coexist as parameters vary? Near a parabolic configuration,
perturbations can separate neutral modes into oscillatory and
exponentially growing or decaying directions. Near a degenerate
elliptic configuration, the question is how oscillatory modes split
and whether their interaction preserves or destroys stability.
Classical Hamiltonian stability theory and the theory of matrix
unfoldings describe these alternatives
\cite{Moser1958,Arnold1971,Galin1982}.
Parameter-dependent KAM theory studies their persistence under
quasiperiodic perturbations, including configurations with multiple
normal frequencies \cite{BHTB90,BHN07}.

For quasiperiodic Schr\"odinger operators, energy plays the role of a bifurcation parameter, whose effect on the dynamics is fixed by the operator. Can a typical perturbation produce distinct dynamical types? Can this coexistence persist throughout a neighborhood of potentials?  Can open and
collapsed spectral gaps coexist after perturbations break the
degeneracy? 
Persistent coexistence would also obstruct the density of operators whose spectrum has only one dynamical type. We study these questions in the regime of positive Lyapunov exponent, and intend to show how the
number and the phases of neutral dual modes govern these  spectral
phenomena.

\subsection{Avila's acceleration problem}

Avila's global theory was a breakthrough in the study of analytic one-frequency Schr\"odinger operators \cite{Avila15}:
\begin{equation}\label{eq:operator}
(H_{v,\alpha,x}u)_n
=u_{n+1}+u_{n-1}+v(x+n\alpha)u_n,
\qquad n\in\Z.
\end{equation}
 It gave a framework for arbitrary analytic potentials and irrational frequencies, distinguishing  the spectral energies into subcritical, critical, and supercritical regimes and describing their variation with energy and potential. 
 The common spectrum is denoted by $\Sigma_{v,\alpha}$. The eigenvalue equation gives the analytic $\mathrm{SL}(2,\R)$ cocycle
$$
S_E^v(x)=\begin{pmatrix}E-v(x)&-1\\1&0\end{pmatrix}.
$$
Write $L_{v,\alpha}(E;y)$ for its Lyapunov exponent at the complexified phase $x+iy$. On the spectrum, an energy is supercritical if $L_{v,\alpha}(E;0)>0$, subcritical if $L_{v,\alpha}(E;y)=0$ on a neighborhood of $y=0$, and critical otherwise. The function $L_{v,\alpha}(E;y)$ is even and convex in $y$.  In particular, the Lyapunov exponent is analytic on appropriate strata even when it is not analytic across their boundaries. A central ingredient is acceleration: 
\begin{equation*}\label{eq:omega}
\omega_v(E)=\frac{1}{2\pi}\partial_y^+ L_{v,\alpha}(E;0)
\end{equation*}
is a nonnegative integer. This discrete invariant makes it possible to study transitions between regimes and to distinguish dynamics within the positive-exponent region.

The almost Mathieu potential $v(x)=2\lambda\cos(2\pi x)$, $\lambda>0$, is the basic example.  Introduced as a model for an electron on a two-dimensional lattice in
a uniform magnetic field \cite{Pe,Ha,R}, this operator has become a
basic model at the intersection of spectral theory, dynamical systems,
and mathematical physics. All its spectral energies are subcritical for $\lambda<1$, critical for $\lambda=1$, and supercritical for $\lambda>1$ \cite{Avila15}. This model has a sharp spectral description, either in the spectral measure side \cite{AYZ,Jit99,JL18,Jit21}, or in the spectrum geometry side \cite{AJ09,AYZ23},  while its single harmonic forces the same acceleration at every spectral energy.

The higher harmonic $2\lambda\cos(2\pi Nx)$, with $\lambda>1$, has the spectral data of the almost Mathieu operator at frequency $N\alpha$, and acceleration $N$ throughout its spectrum. Avila \cite{Avila15} asked whether a typical perturbation breaks this uniform behavior into all the intermediate accelerations:

\begin{problem}\label{prob:Avila}
Fix an integer $N\ge2$, and $\lambda>1$. For a typical real analytic direction $W$, does every sufficiently small positive perturbation
\begin{equation*}\label{eq:family}
v_{\varepsilon,W}(x)=2\lambda\cos(2\pi Nx)+\varepsilon W(x)
\end{equation*}
have a positive-Lebesgue-measure set of spectral energies of each acceleration $1,\ldots,N$?
\end{problem}

We use \emph{typical} in a measure-theoretic sense. A Borel set is
prevalent if some compactly supported probability measure assigns
full measure to every translate of it \cite{HSY92}. Our exceptional
sets satisfy the stronger condition of being Gauss-null: they have
zero measure under every nondegenerate Gaussian measure \cite{BL}.
Their complements are therefore prevalent. We  answer Avila's question for Diophantine frequency $\alpha$,  we recall that an irrational frequency $\alpha$ is Diophantine if there are $\gamma>0$ and $\tau>1$ such that
$$
\|n\alpha\|_{\mathbb{T}}\geqslant\gamma|n|^{-\tau},\qquad n\in\mathbb{Z}\setminus\{0\},
$$
where $\|x\|_{\mathbb{T}}=\inf_{j\in\mathbb{Z}}|x-j|$. We denote the set of these frequencies by $DC$.  Then our main result is as follows:

\begin{theorem}\label{thm:strata}
Fix $\alpha\in DC$, $N\geqslant2$, $\lambda>1$, and $r>0$. There exists a set $\mathcal{G}\subset C_{r}^\omega(\mathbb T,\mathbb R)$ whose complement is a finite union of proper closed real linear hyperplanes with the following property. For every $W_{*}\in\mathcal{G}$, there exist an open neighborhood $\mathcal{U}\subset\mathcal{G}$ of $W_{*}$ and constants $\varepsilon_{0},c_{1},\ldots,c_{N}>0$ such that
\begin{equation}\label{eq:positive-measure}
\operatorname{Leb}\big\{E\in\Sigma_{v_{\varepsilon,W},\alpha}:\omega_{v_{\varepsilon,W}}(E)=k\big\}\geqslant c_{k}\varepsilon
\end{equation}
for all $W\in\mathcal{U}$, $0<\varepsilon<\varepsilon_{0}$, and $k=1,\ldots,N$. In particular, $\mathcal{G}$ is open and dense and has Gauss-null complement.
\end{theorem}

While this paper was in preparation, we learned of work by
Liang--Wang--Xu \cite{LWX26}, who also establish positive-measure
coexistence of all acceleration strata.

\subsection{Whole-spectrum Type I density conjectures}
\label{subsec:typeI-density}
Avila's problem also leads naturally to a complementary question:
how common are operators whose entire spectrum is of Type~I?
This question motivates the density conjectures of
Ge--Jitomirskaya--You \cite{GJY26}.  
We first recall the T-acceleration introduced by
Ge--Jitomirskaya--You \cite{GJY26}.  As $L_{v,\alpha}(E,\eta)$ is 
even and piecewise-affine function, let $\eta_1(E)\ge0$ be the first turning point of $L_{v,\alpha}(E,\eta)$, with
$\eta_1(E)=0$ when $\omega_v(E)>0$.
The T-acceleration is defined by
$$
\overline\omega_{v,\alpha}(E)
=
\lim_{\delta\downarrow0}
\frac{
L_{v,\alpha}(E,\eta_1(E)+\delta)
-
L_{v,\alpha}(E,\eta_1(E))
}{2\pi\delta}.
$$
If there is no turning point in the domain of analyticity,
we set $\overline\omega_{v,\alpha}(E)=0$.
An energy is \emph{Type~I} if
$\overline\omega_{v,\alpha}(E)=1$, and an operator is
Type~I if every spectral energy is Type~I.

The Type~I condition captures a feature shared by all
nonzero almost Mathieu potentials and provides a way to
study general analytic operators through their dual dynamics.
For trigonometric polynomial potentials, quantitative
duality identifies Type~I with simplicity of the lowest nonnegative
dual Lyapunov exponent \cite{GJYZ26,GJY26}.
On the supercritical spectrum, $
\overline\omega_{v,\alpha}(E)=\omega_v(E),$
so Type~I means acceleration one.
This viewpoint has led to robustness of both versions
of the Ten Martini Problem and Dry Ten Martini Problem  \cite{GJY26,LXZ26}.

These results raise a natural question: how common are
operators whose entire spectrum is Type~I?
Ge--Jitomirskaya--You  \cite[Conjecture~1(2)--(3)]{GJY26} proposed the following two density
conjectures.

\begin{conjecture}\leavevmode
\label{conj:whole-spectrum-typeI}
\begin{enumerate}
\item
For every irrational $\alpha$ and every $r>0$, the potentials
whose associated operators are Type~I, with this property
persisting throughout a $ C_{r}^\omega(\mathbb T,\mathbb R)$-neighborhood, form a dense
subset of $ C_{r}^\omega(\mathbb T,\mathbb R)$.
\item
For some $r>0$, the potentials admitting an
$ C_{r}^\omega(\mathbb T,\mathbb R)$-neighborhood on which the Type~I property holds
uniformly over all irrational frequencies form a dense
subset of $ C_{r}^\omega(\mathbb T,\mathbb R)$.
\end{enumerate}
\end{conjecture}

Theorem~\ref{thm:strata} gives an obstruction to both assertions.

\begin{corollary}
\label{cor:non-typeI-open}
For every $\alpha\in DC$ and every $r>0$,
there exist a nonempty open set $\mathcal U\subset C_{r}^\omega(\mathbb T,\mathbb R)$
and a constant $c>0$ such that
$$
\operatorname{Leb}
\left\{
E\in\Sigma_{v,\alpha}:\omega_v(E)\ge2
\right\}
\ge c,
\qquad v\in\mathcal U.
$$
In particular, both assertions of
Conjecture~\ref{conj:whole-spectrum-typeI} are false.
\end{corollary}

Avila's acceleration problem and these density conjectures
concern two aspects of the same perturbative question.
The former asks whether distinct acceleration values
coexist on spectral sets of positive measure; the latter
asks whether operators with Type~I behavior throughout
the spectrum can nevertheless be dense.
Theorem~\ref{thm:strata} connects the two questions:
coexistence along an open set of perturbation directions
yields an open-set obstruction to whole-spectrum
Type~I density.

However, a weaker density question remains open.
Conjecture~1(1) of \cite{GJY26} asks whether, for an open
dense set of potentials, Type~I energies are relatively
open and dense in the spectrum.
Our results do not settle this assertion.
Positive-measure strata of higher acceleration can coexist
with a relatively open dense set of Type~I energies.

\subsection{Gap labels and collapsed gaps}
\label{subsec:gap-labels-collapsed-gaps}

The coexistence of several oscillatory branches leads to a further
question about the gaps: which labels are realized by open gaps?
By the gap-labelling theorem \cite{BLT,JM82}, on each bounded
spectral gap the integrated density of states $\mathcal N_{v,\alpha}$
is constant and equal to $\{\ell\alpha\}$ for a unique
$\ell\in\mathbb Z\setminus\{0\}$. Here $\{\cdot\}$ denotes fractional
part. For each such label, set
$$
 \mathcal I_\ell(v,\alpha)
 :=
 \{E\in\mathbb R:\mathcal N_{v,\alpha}(E)=\{\ell\alpha\}\}.
$$
 If it has positive length, its interior is the
gap with label $\ell$. We call $\ell$ \emph{collapsed} if
$\mathcal I_\ell(v,\alpha)$ is a singleton.

The Ten Martini Problem \cite{Simon} asks whether the spectrum of the almost
Mathieu operator is a Cantor set for every nonzero coupling and
every irrational frequency. It was solved by Avila and Jitomirskaya
\cite{AJ09}. For general analytic potentials, Goldstein and Schlag
\cite{GS11} proved density of gaps in the positive Lyapunov regime
for almost every frequency. Ge--Jitomirskaya--You \cite{GJY26}
established Cantor spectrum for analytic Type~I operators at every
irrational frequency. They subsequently removed the Type~I assumption, proving
Cantor spectrum whenever the Lyapunov exponent is positive
throughout the spectrum \cite{GJY26+}.
Han--Schlag--Yang \cite{HSY26} give an alternative proof
for every fixed Diophantine frequency, using finite-volume
resonances and regularity of the IDS.

The Cantor property concerns the density of gaps. The Dry Ten
Martini Problem \cite{Simon} asks for more: every allowed label
must be realized by an open gap. In the magnetic interpretation of
the almost Mathieu operator, the gap labels correspond to the TKNN
integers that determine the quantized Hall conductance \cite{TKNN}.
Avila--You--Zhou \cite{AYZ23} solved the problem at every
noncritical coupling. For the critical case, one can consult 
 \cite{AYZ26,BF}.

For analytic potentials, every supercritical Type~I energy
satisfying the gap-labelling condition is an endpoint of an open gap
by Li, Xu, and Zhou \cite{LXZ26}. The corresponding subcritical
result was proved by Ge \cite{Ge26} and Li \cite{Li26}. Outside
Type~I, Argentieri and Avila \cite{AA} have announced analytic
perturbations of subcritical almost Mathieu potentials whose
spectra contain intervals. In the supercritical regime, spectral
intervals are ruled out by \cite{GJY26+,HSY26}, individual labels
may still collapsd, actually Puig and
Sim\'o \cite{PS11} found numerical evidence that open and
collapsed gaps coexist.
What mechanism allows some gaps to collapse?

 On the other hand, early back to September 2016, during Q. Zhou's visit to M. Goldstein, Goldstein asked whether large coupling leaves at most finitely many collapsed gaps:

\begin{problem}
\label{prob:goldstein-large-coupling}
Fix $\alpha\in\mathbb R\setminus\mathbb Q$ and a strip width $r>0$.
For a generic $f\in C_r^\omega(\mathbb T,\mathbb R)$ of fundamental
period one, does there exist $\lambda_0(f)>0$ such that
$H_{\lambda f,\alpha,x}$ has at most finitely many collapsed labels
for every $\lambda\geqslant\lambda_0(f)$?
\end{problem}

In the following, we answer  M. Goldstein's question in the sense of prevalence and  show that infinitely many labels can collapse. To state a precise version of the result, recall
a Borel set of a Banach space $X$ is \emph{shy} if some compactly supported
Borel probability measure on $X$ assigns measure zero to every
translate of the set \cite{HSY92}. Set $v_{0}(x)=2\lambda\cos(2\pi Nx)$, and let $\mathcal{S}_{N}$ be the set of potentials in $C_{r}^{\omega}(\mathbb{T},\mathbb{R})$ with fundamental period one and infinitely many collapsed gap labels in each residue class $\ell\equiv a\pmod{N}$, $0\leqslant a<N$.

\begin{theorem}\label{thm:collapsed-gaps}
Fix $\alpha\in DC$, $N\geqslant2$, $\lambda>1$, and $r>0$.
There exists an open dense set $\mathcal G_{\mathrm{cl}}\subset C_r^\omega(\mathbb T,\mathbb R)$,
whose complement is a finite union of proper closed real linear
hyperplanes, such that, for every $W\in\mathcal G_{\mathrm{cl}}$,
$$
 \operatorname{Leb}\{\varepsilon\in(0,\delta):
          v_{\varepsilon,W}\notin\mathcal S_N\}
 =O_W(\delta^3)
 \qquad (\delta\downarrow0).
$$
Moreover, $\mathcal U\cap\mathcal S_N$ is non-shy in $C_r^\omega(\mathbb T,\mathbb R)$
for every open neighborhood $\mathcal U$ of $v_0$.
For every such $\mathcal U$ and every prescribed positive
sequence $\varepsilon_*^j\to0$, the set
\begin{equation}\label{uniform}
    \bigl\{f\in\mathcal U:
       (\varepsilon_*^j)^{-1}f\in\mathcal S_N
       \text{ for infinitely many }j\bigr\}
\end{equation}
is also non-shy in $C_r^\omega(\mathbb T,\mathbb R)$.
\end{theorem}

The class $a=0$ is particularly revealing. At $\varepsilon=0$,
the model is the almost Mathieu operator at frequency $N\alpha$,
so every nonzero label in $N\mathbb Z$ corresponds to an open gap
by \cite{AYZ23}. The other labels are collapsed because of the
smaller period of the sampling function. Thus the theorem closes infinitely many gaps
already present in the unperturbed spectrum. Moreover, infinitely many collapsed gaps can occur at arbitrarily
large couplings for the same sampling function. In particular,
the eventual-finiteness property in Problem \ref{prob:goldstein-large-coupling}
does not hold on a prevalent set of sampling functions.
When the result was finished, we learned a supercritical Type~II
counterexample to the dry statement was announced 
\cite{GWXZ}.

Theorem~\ref{thm:collapsed-gaps} gives infinitely many collapsed
labels at each retained amplitude. We now fix the label in advance.
Even after the perturbation breaks the $1/N$-periodicity of $v_0$,
each initially collapsed label remains collapsed for all sufficiently
small amplitudes in a common set of density one at zero. The smallness
threshold may depend on the label.

\begin{theorem}
\label{thm:persistence-collapsed-labels}Fix $\alpha\in DC$, $N\geqslant2$, $\lambda>1$, and $r>0$.
There exists a set
$\mathcal G_{\mathrm{all}}\subset C_r^\omega(\mathbb T,\mathbb R)$
whose complement is a countable union of proper closed real linear
hyperplanes with the following property.
For every $W\in\mathcal G_{\mathrm{all}}$, there is a Borel set
$\mathcal E_W\subset(0,1)$ satisfying
$$
  \operatorname{Leb}\bigl((0,\delta)\setminus\mathcal E_W\bigr)
  =O_W(\delta^2)
  \qquad (\delta\downarrow0).
$$
For each $m\in\mathbb Z\setminus N\mathbb Z$, there exists
$\varepsilon_m(W)>0$ such that, for every
$\varepsilon\in\mathcal E_W\cap(0,\varepsilon_m(W))$, the label $m$
is collapsed for $H_{v_{\varepsilon,W},\alpha,x}$.
\end{theorem}

As a consequence, for $W\in\mathcal G_{\mathrm{all}}$, any prescribed finite set of
nonzero labels retains its unperturbed gap status for all sufficiently
small $\varepsilon\in\mathcal E_W$. The nonzero labels in $N\mathbb Z$
have open gaps by Avila--You--Zhou \cite{AYZ23}, with sharp gap-width
asymptotics due to Ge--You--Zhou \cite{GYZ24}. Each such gap persists
under sufficiently small perturbations. For labels outside $N\mathbb Z$,
persistence follows from Theorem~\ref{thm:persistence-collapsed-labels}.

\subsection{Effective  dispersion and joint spectral asymptotics}
\label{subsec:dispersion-landscape}

The preceding results establish coexistence and its stability.
We next ask where the transitions occur and how the spectral
data change across them. Our answer lies in a family of $N$ local dispersion
branches that meet at the unperturbed outer almost Mathieu edge
and separate under perturbation. Their separation organizes these
transitions through the local energy-quasimomentum  dispersion. 

For periodic operators, the dispersion relation between energy
and the Floquet exponent is central to band theory in solid-state
physics, linking spectral structure to wave propagation and state
counting \cite{Bloch29}. Quasiperiodicity removes the translation
symmetry behind this description, although Bloch waves remain
available in suitable regimes \cite{Eliasson92}.
Lesser and Lifshitz \cite{LL22} asked whether an
\textquotedblleft energy-momentum dispersion curve, or effective
band structure\textquotedblright\
could be associated with quasiperiodic Bloch waves.
They emphasized that such a description could help interpret
experimental signatures of dispersion in quasicrystals, as well
as observations that weak disorder can enhance conductivity and
increase the spatial extent of wave functions.
We obtain a local version of this picture for the analytic Aubry
dual, together with its consequences for the spectral data of
the original operator.

Fix the parameters $\alpha,N,\lambda,r$ of Theorem~\ref{thm:strata}.
Let
$$
 E_0=\sup\Sigma_{v_0,\alpha},
 \qquad E_\varepsilon(\xi)=E_0-\varepsilon\xi .
$$
Write $\Sigma_{\varepsilon,W}=\Sigma_{v_{\varepsilon,W},\alpha}$.
We set $L_{\varepsilon,W}$ and $\mathcal N_{\varepsilon,W}$
for the Lyapunov exponent and the IDS of the original
Schr\"odinger operator, and $\widehat H_{\varepsilon,x}$ for
its Aubry dual, defined by \begin{equation*}\label{eq:dual-operator}
 \bigl(\widehat H_{\varepsilon,x}u\bigr)_j
 =
 \lambda(u_{j+N}+u_{j-N})
 +2\cos\bigl(2\pi(x+j\alpha)\bigr)u_j
 +\varepsilon\sum_{q}w_qu_{j+q},
\end{equation*}
A real number $k$ is a Bloch quasimomentum at energy $E$ if
the dual equation admits a nonzero analytic block Bloch solution
\begin{equation}\label{eq:intro-block-Bloch}
 u_{mN+s}(x)=e^{imNk}F_s(x+mN\alpha),
 \qquad m\in\mathbb Z,\quad 0\le s<N,
\end{equation}
with $F\in C^\omega(\mathbb T,\mathbb C^N)$.
Thus $e^{iNk}$ is the block Floquet multiplier and $k$ is its
Floquet exponent normalized per lattice site.  For $W\in\mathcal G$,
there exist
$\xi_1(W)<\cdots<\xi_N(W)$ which are determined by $W$\footnote{These transition values are defined in \eqref{eq:transition-values}.}. The constant
$\kappa>0$ below depends only on the unperturbed operator.
We write $s_{+}=\max\{s,0\}$ and take $\sqrt{s}$ to be nonnegative for $s\ge0$ and positive
imaginary for $s<0$.

\begin{theorem}\label{thm:landscape}
	Under the assumption of Theorem \ref{thm:strata}, for every $W_{*}\in\mathcal{G}$ and nondegenerate compact interval $I\subset\mathbb{R}\setminus\{\xi_{1}(W_{*}),\ldots,\xi_{N}(W_{*})\}$, there are $\delta,K>0$ and
	$\varepsilon_{0}=\varepsilon_{0}(\alpha,N,\lambda,r,W_{*},I)>0$
	such that the following hold for $\|W-W_{*}\|_{r}<\delta$ and $0<\varepsilon<\varepsilon_{0}$. 
	\begin{enumerate}
		\item[\textup{(i)}]
		If $I\subset(\xi_{1}(W_{*}),\infty)$, there is a compact set $T_{\varepsilon,W}\subset I$ satisfying $\operatorname{Leb}(I\setminus T_{\varepsilon,W})\leqslant K\varepsilon^{\frac{5}{2}}.$	For $\xi\in T_{\varepsilon,W}$, one has $E_{0}-\varepsilon\xi\in\Sigma_{\varepsilon,W}$ and
		\begin{equation}\label{eq:intro-acceleration-count}
			\omega_{\varepsilon,W}(E_{0}-\varepsilon\xi)=\#\{j:\xi_{j}(W)<\xi\}.
		\end{equation}
		For each $j$ with $\xi_{j}(W)<\min I$, the dual equation admits analytic block Bloch solutions with positive quasimomenta $k_{\varepsilon,j}(W,\xi)$ satisfying
		\begin{equation}\label{eq:intro-dispersion}
			E_{0}-\varepsilon\xi=E_{0}-\varepsilon\xi_{j}(W)-\kappa^{-2}k_{\varepsilon,j}(W,\xi)^{2}+O(\varepsilon^{\frac{3}{2}}),
		\end{equation}
	uniformly for $\|W-W_{*}\|_{r}<\delta$ and $\xi\in T_{\varepsilon,W}$.
		\item[\textup{(ii)}]
		As $\varepsilon\to0$, the following estimate holds uniformly for $\|W-W_{*}\|_{r}<\delta$, $\xi\in I$, and $\eta$ in bounded subsets of $\mathbb{R}$:
		\begin{equation}\label{eq:intro-uniform-landscape}
			\begin{aligned}
				&\quad L_{\varepsilon,W}(E_{0}-\varepsilon\xi;\sqrt{\varepsilon}\eta)+i\pi\big[1-\mathcal{N}_{\varepsilon,W}(E_{0}-\varepsilon\xi)\big]\\
				&=\log\lambda+\sqrt{\varepsilon}\sum_{j=1}^{N}\max\big\{\kappa\sqrt{(\xi_{j}(W)-\xi)_{+}},2\pi|\eta|\big\}+i\kappa\sqrt{\varepsilon}\sum_{j=1}^{N}\sqrt{(\xi-\xi_{j}(W))_{+}}+O(\varepsilon).
			\end{aligned}
		\end{equation}

	\end{enumerate}
\end{theorem}

 Here the
coincident edge modes separate into branches with distinct
thresholds and the same leading quadratic term in
quasimomentum, as expressed by \eqref{eq:intro-dispersion}.
The joint estimate \eqref{eq:intro-uniform-landscape} gives a
quantitative local refinement of the multiplicative Jensen
formula \cite{GJYZ26}; see also \cite[Theorem~1.9]{HS26Pi}
for an independent proof in the trigonometric polynomial case.
It resolves the dual growth rates into explicit dispersion
branches. Their competition with the complexification scale
determines the Lyapunov profile, while their oscillatory phases
determine the IDS.
At $\eta=0$, this takes the form
\begin{equation*}
\begin{aligned}
 L_{\varepsilon,W}(E_\varepsilon(\xi);0)
 +i\pi\bigl[1-\mathcal N_{\varepsilon,W}(E_\varepsilon(\xi))\bigr]
 =\log\lambda+
 \kappa\sqrt\varepsilon\sum_{j=1}^N\sqrt{\xi_j-\xi}
 +O(\varepsilon),
\end{aligned}
\end{equation*}
where the square root is positive imaginary for a negative
argument. After subtracting $\log\lambda$, the leading real
part sums the decay rates of the nonoscillatory branches,
and the leading imaginary part sums the oscillatory phases,
giving the spectral mass above $E_\varepsilon(\xi)$.
The same square-root branch describes both: crossing its
threshold changes its contribution from real to imaginary.
Thus Lyapunov growth and spectral mass are determined by
a common effective dispersion relation. On $T_{\varepsilon,W}$, acceleration
counts the oscillatory branches.
The breakpoints of the limiting complexified Lyapunov profile
record the individual dual decay rates. All these transitions
occur while the original real-axis Lyapunov exponent remains
close to $\log\lambda$.

The quadratic edge dispersion has a linear counterpart at the
collapsed energies, where all $2N$ neutral dual modes have
nonzero velocities.

\begin{corollary}
\label{cor:linear-dispersion}
The infinitely many collapsed energies in
Theorem~\ref{thm:collapsed-gaps} may be chosen so that the Aubry dual
has $2N$ first-order Bloch branches at each such energy $E_\ell$.
The branch at quasimomentum $k_j$ has a real analytic energy
$E_j(k)$ satisfying
\begin{equation}\label{eq:collapsed-linear-dispersion}
 E_j(k)=E_\ell+s_j(k-k_j)+O\bigl(|k-k_j|^2\bigr),
 \qquad s_j\ne0,\qquad j=1,\ldots,2N.
\end{equation}
Here a first-order Bloch branch means that the ansatz
$u_m(k,x)=e^{imk}f_j(k,x+m\alpha)$, with $f_j$ analytic and
$f_j(k_j,\cdot)\not\equiv0$, solves the dual equation at energy
$E_j(k)$ up to an error $O(|k-k_j|^2)$, uniformly in $x$ and $m$.
\end{corollary}

The two mechanisms have different local dispersion laws.
At the outer edge, the common quadratic term produces square-root
phases on the oscillatory side and square-root growth rates on
the hyperbolic side. Quantitative duality converts these branch
data into the IDS and the real and complexified Lyapunov
profiles. At the constructed collapsed energies, the first-order
Bloch branches have nonzero slopes. Their inverse absolute
velocities, with the rotation-number normalization, contribute
positively to the derivative of the IDS. The same slow dynamics
thus describes both the onset of oscillation and the strict
spectral crossing at a collapsed label.

\subsection{Main ingredients and local dynamics}
\label{sec:main-ingredients}

Our approach rests on a dynamical interpretation of acceleration.
In the positive Lyapunov exponent regime, acceleration is related
to the asymptotic distribution of zeros of finite-volume Dirichlet
determinants \cite{HS26,WWYZ24}.
Avila \cite{Avila15} suggested relating it to the number of dynamical
critical points, with collisions producing intermediate
accelerations. Liang--Wang--Xu \cite{LWX26} develop this viewpoint
through the real zeros of overlap coefficients for first-return
cocycles. For our purposes, the essential step is quantitative
global theory \cite{GJYZ26}, which identifies acceleration with
the number of neutral pairs in the Aubry dual. Changes in
acceleration can thus be studied through changes in the dual
dynamics.

Building on canonical center dynamics \cite{CCB26}, we construct
a fixed-dimensional slow cocycle that preserves the
Hermitian-symplectic form. Energy monotonicity
\cite{WXZ,LXZ26} then controls the motion of its modes.
The two coexistence problems concern the number of neutral
pairs and the relations among their phases.
Successive collisions in a multiple parabolic unfolding change
the number of elliptic pairs and hence the acceleration.
For several elliptic pairs, relative phase variation allows
a prescribed gap label to be realized while avoiding pairwise
quasiperiodic resonances. This gives the mechanism for gap collapse.

\subsubsection{Quantitative duality and a uniform slow cocycle}
\label{sec:slow}

Aubry duality links the quasiperiodic Schr\"odinger family
\eqref{eq:operator} to the long-range family 
\begin{equation*}\label{infi-1}
 (\widehat H_{v,\alpha,\theta}u)_n
 =
 \sum_{k\in\mathbb Z}\widehat v_k\,u_{n+k}
 +2 \cos 2 \pi (\theta+n\alpha)u_n.
\end{equation*}
When $v$ is a trigonometric polynomial of degree $m$,
the dual dynamics admits a $2m$-dimensional
Hermitian-symplectic cocycle $(\alpha,\mathcal A_E)$ \cite{HP13}.\footnote{Consult Section \ref{sec:aubry-duality} for details.}
Quantitative global theory \cite{GJYZ26} relates
acceleration to its partially hyperbolic splitting:
\begin{equation*}\label{eq:qgt-splitting}
 \mathbb C^{2m}=E^s\oplus E^c\oplus E^u,
 \qquad
 \dim E^c=2\omega(E).
\end{equation*}
The center accounts precisely for the zero Lyapunov
exponents, counted with multiplicity. Acceleration
therefore measures the dimension of the neutral
dual dynamics.

For general analytic potentials, the dual equation
has infinite range, and the finite-dimensional
transfer-matrix representation is no longer available.
A breakthrough of Ge--Jitomirskaya \cite{GJ26}
is the construction of an intrinsic center cocycle
for infinite-range dual operators. By combining
holomorphic trivialization via the Oka--Grauert
principle \cite{Forstneric2017} with convergence of finite-range
approximants, they recover a finite-dimensional
description of the neutral dynamics.

These results explain the dynamical meaning of the acceleration.
To follow a branch as it enters or leaves the neutral
dynamics, we need a bundle that contains it on both sides
of the transition. We therefore construct a fixed-dimensional
slow bundle containing the center and the directions that
become weakly hyperbolic.

For $0<a<2\pi r$, set the weighted space\footnote{The role of this weight was first recognized in \cite{LXZ26}, then the weighted space was formally introduced in \cite{CCB26}.}
$$
 X_a=
 \left\{u:\sum_{j\in\mathbb Z}|u_j|e^{-a|j|}<\infty\right\},
$$
and define the shift on sequences by
$
(Su)_n := u_{n+1},  n\in\Z.
$
The slow cocycle is obtained from translation on these
weighted solution bundles.

\begin{theorem}\label{thm:slow}
For any $\alpha\in \R\backslash \Q$, fix $W\in  C_{r}^\omega(\mathbb T,\mathbb R)$, $0<a<2\pi r$, and $E_*\in \Sigma_{v_0,\alpha}$. There exist $\varepsilon_0,h>0$ such that,
for $|\varepsilon|,|E-E_*|<\varepsilon_0$, the following hold.
\begin{enumerate}
\item
The weighted solution spaces $\mathcal M_{\varepsilon,E}(x)
 =
 \ker_{X_a}\bigl(\widehat H_{\varepsilon,x}-E\bigr) $
form a globally analytically trivial bundle of dimension
$2N$.
\item
There is a global frame $U_{\varepsilon,E}(x)$, jointly analytic in
$(x,\varepsilon,E)$,
such that
$$
 S^{N}U_{\varepsilon,E}(x)
 =U_{\varepsilon,E}(x+N\alpha)C_{\varepsilon,E}(x).
$$
The induced slow cocycle $(N\alpha,C_{\varepsilon,E})$ is
Hermitian-symplectic for real $\varepsilon,E$.
\end{enumerate}
\end{theorem}

A precise version is given in Proposition~\ref{prop:slow1} and Proposition~\ref{prop:slow};
the first-order response to energy and potential is given in
Proposition~\ref{energyvar}.
The construction develops the transfer-operator method of
canonical center bundles introduced by
Li--Wang--You--Zhou \cite{CCB26}. Their reference problem is
a translation-invariant hopping operator with a small
on-site perturbation,
\begin{equation}\label{infi}
 (\widehat H_{v,w,\alpha,\theta}u)_n
 =
 \sum_{k\in\mathbb Z}\widehat v_k\,u_{n+k}
 +\epsilon w(\theta+n\alpha)u_n.
\end{equation}
The zeros of the free hopping symbol determine a splitting
in weighted spaces; transfer operators and their spectral
projections continue the corresponding solution dynamics.
This replaces companion matrices by canonical center dynamics for the
infinite-range operator \eqref{infi}.
The smallness required by this reference problem is not
available throughout the present coupling range. Dividing
the dual equation by $\lambda$ leaves an on-site term of
size $2/\lambda$, which need not be perturbatively small
when $\lambda>1$. We instead include this term in the
reference operator and consider $\widehat H_{\varepsilon,x}$ as a perturbation of $\widehat H_{0,x}$, while $\widehat H_{0,x}$ corresponding to $N$-folded subcritical almost Mathieu cocycle.

The slow bundle is defined by a growth condition that applies
uniformly at every $E_*\in\Sigma_0$. Subcriticality of the dual
almost Mathieu cocycle gives uniform subexponential bounds on its
forward and backward iterates. On contours inside and outside the
unit circle, exponential weights dominate this growth and give
uniform inverses for the operator-valued symbol. These inverses
persist under small analytic perturbations, uniformly in the
truncation range (Lemma~\ref{lem:inverseD}). The contour projection
of \cite{CCB26} then selects the full space of solutions in $X_a$.
Its rank remains $2N$, and translation on its range defines the
intrinsic slow cocycle (Proposition~\ref{prop:slow1}).

Translation preserves the symplectic form, whose restriction to the
slow bundle remains nondegenerate. Normalizing this form gives
Hermitian-symplectic coordinates in which to follow the slow
multipliers (Lemma~\ref{lem:sympD}). The distinction between the
two applications enters only at this stage. At the upper edge,
translation by $N$ sites exhibits $N$ identical parabolic blocks
at $+1$ (Proposition~\ref{prop:center}). At the interior reference
energy used for gap collapse, translation by one site distinguishes
the lifts of the scalar Floquet phases. In this frame the reference
matrix is diagonal unitary, with simple, pairwise nonresonant
multipliers (Lemmas~\ref{cor:ell} and~\ref{lem:seed}). Thus the
upper edge carries a multiple parabolic multiplier, while at the
interior energy distinct elliptic modes share the same zero
Lyapunov exponent.

For a polynomial perturbation of actual degree $m\ge N$,
the slow reduction refines the quantitative-duality picture \cite{GJYZ26}.
On the sets with acceleration $k$,
\begin{equation*}\label{eq:five-splitting}
 \mathbb C^{2m}
 =
 E^{ss}\oplus E^{sc}\oplus E^c\oplus E^{cu}\oplus E^{uu},
 \qquad
 E^{\mathrm{slow}}=E^{sc}\oplus E^c\oplus E^{cu}.
\end{equation*}
Here $E^{ss}$ and $E^{uu}$ are the strongly contracting
and expanding bundles, while $E^{sc}$ and $E^{cu}$
describe weak contraction and expansion within the
slow dynamics. 
The respective dimensions are
$m-N,N-k,2k,N-k,m-N$.
Away from the transition values, the nonzero slow exponents
have size $\sqrt{\varepsilon}$, whereas the fast exponents remain
separated from zero. The fixed slow bundle therefore contains
the entire changing neutral center, together with the weak
hyperbolic directions created by the unfolding.

\subsubsection{Hermitian-symplectic monotonicity and discrete Hellmann--Feynman identity}
\label{sec:unfolding}

At fixed potential the slow cocycle is monotonic in energy after a
fixed iterate (Corollary~\ref{parmono}). To understand how the same
family can exhibit different accelerations or a collapsed gap, we
also need to know how a perturbation of the potential moves its
individual modes. Proposition~\ref{energyvar} gives both variations
in the same Hermitian-symplectic form.

Avila \cite{Avila2009LP} made explicit the monotonicity in energy
inherent in the Schr\"odinger structure. After two transfer steps,
every real projective direction turns strictly clockwise as energy
increases. Avila--Krikorian \cite{AK15}
developed  monotonic theory for general $\mathrm{SL}(2,\mathbb R)$
cocycles. The theory extends to symplectic cocycles \cite{Xu19}
and to Hermitian-symplectic cocycles \cite{WXZ}.

\begin{definition}[\cite{WXZ}]
Let $\mathcal J$ be a nondegenerate skew-Hermitian form, and let $\tau\mapsto A_\tau\in C(\mathbb T,\mathrm{HSp}(\mathcal J))$
be a $C^1$ family. The family $(\alpha,A_\tau)$ is \emph{positively
monotonic} in $\tau$ if, for every $x,\tau$ and every nonzero $z$
with $z^*\mathcal J z=0$,
\begin{equation*}\label{eq:slow-monotonicity}
 z^*A_\tau(x)^*\mathcal J\partial_\tau A_\tau(x)z>0 .
\end{equation*}
\end{definition}

This viewpoint reveals the dynamical content of the classical
Moser--P\"oschel argument \cite{MP84}. A nontrivial parabolic
normal form separates elliptic and hyperbolic dynamics, and an
energy perturbation toward the hyperbolic side opens a spectral
gap \cite{Puig2,Puig1}. Although the classical argument uses
the Schr\"odinger form, Avila--You--Zhou \cite{AYZ23} bring out
the monotonicity underlying this use of the structure through
the projective action and the associated fibred rotation number.
They make this principle effective in the almost-reducible
setting and combine it with a cone argument to resolve the
noncritical DTMP. Li--Wu \cite{LW} extend
this rotation-number framework to 
Hermitian-symplectic cocycles.

For the Aubry dual $\widehat H_{v,2\cos,\alpha,\theta}$ with polynomial $v$, the center bundle itself varies with energy implicitly.
Wang--Xu--Zhou \cite[Theorem~4.2]{WXZ} showed that a partially
hyperbolic center inherits the energy monotonicity locally in suitable
coordinates \footnote{This was first realized by A.Avila.}.  Li--Xu--Zhou \cite{LXZ26} obtained
a piecewise monotonic reduction over global compact energy intervals.
Together with a discrete Hellmann--Feynman identity \cite[Lemma~5.9]{LXZ26}, this extends
parabolic gap opening to reduced cocycles and yields the supercritical
Type~I DTMP.
Ge--Jitomirskaya--You also use this energy monotonicity in arbitrary finite
center dimension to  prove the supercritical Cantor-spectrum
\cite{GJY26+}.

Li--Xu--Zhou \cite[\S5]{LXZ26} express the energy monotonicity through the
motion of Lagrangian subspaces. This allows us to connect the spectral problem to
Krein's theory of positive symplectic paths and to the Maslov index
\cite{Ekeland90,Long02,FJN03}. This geometric viewpoint will be crucial to the mechanisms
developed below. Along a continuous path of Hermitian-symplectic matrices, a simple
multiplier on the unit circle can leave it only at a collision
involving both Krein signs
\cite{LM97}.

Another main input is the discrete Hellmann--Feynman identity in
Proposition~\ref{energyvar}, which extends the finite-range identity
of Li--Xu--Zhou \cite[Lemma~5.9]{LXZ26} to the slow solution bundle
of the infinite-range Aubry dual. It explains the variation relation  with respect to both energy and the potential. 
For both applications below we use essentially the phase average of this identity, it gives all the various  monotonicity we need, which we try to explain below.

\medskip
\noindent\emph{Parabolic collisions and coexistence of accelerations.}
At the upper spectral edge $E_0$, the $N$ residue classes carry
translates of the same positive edge state $p$  (Lemma~\ref{lem:edge}),
 obtained by combining
almost reducibility \cite{AJ10,Eliasson92} with
a variational characterization of the spectral edge.
Such a solution plays a central role
in reaction--diffusion theory, where principal
eigenvalues govern persistence and enter variational
formulas for propagation speeds
\cite{BNV94,BHN05,LWZT,NR17}.
On the common parabolic eigenspace, the averaged
Hellmann--Feynman identity (Corollary \ref{cor:modevar}) separates the roles of energy and
potential: energy shifts every mode equally, while the potential
acts through a Hermitian matrix $M_W$
(Lemma~\ref{lem:average}).
The positive edge state of Lemma~\ref{lem:edge} has positive
overlap with its translate by $\alpha$. This allows the first
Fourier harmonic of $W$ to distinguish any two Fourier modes.
Hence $M_W$ has simple spectrum for generic $W$
(Lemma~\ref{lem:generic}), splitting the common edge into
distinct leading thresholds.

In the scale $E=E_0-\varepsilon\xi$, fix $\varepsilon>0$.
Increasing $\xi$ lowers the energy and strictly increases
every ordered eigenvalue of the averaged Hermitian block
(Lemma~\ref{lem:average}).
Finite-order averaging retains this leading structure
(Proposition~\ref{prop:center}).
Together with the negative parabolic shear from
Lemma~\ref{lem:edge}, it determines the unfolding shown in
Figure~\ref{fig:positive-parabolic}: a contracting and an
expanding multiplier meet at $+1$ and enter the unit circle
with opposite Krein signs, the upper branch turning
counterclockwise and the lower branch clockwise.
Away from the thresholds, this directed phase motion persists
in the constant approximation
(Lemma~\ref{lem:mp-transversality}).
Each threshold crossed adds one elliptic pair; their separation
allows elliptic and hyperbolic modes to coexist.

The threshold separation also controls the quasiperiodic
perturbation. Each elliptic angle follows a square-root law from its own
threshold, so even branches rotating in the same direction
have distinct angular speeds. Their relative motion makes
the resonance conditions transverse to energy variation away
from the thresholds (Lemma~\ref{lem:mp-transversality}).
This provides the small-divisor control needed to obtain
reducibility of the elliptic dynamics on positive-measure
sets of energies. On these sets, $k$ elliptic pairs give
exactly $2k$ zero slow Lyapunov exponents and hence
acceleration $k$ (Lemma~\ref{lem:par}).
The quantitative estimates retain a spectral set of measure
at least $c_k\varepsilon$ for every $k=1,\ldots,N$, all for
the same perturbed potential (Theorem~\ref{thm:strata}).

The collision at each threshold is the familiar parabolic
unfolding of Hamiltonian stability theory
\cite{Galin1982,BHTB90,CLW19}. Here the potential separates the
collision energies, and the acceleration records how many pairs
have reached the elliptic side.

\begin{figure}[tbp]
\centering
\definecolor{pathblue}{RGB}{27,85,135}
\definecolor{pathochre}{RGB}{155,72,25}
\tikzset{
  slowaxis/.style={line width=.32pt,black!32},
  slowcircle/.style={line width=.56pt,black!62},
  slowarrow/.style={-{Stealth[length=1.7mm,width=1.25mm]},line width=.95pt},
  slowpoint/.style={circle,inner sep=1.65pt,fill},
  slowtext/.style={font=\scriptsize,align=center}
}
\begin{minipage}[t]{.41\linewidth}\centering
\begin{tikzpicture}[x=1.35cm,y=1.35cm]
  \path[use as bounding box] (-1.58,-1.63) rectangle (1.78,1.56);
  \node[slowtext] at (.08,1.40) {$E=E_0-\varepsilon\xi$, \ $\xi\uparrow$};
  \draw[slowaxis] (-1.22,0)--(1.75,0);
  \draw[slowaxis] (0,-1.18)--(0,1.18);
  \draw[slowcircle] (0,0) circle (1);
  \draw[slowarrow] (.54,0)--(.92,0);
  \draw[slowarrow] (1.65,0)--(1.08,0);
  \node[slowpoint] at (.60,0) {};
  \node[slowpoint] at (1.667,0) {};
  \node[circle,draw,fill=white,line width=.65pt,inner sep=1.75pt] at (1,0) {};
  \node[slowtext,anchor=north west] at (1.005,-.06) {$+1$};
  \foreach \sign in {1,-1}{
    \draw[pathblue!60,line width=.65pt] (0:1)
      arc[start angle=0,end angle={\sign*44},radius=1];
    \draw[slowarrow,pathblue] ({\sign*44}:1)
      arc[start angle={\sign*44},end angle={\sign*73},radius=1];
    \node[slowpoint,pathblue] at ({\sign*52}:1) {};
    \draw[slowarrow,pathochre] ({\sign*6}:1)
      arc[start angle={\sign*6},end angle={\sign*33},radius=1];
    \node[slowpoint,pathochre] at ({\sign*17}:1) {};
  }
  \node[slowtext] at (.05,-1.39) {successive thresholds};
\end{tikzpicture}
\caption{Parabolic collisions.}
\label{fig:positive-parabolic}
\end{minipage}\hfill
\begin{minipage}[t]{.57\linewidth}\centering
\begin{tikzpicture}[x=1cm,y=1cm]
  \path[use as bounding box] (-1.54,-2.20) rectangle (4.69,2.11);
  \node[slowtext] at (1.57,1.89)
    {$R_{\rho_1}\oplus R_{\rho_2}$, \ $0<\rho_1<\rho_2<\tfrac12$};
  \foreach \pos in {0,3.14}{
    \begin{scope}[xshift=\pos cm]
      \draw[slowaxis] (-1.20,0)--(1.20,0);
      \draw[slowaxis] (0,-1.20)--(0,1.20);
      \draw[slowcircle] (0,0) circle (1.07);
      \foreach \sign in {1,-1}{
        \node[slowpoint,pathblue] at ({\sign*42}:1.07) {};
        \node[slowpoint,pathochre] at ({\sign*122}:1.07) {};
      }
    \end{scope}
  }
  \foreach \sign in {1,-1}{
    \draw[slowarrow,pathblue] ({\sign*45}:1.07)
      arc[start angle={\sign*45},end angle={\sign*15},radius=1.07];
    \draw[slowarrow,pathochre] ({\sign*125}:1.07)
      arc[start angle={\sign*125},end angle={\sign*95},radius=1.07];
    \begin{scope}[xshift=3.14cm]
      \draw[slowarrow,pathblue] ({\sign*39}:1.07)
        arc[start angle={\sign*39},end angle={\sign*69},radius=1.07];
      \draw[slowarrow,pathochre] ({\sign*125}:1.07)
        arc[start angle={\sign*125},end angle={\sign*95},radius=1.07];
    \end{scope}
  }
  \node[slowtext] at (0,-1.62) {$E\uparrow$\\fixed potential};
  \node[slowtext] at (3.14,-1.62)
    {$\rho_1\uparrow,\ \rho_2\downarrow$\\fixed total rotation};
\end{tikzpicture}
\caption{Elliptic motion at a prescribed label.}
\label{fig:positive-elliptic}
\end{minipage}
\end{figure}

\medskip
\noindent\emph{Elliptic modes and collapsed gaps.}
At the interior reference energy, all $N$ slow pairs are elliptic.
A gap label fixes a signed sum of their reducible phases, whereas a
resonance involves the difference of two phases. With a single pair,
the signed sum is itself the phase difference, up to sign. The label
therefore forces a resonance, as in the supercritical Type~I
gap-opening mechanism \cite{LXZ26}. With several pairs, the same
constraint leaves the relative phases free.
Figure~\ref{fig:positive-elliptic} illustrates this for
$R_{\rho_1}\oplus R_{\rho_2}$, where one angle can increase while
the other decreases, keeping their sum fixed.

The perturbation makes this freedom effective. By the averaged
Hellmann--Feynman identity (Corollary~\ref{cor:modevar}), each mode
responds to the potential through its own solution profile.
Adjusting the energy compensates the change in total rotation.
Even after this adjustment, every pairwise phase difference has
nonzero first variation for a generic perturbation direction
(Lemma~\ref{cl:variation}).  Varying the amplitude at fixed IDS allows us to control the small
divisors by excluding amplitudes near the resonances
$$
  \theta_{\mathbf a}-\theta_{\mathbf b}=n\alpha\pmod 1,
  \qquad \mathbf a\ne\mathbf b,\quad n\in\mathbb Z,
$$
while such a resonance can open a spectral gap\footnote{ This was first realized by Q.~Zhou and S.~Jitomirskaya  in
2015.}.
For $n=0$, avoiding resonance means that the multipliers themselves
are distinct, as in Figure~\ref{fig:positive-elliptic}. This is why
we construct the elliptic slow cocycle from the one-site shift (Lemma \ref{cor:ell}).
At the reference energy, the $N$-step return used at the parabolic
edge would identify the $N$ distinct roots of each scalar Floquet
multiplier.

Lemma~\ref{cl:KAM} controls the remaining resonances while preserving
the label. At each step, the diagonal corrections shift the phases.
Their signed sum remains strictly monotone in energy, so a small
energy adjustment restores the prescribed rotation. Along this
constraint, each phase difference still varies with the perturbation
amplitude. One can therefore exclude the amplitudes at which a
difference comes too close to $n\alpha$, with a lower bound that
depends on $|n|$. For the remaining amplitudes, the limiting energy
$E_\infty$ has the prescribed IDS and the slow cocycle is reducible
to a diagonal unitary constant with quantitatively nonresonant phases.

Now fix the potential and vary only the energy. As the energy arrows
in Figure~\ref{fig:positive-elliptic} indicate, eigenphases of
opposite Krein sign rotate in opposite directions. Those same signs
enter the rotation number, so the angular speeds add in the IDS
derivative. By the averaged identity, each mode contributes the
positive mass of its solution profile. The IDS therefore crosses the
prescribed level with positive derivative at $E_\infty$
(Lemma~\ref{cl:collapse}). Its level set consists of a single energy,
so the labelled gap is collapsed while all $N$ neutral pairs remain
present.
The local stability in Figure~\ref{fig:positive-elliptic} is the
Krein stability of separated elliptic modes, also encountered in
the normal dynamics of Hamiltonian tori \cite{LM97,BHN07}.

\subsubsection{Effective dispersion}
\label{subsubsec:spectral-response}

The two geometries have different dispersion laws. Near the outer edge,
the modes inherit a common quadratic leading dependence on quasimomentum,
while the potential separates their threshold energies
(Theorem~\ref{thm:landscape}). In the leading normal form, real
quasimomentum describes oscillation and its imaginary continuation
describes weak exponential growth. These are
the two sides of the parabolic collision in
Figure~\ref{fig:positive-parabolic}. Quantitative duality turns the same
slow spectral data into the IDS and the real and complexified Lyapunov
profiles. The uniform slow reduction makes this description independent
of the hopping range. The resulting asymptotics hold on the whole energy
window, including gaps, even though exact Bloch waves are obtained on
the reducibility set.

At the collapsed energies, the first-order Bloch branches have nonzero
linear slopes (Corollary~\ref{cor:linear-dispersion}). These local laws
remain compatible with gaps accumulating at the reference energy. The energy and
potential variations in Lemma~\ref{cl:variation} have a direct meaning:
energy variation gives the inverse branch velocity, and the expectation of the
perturbation shifts the branch. The density of states adds the inverse
absolute velocities. Oppositely moving branches therefore reinforce the
spectral response, explaining the positive IDS derivative at a collapsed
label. Thus the distinction between the two mechanisms is visible in
dispersion itself: quadratic turning points permit changes of dynamical
type, while the nonzero first-order slopes give a strict IDS crossing at
a collapsed label.

\section{Preliminaries}
\subsection{Notations}
Set
$$
\mathbb{T}_{h}=\{x\in\mathbb{C}/\mathbb{Z}:|\operatorname{Im}x|<h\},\quad \mathbb{D}_{s}=\{z\in\mathbb{C}:|z|<s\}.
$$ 
For a complex Banach space $\mathscr Y$, let $\mathscr H_{h,s}(\mathscr Y)$ be the space of bounded holomorphic maps from $\mathbb T_h\times\mathbb D_s^2$ into $\mathscr Y$, with norm
$$
\|F\|_{h,s;\mathscr Y}=\sup_{\substack{x\in\mathbb T_h\\|\varepsilon|,|t|<s}}\|F(x,\varepsilon,t)\|_{\mathscr Y}.
$$
We write $F=O_{h,s}(s^p)$ if $\|F\|_{h,s;\mathscr Y}\leqslant Cs^p$. We omit $\mathscr Y$ when it is clear, and $h$ or $s$ for functions independent of $x$ or $(\varepsilon,t)$, respectively.
We write $\langle F\rangle=\int_{\mathbb T}F(x)\,dx$ and
$\langle F\rangle_{\mathrm{diag}}$ for the diagonal part of this average.

\subsection{Linear cocycles and Lyapunov exponents}
 Given $ \alpha \in \mathbb{R} \setminus \mathbb{Q} $, $A\in C(\mathbb{T},\mathrm{GL}(m,\mathbb{C}))$, we define the one-frequency cocycle $ (\alpha, A) $ as follows:
$$
(\alpha, A) \colon \left\{
\begin{array}{rcl}
\mathbb{T} \times \mathbb{C}^{m} & \to & \mathbb{T} \times \mathbb{C}^{m} \\
(x, v) & \mapsto & (x + \alpha, A(x) \cdot v)
\end{array}
\right. .
$$
The iterates of $ (\alpha, A) $ take the form $ (\alpha, A)^n = (n\alpha, A_n) $, where
$$
A^{(n)}(x) :=
\begin{cases}
A(x + (n-1)\alpha) \cdots A(x + \alpha) A(x), & n \geqslant 0 \\
A^{-1}(x + n\alpha) A^{-1}(x + (n+1)\alpha) \cdots A^{-1}(x - \alpha), & n < 0
\end{cases} .
$$
We denote the Lyapunov exponents of $ (\alpha, A) $ by $ L_1(A) \geqslant L_2(A) \geqslant \ldots \geqslant L_m(A) $, ordered according to their multiplicities. These exponents are defined as follows:
$$
L_k(A) = \lim_{n \to \infty} \frac{1}{n} \int_{\mathbb{T}} \ln \sigma_k(A^{(n)}(x)) \, dx,
$$
where $ \sigma_1(A^{(n)}(x)) \geqslant \ldots \geqslant \sigma_m(A^{(n)}(x)) $ denote the singular values.

Let $(\alpha, A)$ be a cocycle. An invariant continuous splitting
$
\mathbb{C}^{d} = E^1(x) \oplus \cdots \oplus E^k(x)$ 
is called \textit{dominated} if there exists an integer $n \geqslant 1$ such that for any $x$, any $1 \leqslant j < k$, and any unit vectors $w_j \in E^j(x)$ and $w_{j+1} \in E^{j+1}(x)$, we have:
$\|A^{(n)}(x) w_j \| > \| A^{(n)}(x) w_{j+1} \|.$
A cocycle $(\alpha, A)$ is said to be \textit{partially hyperbolic} if there exists an invariant dominated continuous splitting
$$
\mathbb{C}^{d} = E^u(x) \oplus E^c(x) \oplus E^s(x),
$$
where $E^u$, $E^c$, and $E^s$ denote the unstable, center, and stable bundles, respectively. Moreover, there exist constants $C > 0$ and $c > 0$ such that for every $n \geqslant 0$, the following contraction properties hold:
\begin{align*}
    \|A^{(n)}(x)v\| &\leqslant C e^{-cn}\|v\|, \quad \text{for } v \in E^s(x), \\
    \|A^{(n)}(x)^{-1}v\| &\leqslant C e^{-cn}\|v\|, \quad  \text{for } v \in E^u(x+n\alpha).
\end{align*}

\subsection{Hermitian symplecticity and fibered rotation number}\label{sec:rotation-prelim}
Let $\mathcal J$ be a nondegenerate skew-Hermitian matrix defining the form
$(u,v)\mapsto u^*\mathcal J v$ on $\mathbb C^{2N}$.
Its preserving Lie group and Lie algebra are
$$
\begin{aligned}
\mathrm{HSp}(\mathcal J)&=\{C\in\mathrm{GL}(2N,\mathbb C):C^{*}\mathcal J C=\mathcal J\},\quad
\mathfrak{hsp}(\mathcal J)&=\{A\in\mathfrak{gl}(2N,\mathbb C):A^{*}\mathcal J+\mathcal J A=0\}.
\end{aligned}
$$
Assume $i\mathcal J$ has signature $(N,N)$. Two canonical representatives are
$$
J_{2N}=\begin{pmatrix}0&-I_N\\I_N&0\end{pmatrix},\qquad
\widehat J_{2N}=\begin{pmatrix}-iI_N&0\\0&iI_N\end{pmatrix}.
$$

Let $A\in C(\mathbb T,\mathrm{HSp}(\mathcal J))$ be homotopic to a constant.
The Hermitian symplectic cocycle $(\alpha,A)$ admits a fibered rotation
number \cite{LW}. It is convenient to work with the form $\widehat J_{2N}$,
for which every Lagrangian subspace is the graph of a matrix $Z\in\mathrm U(N)$.
Writing $A(x)=(A_{ij}(x))_{i,j=1}^2\in\mathrm{HSp}(\widehat J_{2N})$,
its action is
$$
 A(x)\cdot Z=(A_{21}(x)+A_{22}(x)Z)(A_{11}(x)+A_{12}(x)Z)^{-1}.
$$
There is a continuous lift $\psi_A:\mathbb T\times\mathrm U(N)\to\mathbb R$
such that
$$
 e^{2\pi i\psi_A(x,Z)}=\frac{\det(A(x)\cdot Z)}{\det Z}.
$$
Let $\mu$ be a probability measure invariant under
$(x,Z)\mapsto(x+\alpha,A(x)\cdot Z)$ and projecting to Lebesgue measure
on $\mathbb T$. The fibered rotation number is
$$
 \rho(\alpha,A)=\int_{\mathbb T\times\mathrm U(N)}\psi_A(x,Z)\,d\mu(x,Z)
 \pmod{\mathbb Z}.
$$
It is independent of $\mu$ and the choice of lift.
For local comparisons, we choose compatible continuous real lifts.

\begin{lemma}\label{lem:unitary-rotation}
Let $D=\operatorname{diag}(D_-,D_+)$, where
$D_-,D_+\in\mathrm U(N)$ are constant. Then
\begin{equation}\label{eq:unitary-rotation}
 \rho(\alpha,D)
 =\frac{\arg\det D_+-\arg\det D_-}{2\pi}
 \pmod{\mathbb Z}.
\end{equation}
Moreover, there are $c_N,C_N>0$, such that every
$F\in C(\mathbb T,\mathfrak{hsp}(\widehat J_{2N}))$ with
$\|F\|_{C^0}<c_N$,
\begin{equation}\label{cl:compact-rotation-bound}
 |\rho(\alpha,De^F)-\rho(\alpha,D)|
 \le C_N\|F\|_{C^0}
\end{equation}
for compatible real lifts. The constants are uniform in $D$ and
the irrational frequency $\alpha$.
\end{lemma}

\begin{proof}
For $D$, the graph action is $Z\mapsto D_+ZD_-^{-1}$, so
$\psi_D=(\arg\det D_+-\arg\det D_-)/(2\pi)$ is constant.
Integration gives \eqref{eq:unitary-rotation}.

Write $e^{F(x)}=(P_{ij}(x))_{i,j=1}^2$. Since
$P_{21}+P_{22}Z=(P_{22}+P_{21}Z^{-1})Z$, we have
$$
 \frac{\det((De^{F(x)})\cdot Z)}{\det Z}=
 \frac{\det D_+}{\det D_-}
 \frac{\det(P_{22}+P_{21}Z^{-1})}{\det(P_{11}+P_{12}Z)}.
$$
Both matrices in the last quotient are within
$2\|e^F-I_{2N}\|_{C^0}$ of $I_N$, uniformly for $x\in\mathbb T$
and $Z\in\mathrm U(N)$. Using their logarithms near $I_N$ gives
$$
\begin{aligned}
 \psi_{De^F}(x,Z)-\psi_D
 &=\frac1{2\pi}\operatorname{Im}\operatorname{tr}
 \bigl(\log(P_{22}+P_{21}Z^{-1})-\log(P_{11}+P_{12}Z)\bigr).
\end{aligned}
$$
The estimates $\|e^F-I_{2N}\|_{C^0}\le2\|F\|_{C^0}$ and
$\|\log(I_N+Q)\|\le2\|Q\|$ for $\|Q\|\le1/2$ imply
$\|\psi_{De^F}-\psi_D\|_{C^0}\le C_N\|F\|_{C^0}$.
Integrating against an invariant probability measure for
$(\alpha,De^F)$ proves \eqref{cl:compact-rotation-bound}.
\end{proof}

\begin{remark}[\cite{LW}]\label{rem:rotation-properties}
The rotation number is continuous in the $C^0$ topology. Its definition gives $\rho(q\alpha,A_q)=q\rho(\alpha,A) \pmod \Z$ for every integer $q\ge1$ and
$$
 \rho(\alpha,A_1\oplus A_2)=\rho(\alpha,A_1)+\rho(\alpha,A_2) \pmod \Z,
$$
where the direct sum is compatible
with the form $\mathcal J_1\oplus\mathcal J_2$. If for all $x\in\T$, $A(x)$ preserve a fixed Lagrangian
subspace, then $\rho(\alpha,A)=0\pmod{\mathbb Z}$.
\end{remark}

For $B\in C(\mathbb T,\mathrm{HSp}(\mathcal J))$, its degree $\deg B$ is the winding of
$\det(B(x)\cdot I_N)$,
with the graph action computed in the $\widehat J_{2N}$ coordinates.

\begin{lemma}[{\cite{LW}}]\label{lem:rotation-conjugacy}
Let $A,B\in C(\mathbb T,\mathrm{HSp}(\mathcal J))$, with $A$ homotopic to a constant. Then
\begin{equation*}\label{eq:rotation-conjugacy}
 \rho\bigl(\alpha,B(\cdot+\alpha)^{-1}AB\bigr)
 =\rho(\alpha,A)-\alpha\deg B\pmod{\mathbb Z}.
\end{equation*}
\end{lemma}

\begin{remark}[\cite{LW}]\label{rem:rotation-conjugacy}
If $B$ is $2$-periodic and $B(\cdot+\alpha)^{-1}AB$ is $1$-periodic,
the degree term is $-\frac{\alpha}{2}\deg[B(2\cdot)]$.
For continuous one-periodic families $A_t,B_t$ on a connected parameter
set, $\deg B_t$ is constant. Thus, writing
$\widetilde A_t=B_t(\cdot+\alpha)^{-1}A_tB_t$, continuous real lifts satisfy
\begin{equation}\label{eq:rotation-conjugacy-difference}
 \rho(\alpha,\widetilde A_t)-\rho(\alpha,\widetilde A_s)
 =\rho(\alpha,A_t)-\rho(\alpha,A_s).
\end{equation}
\end{remark}
\subsection{Schr\"odinger operators, IDS and Aubry duality}\label{sec:aubry-duality}
For real analytic $v$, recall the operator $H_{v,\alpha,x}$ in \eqref{eq:operator}.  Its integrated density of states is
\begin{equation*}\label{eq:idsdef}
\mathcal{N}_{v,\alpha}(E)=\int_{\mathbb{T}}\langle\delta_{0},\mathbf{1}_{(-\infty,E]}(H_{v,\alpha,x})\delta_{0}\rangle\,\mathrm{d}x.
\end{equation*}

Writing $v(x)=\sum_{j\in\mathbb{Z}}\widehat v_{j}e^{2\pi ijx}$, its Aubry dual is the operator:
\begin{equation*}\label{eq:Lx_def}
(\widehat H_{v,2\cos,\alpha,x}u)_{n}=\sum_{j\in\mathbb{Z}}\widehat v_{j}u_{n+j}+2\cos\big(2\pi(x+n\alpha)\big)u_{n}.
\end{equation*}
The original and dual families have the same spectrum and IDS \cite{HP13}. For the family $v=v_{0}+\varepsilon W$, we use the abbreviations $\Sigma_{\varepsilon,W}$, $\mathcal{N}_{\varepsilon,W}$ and $\widehat H_{\varepsilon,x}$, as in the introduction.

Suppose now that $v$ is a real trigonometric polynomial of degree $m\geqslant1$. Write $v(x)=\sum_{j=-m}^{m}\widehat v_{j}e^{2\pi ijx}$, with $\widehat v_{m}\neq0$.
The equation $\widehat H_{v,2\cos,\alpha,x}u=Eu$ gives the one-step companion cocycle
$$
\mathcal A_E^{v}(x)=\frac{1}{\widehat v_m}
\begin{pmatrix}
 -\widehat v_{m-1}&\cdots&-\widehat v_1&
 E-2\cos(2\pi x)-\widehat v_0&
 -\widehat v_{-1}&\cdots&-\widehat v_{-m+1}&-\widehat v_{-m}\\
 \widehat v_m&&&&&&&\\
 &&&&&&&\\
 &&&&&&&\\
 &&&\ddots&&&&\\
 &&&&&&&\\
 &&&&&&\widehat v_m&
\end{pmatrix}.
$$ The associated Hermitian-symplectic form is $z^*J_vw$, where
$$
C_v=\begin{pmatrix}
 \widehat v_m&\cdots&\widehat v_1\\
 0&\ddots&\vdots\\
 0&0&\widehat v_m
\end{pmatrix},
\qquad
J_v=\begin{pmatrix}0&-C_v^*\\ C_v&0\end{pmatrix}.
$$
The matrix $J_v$ is nondegenerate and skew-Hermitian, and $iJ_v$
is Hermitian with signature $(m,m)$. For real $E$ and $x$,
$\mathcal A_E^v(x)^*J_v\mathcal A_E^v(x)=J_v$.

Choose a constant $T_{v}$ satisfying $T_{v}^{*}J_{v}T_{v}=J_{2m}$ and set $\widehat{\mathcal A}_E^v=T_{v}^{-1}\mathcal{A}_{E}^{v}T_{v}$, so that $\widehat{\mathcal A}_E^v(x)\in\mathrm{HSp}(J_{2m})$ for real variables. For example, one may take $T_{v}=\operatorname{diag}(C_{v}^{-1},I_{m})$. The rotation number is independent of the matrix representation of the symplectic form; in particular, $\rho(\alpha,\mathcal A_E^v)=\rho(\alpha,\widehat{\mathcal A}_E^v)$.

\begin{lemma}[\cite{LW}]\label{lem:idsrot}
There exists a continuous real lift of $E\mapsto\rho(\alpha,\mathcal A_E^v)$ such that, for every $E\in\mathbb R$,
\begin{equation*}\label{idsrot}
1-\mathcal N_{v,\alpha}(E)=\rho(\alpha,\mathcal A_E^v).
\end{equation*}
\end{lemma}

\begin{proof}
By \cite[Theorem~1.7 and Appendix~B]{LW} and the iteration property,
$$
m\big(1-\mathcal N_{v,\alpha}(E)-\rho(\alpha,\mathcal A_E^v)\big)\in\mathbb Z.
$$
In the upper unbounded spectral gap, the cocycle is uniformly hyperbolic and its stable Lagrangian bundle is homotopic to a constant. Choose the continuous real lift to be zero there, where $\mathcal N_{v,\alpha}=1$. The displayed expression is continuous in $E$ and therefore vanishes identically.
\end{proof}

\section{The solution bundle and intrinsic cocycle}\label{sec:solution-bundle}
Let $v_0(x)=2\lambda\cos(2\pi Nx)$ with $\lambda>1$, and write
$\Sigma_0=\Sigma_{v_0,\alpha}$. Fix $W\in C_r^\omega(\mathbb T,\mathbb R)$
with Fourier coefficients $w_j$. For $E_*\in\Sigma_0$ and
$l\in[N,\infty]:=\{N,N+1,\ldots\}\cup\{\infty\}$, the dual equation for
$v_0+\varepsilon W^{(l)}$ at energy $E_*+t$ is
$\mathbf L_{\varepsilon,t}^{(l)}(x)u=0$, where
$$
 \bigl(\mathbf L_{\varepsilon,t}^{(l)}(x)u\bigr)_m
 =\lambda(u_{m+N}+u_{m-N})
  +\bigl(2\cos(2\pi(x+m\alpha))-E_*-t\bigr)u_m
  +\varepsilon\sum_{|j|\leqslant l}w_j u_{m+j}.
$$
Here $W^{(l)}$ is the $l$-th Fourier truncation of $W$, with
$W^{(\infty)}=W$.

At $(\varepsilon,t)=(0,0)$, the equation decouples on the residue
classes modulo $N$. Thus the $N$-th iterate of the companion cocycle
$(\alpha,A_{*,1})$, where $A_{*,1}=\mathcal A_{E_*}^{v_0}$, is the
transfer cocycle of a strip with $N$ independent channels. 
In residue-class coordinates,
\begin{equation}\label{equ:direct}
(N\alpha,A_{*,N})=\bigoplus_{j=0}^{N-1}\big(N\alpha,S_{\lambda^{-1},E_{*}/\lambda}(\,\cdot+j\alpha)\big),
\end{equation}
where
$
S_{\lambda,E}(x)=\begin{pmatrix}E-2\lambda\cos(2\pi x)&-1\\1&0\end{pmatrix}
$
is the almost Mathieu cocycle.

For every $E_*\in\Sigma_0$, each summand has zero Lyapunov exponent
throughout $|\Im x|<(\log\lambda)/(2\pi)$
\cite{Avila15}. Hence, for $0<h<(\log\lambda)/(2\pi)$ and
$\eta>0$, the uniform subadditive bound gives
\begin{equation}\label{eq:reference-growth}
 \sup_{E_*\in\Sigma_0}\ \sup_{|\Im x|\leqslant h}
 \|A_{*,n}(x)\|\leqslant C_{h,\eta}e^{\eta|n|},
 \qquad n\in\mathbb Z.
\end{equation}The bound is uniform in $E_{*}$ and $|\operatorname{Im}x|\leqslant h$
by a finite-cover argument using continuity of the finite iterates
and compactness of $\Sigma_{0}\times[-h,h]$.

Fix $0<a<2\pi r$ and set
$$
 X_a=\Bigl\{u:\ \|u\|_a=\sum_{m\in\mathbb Z}|u_m|e^{-a|m|}<\infty\Bigr\},
 \qquad
 \mathcal M_{\varepsilon,t}^{(l)}(x)
 =\ker_{X_a}\mathbf L_{\varepsilon,t}^{(l)}(x).
$$
For finite $l$, let
$\mathcal E_lu=(u_{l-1},\ldots,u_0,u_{-1},\ldots,u_{-l})^{\mathsf T}$
and write $\mathcal E=\mathcal E_N$. At the unperturbed equation,
$\mathcal E u$ specifies two successive values on each channel and
hence determines a unique solution. By \eqref{eq:reference-growth},
every such solution belongs to $X_a$. Denote the inverse of
$\mathcal E$ on $\mathcal M_{0,0}^{(l)}(x)$ by $\widetilde U_*(x)$.
With $(Su)_m=u_{m+1}$, the $q$-step reference matrix is
$A_{*,q}(x)=\mathcal E S^q\widetilde U_*(x)$.

Proposition~\ref{prop:slow1} shows that the same $2N$ coordinates
parametrize the weighted solutions after perturbation, uniformly in
the hopping range. Translation on this bundle defines the slow
cocycle. 

\begin{proposition}\label{prop:slow1}
Let $\alpha\in\mathbb{R}\backslash\mathbb{Q}$ and fix an integer $q\geqslant1$. There exist $s_{*},h,\mu,C>0$, independent of $E_{*}\in\Sigma_{0}$ and $l\in\{N,N+1,\ldots\}\cup\{\infty\}$, such that
for $|\Im x|<h$ and $|\varepsilon|,|t|<s_*$ the map
$
 \mathcal E:\mathcal M_{\varepsilon,t}^{(l)}(x)
       \longrightarrow\mathbb C^{2N}$
is an isomorphism. Its inverse $\widetilde U_{\varepsilon,t}^{(l)}(x)$
is jointly analytic in $(x,\varepsilon,t)$. In this frame the slow
cocycle is $(q\alpha,A_{\varepsilon,t}^{(l)})$, where
$
 S^q\widetilde U_{\varepsilon,t}^{(l)}(x)
   =\widetilde U_{\varepsilon,t}^{(l)}(x+q\alpha)
      A_{\varepsilon,t}^{(l)}(x).
$
For every $0<s\leqslant s_*$,
\begin{align*}
 \|\widetilde U_{\varepsilon,t}^{(l)}-\widetilde U_*\|_{h,s}
 +\|A_{\varepsilon,t}^{(l)}-A_{*,q}\|_{h,s}
 &\leqslant Cs,\\
 \|\widetilde U_{\varepsilon,t}^{(l)}
       -\widetilde U_{\varepsilon,t}^{(\infty)}\|_{h,s}
 +\|A_{\varepsilon,t}^{(l)}
       -A_{\varepsilon,t}^{(\infty)}\|_{h,s}
 &\leqslant Cs e^{-\mu l}\qquad(l<\infty).
\end{align*}
\end{proposition}
\begin{remark}
 The slow cocycle $A_{\varepsilon,t}^{(l)}$ exists for every
irrational frequency, while our main applications require
$\alpha$ to be Diophantine.
\end{remark}

\begin{proof}
We divide the proof into three steps.

\noindent\textbf{Step 1. Inverting the operator-valued symbol.}
Fix $0<\chi<a/2$. Applying \eqref{eq:reference-growth} with
$\eta=\chi$ gives
$$
\sup_{E_{*}\in\Sigma_{0},\,x\in\mathbb{T}_{h}}
\|A_{*,n}(x)\|\leqslant Ce^{\chi|n|},\qquad n\in\mathbb{Z}.
$$
Since $A_{*,1}$ propagates the initial coordinates $\mathcal{E}u$,
every unperturbed solution satisfies
$$
|u_{m}|\leqslant Ce^{\chi|m|}\|\mathcal{E}u\|,
\qquad m\in\mathbb{Z}.
$$
Consequently,
$$
\sup_{E_{*}\in\Sigma_{0},\,x\in\mathbb{T}_{h}}
\|\widetilde U_{*}(x)\|_{\mathcal{L}(\mathbb{C}^{2N},X_{a})}
\leqslant C\sum_{m\in\mathbb{Z}}e^{-(a-\chi)|m|}<\infty.
$$
Each coordinate of $\widetilde U_{*}$ is holomorphic by the
recurrence. The preceding summable bound gives locally uniform
convergence in $X_{a}$, so $\widetilde U_{*}$ is holomorphic with
values in $\mathcal{L}(\mathbb{C}^{2N},X_{a})$.

For a function $f$, write $(T_{\gamma}f)(x)=f(x+\gamma)$ and let $M_{f}$ denote multiplication by $f$. We use either $\mathscr{A}=B(\mathbb{T},\mathbb{C})$ or $\mathscr{A}=C_{h}^{\omega}(\mathbb{T},\mathbb{C})$, with the supremum norm. For $R,\sigma>0$, let
$$
\mathcal{A}_{R,\sigma}=\{z\in\mathbb{C}:Re^{-\sigma}<|z|<Re^{\sigma}\}.
$$
Let $\mathscr{W}_{R,\sigma,s}(\mathscr{A})$ be the algebra of series
$$
\begin{aligned}
K(z,\varepsilon,t)&=\sum_{k\in\mathbb{Z}}z^{k}M_{a_{k}(\cdot,\varepsilon,t)}T_{k\alpha},\\
\|K\|_{\mathscr{W}_{R,\sigma,s}(\mathscr{A})}&=\sum_{k\in\mathbb{Z}}R^{k}e^{\sigma|k|}\sup_{|\varepsilon|,|t|<s}\|a_{k}(\cdot,\varepsilon,t)\|_{\mathscr{A}}<\infty,
\end{aligned}
$$
where $a_{k}(\cdot,\varepsilon,t)\in \mathscr A$ are holomorphic in $(\varepsilon,t)$ and independent of $z$. The series converges normally on $\mathcal{A}_{R,\sigma}\times\mathbb{D}_{s}^{2}$. The product rule
$$
(M_{a_{k}}T_{k\alpha})(M_{b_{j}}T_{j\alpha})=M_{a_{k}b_{j}(\cdot+k\alpha)}T_{(k+j)\alpha}
$$
and the submultiplicativity of $R^{k}e^{\sigma|k|}$ make this a unital Banach algebra. We use the corresponding matrix algebra with the same notation.

Define
\begin{equation*}
D_{\varepsilon,t}^{(l)}(z)=\lambda z^{N}T_{N\alpha}+\lambda z^{-N}T_{-N\alpha}+M_{2\cos(2\pi x)-E_{*}}-tI+\varepsilon\sum_{|j|\leqslant l}w_jz^{j}T_{j\alpha},
\end{equation*}
and write $D_{0,0}=D_{0,0}^{(l)}$, which is independent of $l$.

\begin{lemma}\label{lem:inverseD}
	Fix $\chi<\delta<a<b<\varrho<2\pi r$, $0<\sigma<\min\{\delta-\chi,\varrho-b\}$, and $0<\mu<2\pi r-\varrho$. For sufficiently small $s_{*}$, $D_{\varepsilon,t}^{(l)}$ is invertible in $\mathscr{W}_{\rho,\sigma,s}$ for $\rho\in[e^{-b},e^{-\delta}]\cup[e^{\delta},e^{b}]$. Uniformly in $E_{*},l,s,\rho$,
\begin{enumerate}
\item $\|(D_{\varepsilon,t}^{(l)})^{-1}\|_{\mathscr{W}_{\rho,\sigma,s}(\mathscr{A})}\leqslant C,$
\item $\|(D_{\varepsilon,t}^{(l)})^{-1}-D_{0,0}^{-1}\|_{\mathscr{W}_{\rho,\sigma,s}(\mathscr{A})}\leqslant Cs,$
\item $\|(D_{\varepsilon,t}^{(l)})^{-1}-(D_{\varepsilon,t}^{(\infty)})^{-1}\|_{\mathscr{W}_{\rho,\sigma,s}(\mathscr{A})}\leqslant Cs e^{-\mu l}.$
\end{enumerate}
\end{lemma}

\begin{proof}
Let $\mathcal{S}(z)=z^{N}T_{N\alpha}$, $a_{*}(x)=(E_{*}-2\cos(2\pi x))/\lambda$, and
$$
\mathcal{A}(z)=\mathcal{S}(z)I_{2}-M_{S_{\lambda^{-1},E_{*}/\lambda}}
=\begin{pmatrix}\mathcal{S}(z)-M_{a_{*}}&I\\-I&\mathcal{S}(z)\end{pmatrix}.
$$
The product rule gives, for $j\geqslant0$,
$$
\begin{aligned}
(\mathcal{S}(z)^{-1}M_{S_{\lambda^{-1},E_{*}/\lambda}})^{j}
&=z^{-Nj}M_{(S_{\lambda^{-1},E_{*}/\lambda})^{(j)}(\cdot-jN\alpha)}T_{-jN\alpha},\\
(M_{(S_{\lambda^{-1},E_{*}/\lambda})^{-1}}\mathcal{S}(z))^{j}
&=z^{Nj}M_{((S_{\lambda^{-1},E_{*}/\lambda})^{(j)})^{-1}}T_{jN\alpha}.
\end{aligned}
$$

For $\rho\geqslant e^{\delta}$ the first power has norm at most $C\rho^{-Nj}e^{N(\sigma+\chi)j}$. For $\rho\leqslant e^{-\delta}$ the second has norm at most $C\rho^{Nj}e^{N(\sigma+\chi)j}$. Both bounds are at most $C e^{-N(\delta-\sigma-\chi)j}$ on the respective annuli. Factoring $\mathcal{A}(z)$ on the left by $\mathcal{S}(z)I_{2}$ on the outer annulus and by $-M_{S_{\lambda^{-1},E_{*}/\lambda}}$ on the inner annulus therefore gives
$$
\mathcal{A}(z)^{-1}=
\begin{cases}
\displaystyle\sum_{j\geqslant0}(\mathcal{S}(z)^{-1}M_{S_{\lambda^{-1},E_{*}/\lambda}})^{j}\mathcal{S}(z)^{-1},&\rho\geqslant e^{\delta},\\
\displaystyle-\sum_{j\geqslant0}(M_{(S_{\lambda^{-1},E_{*}/\lambda})^{-1}}\mathcal{S}(z))^{j}M_{(S_{\lambda^{-1},E_{*}/\lambda})^{-1}},&\rho\leqslant e^{-\delta}.
\end{cases}
$$
	The triangular identity
	$$
	\mathcal{A}(z)\begin{pmatrix}I&0\\\mathcal{S}(z)^{-1}&I\end{pmatrix}
	=\begin{pmatrix}\lambda^{-1}D_{0,0}(z)&I\\0&\mathcal{S}(z)\end{pmatrix}
	$$
	gives $D_{0,0}(z)^{-1}=\lambda^{-1}(1,0)\mathcal{A}(z)^{-1}\binom{1}{0}$, hence the uniform inverse bound at $(0,0)$. Since $\rho^{j}e^{\sigma|j|}\leqslant e^{\varrho|j|}$, Fourier decay gives
	$$
	\|D_{\varepsilon,t}^{(l)}-D_{0,0}\|_{\mathscr{W}_{\rho,\sigma,s}}\leqslant Cs,\qquad
	\|D_{\varepsilon,t}^{(l)}-D_{\varepsilon,t}^{(\infty)}\|_{\mathscr{W}_{\rho,\sigma,s}}\leqslant Cs e^{-\mu l}.
	$$
	The claimed estimates follow from the Neumann series and inverse difference identity. 
\end{proof}

\medskip
\noindent\textbf{Step 2. The transfer operator and the full fiber.}
Choose
$$
\chi<\delta<a<b<\varrho<2\pi r,\qquad 0<\mu<2\pi r-\varrho.
$$
For $\gamma\in\{a,\varrho\}$, let $\mathscr{B}_{\gamma}=B(\mathbb{T},X_{\gamma})$ with the supremum norm, and define
$$
\widehat{\mathcal{M}}_{\varepsilon,t,a}^{(l)}
=\{\Psi\in\mathscr{B}_{a}:\mathbf{L}_{\varepsilon,t}^{(l)}(x)\Psi(x)=0\text{ for every }x\in\mathbb{T}\}.
$$
The transfer operator is $(\mathscr{T}\Psi)(x)=S\Psi(x-\alpha)$. For $e^{-\varrho}<|z|<e^{\varrho}$, let
$$
(P_{z}f)_{m}(x)=z^{m}f(x+m\alpha)
$$
and
\begin{equation*}
(Q_{z}\Psi)_{m}(x)=
\begin{cases}
-\displaystyle\sum_{j=0}^{m-1}z^{m-1-j}\Psi_{j}(x+(m-j)\alpha),&m\geqslant1,\\
0,&m=0,\\
\displaystyle\sum_{j=1}^{-m}z^{m-1+j}\Psi_{-j}(x+(m+j)\alpha),&m\leqslant-1.
\end{cases}
\end{equation*}
For $k\neq0$, write
\begin{equation}\label{equ:part}
t_k=\lambda\mathbf{1}_{\{|k|=N\}}+\varepsilon w_k.
\end{equation}
Define\begin{equation}\label{eq:GzD}
\mathscr{G}_{z}^{(l)}\Psi(x)
=\sum_{k=1}^l t_k\sum_{j=0}^{k-1}z^{k-j-1}\Psi_{j}(x+(k-j)\alpha)-\sum_{k=1}^l t_{-k}\sum_{j=1}^{k}z^{-k-1+j}\Psi_{-j}(x-(k-j)\alpha).
\end{equation}

 Lemma~\ref{lem:inverseD} defines the extended resolvent
\begin{equation*}
\mathscr{R}_{\varepsilon,t}^{(l)}(z)
=P_{z}(D_{\varepsilon,t}^{(l)}(z))^{-1}\mathscr{G}_{z}^{(l)}+Q_{z}
\end{equation*}
on the two closed annuli $e^{-b}\leqslant|z|\leqslant e^{-\delta}$ and $e^{\delta}\leqslant|z|\leqslant e^{b}$. It is a bounded operator on $\mathscr{B}_{\varrho}$, holomorphic in $z,\varepsilon,t$. For $\kappa\in\{\delta,b\}$, let $\Gamma_{\kappa,\pm}$ be the positively oriented circles $|z|=e^{\pm\kappa}$ and set $\Gamma_{\kappa}=\Gamma_{\kappa,+}-\Gamma_{\kappa,-}$. Define
\begin{equation*}
\mathbf{P}_{\varepsilon,t}^{(l),\kappa}
=\frac{1}{2\pi i}\int_{\Gamma_{\kappa}}\mathscr{R}_{\varepsilon,t}^{(l)}(z)\,\mathrm{d}z.
\end{equation*}
Cauchy's theorem gives $\mathbf{P}_{\varepsilon,t}^{(l),\delta}=\mathbf{P}_{\varepsilon,t}^{(l),b}$; denote this operator by $\mathbf{P}_{\varepsilon,t}^{(l)}$.

\begin{lemma}\label{lem:fiberD}
For sufficiently small $s_{*}>0$, $\mathbf{P}_{\varepsilon,t}^{(l)}$ is a bounded projection and acts fiberwise:
\begin{equation*}
(\mathbf{P}_{\varepsilon,t}^{(l)}\Psi)(x)=P_{\varepsilon,t}^{(l)}(x)\Psi(x),\qquad P_{\varepsilon,t}^{(l)}(x)\in\mathcal{L}(X_{\varrho},X_{a}).
\end{equation*}
Its range is the full weighted solution space:
\begin{equation*}\label{eq:rankMcD}
\operatorname{Ran}P_{\varepsilon,t}^{(l)}(x)=\mathcal{M}_{\varepsilon,t}^{(l)}(x)=\ker_{X_{a}}\mathbf{L}_{\varepsilon,t}^{(l)}(x),\qquad \dim\mathcal{M}_{\varepsilon,t}^{(l)}(x)=2N.
\end{equation*}
The fiber projections are jointly analytic and, with $P_{*}=P_{0,0}^{(l)}$, satisfy
\begin{equation}\label{eq:projbdD}
\|P_{\varepsilon,t}^{(l)}-P_{*}\|_{h,s;\mathcal{L}(X_{\varrho},X_{a})}\leqslant Cs,\qquad
\|P_{\varepsilon,t}^{(l)}-P_{\varepsilon,t}^{(\infty)}\|_{h,s;\mathcal{L}(X_{\varrho},X_{a})}\leqslant Cs e^{-\mu l},
\end{equation}
uniformly in $0<s\leqslant s_{*}$ and $l\in[N,\infty]$.
\end{lemma}

\begin{proof}
We use the extended-resolvent construction of \cite{CCB26}, whose argument depends only on symbol invertibility, and thus applies also to the case considered here: $(z-\mathscr T)\mathscr R_{\varepsilon,t}^{(l)}(z)=I$ on $\mathscr B_{\varrho}$; moreover, on solution sections $\mathscr R_{\varepsilon,t}^{(l)}(z)$ is the genuine resolvent of $z-\mathscr T$.

Since $Q_{z}$ is holomorphic throughout $e^{-\varrho}<|z|<e^{\varrho}$, its integral over $\Gamma_{\delta}$ vanishes. Therefore
\begin{equation}\label{eq:projcoordD}
(\mathbf{P}_{\varepsilon,t}^{(l)}\Psi)_{m}(x)
=\frac{1}{2\pi i}\int_{\Gamma_{\delta}}z^{m}\big[(D_{\varepsilon,t}^{(l)}(z))^{-1}\mathscr{G}_{z}^{(l)}\Psi\big](x+m\alpha)\,\mathrm{d}z.
\end{equation}
The inverse symbol has the Laurent expansions
$$
(D_{\varepsilon,t}^{(l)}(z))^{-1}
=\sum_{k\in\mathbb{Z}}z^{k}M_{d_{k,\pm}^{(l)}(\cdot,\varepsilon,t)}T_{k\alpha}
$$
on the inner and outer annuli. Write each term of $\mathscr{G}_{z}^{(l)}\Psi$ as $g_{j,v}^{(l)}(\varepsilon)z^{v-1}\Psi_{j}(x+v\alpha)$. The contour integral keeps only $m+k+v=0$, giving
\begin{equation*}
(\mathbf{P}_{\varepsilon,t}^{(l)}\Psi)_{m}(x)
=\sum_{j,v\in\mathbb{Z}}g_{j,v}^{(l)}(\varepsilon)
\big[d_{-m-v,+}^{(l)}(x+m\alpha,\varepsilon,t)-d_{-m-v,-}^{(l)}(x+m\alpha,\varepsilon,t)\big]\Psi_{j}(x).
\end{equation*}
Absolute convergence follows from Lemma~\ref{lem:inverseD} and the Fourier decay of $W$. This proves that the operator acts fiberwise, with analytic fiber coefficients.
By \cite[Propositions~5.1 and~5.3]{CCB26}, the extended
resolvent restricts to the genuine resolvent of $z-\mathscr T$
on the solution section space, and
$\mathbf P_{\varepsilon,t}^{(l)}$ is a bounded projection
on $\mathscr B_{\varrho}$ whose range consists of solution
sections. Moreover, the contour estimates of \eqref{eq:projcoordD} give
$$
\sup_{x\in\mathbb T}
|(\mathbf P_{\varepsilon,t}^{(l)}\Psi)_m(x)|
\leqslant Ce^{\delta|m|}\|\Psi\|_{\mathscr B_{\varrho}},
\qquad m\in\mathbb Z.
$$
Since $\delta<a$, we obtain
$
\operatorname{Ran}\mathbf P_{\varepsilon,t}^{(l)}
\subseteq\widehat{\mathcal M}_{\varepsilon,t,a}^{(l)},
\|\mathbf P_{\varepsilon,t}^{(l)}\|_
{\mathcal L(\mathscr B_{\varrho},\mathscr B_a)}
\leqslant C.
$

Conversely, let
$\Phi\in\widehat{\mathcal M}_{\varepsilon,t,a}^{(l)}$.
Since $a<b<\varrho$, we have
$
\|\mathscr T^n\Phi\|_{\mathscr B_{\varrho}}
\leqslant e^{a|n|}\|\Phi\|_{\mathscr B_a}, n\in\mathbb Z.
$
Thus the genuine resolvent on the solution section space
admits the expansions
$$
\mathscr R_{\varepsilon,t}^{(l)}(z)\Phi
=
\begin{cases}
\displaystyle\sum_{n\geqslant0}z^{-n-1}\mathscr T^n\Phi,
&|z|=e^b,\\
\displaystyle-\sum_{n\geqslant0}z^n\mathscr T^{-n-1}\Phi,
&|z|=e^{-b}.
\end{cases}
$$
Integrating over $\Gamma_b$ gives
$\mathbf P_{\varepsilon,t}^{(l),b}\Phi=\Phi.$
Since $\mathbf P_{\varepsilon,t}^{(l),b}
=\mathbf P_{\varepsilon,t}^{(l)}$, we conclude that
$
\operatorname{Ran}\mathbf P_{\varepsilon,t}^{(l)}
=\widehat{\mathcal M}_{\varepsilon,t,a}^{(l)}.
$
By the fiberwise property, applying this identity to
sections supported at an arbitrary phase $x$ yields
$\operatorname{Ran}P_{\varepsilon,t}^{(l)}(x)=\ker_{X_a}\mathbf L_{\varepsilon,t}^{(l)}(x).$

At $(\varepsilon,t)=(0,0)$, \eqref{eq:reference-growth} gives
$|u_m|\leqslant C e^{\chi|m|}\|\mathcal E u\|$ with $\chi<a$,
so the solutions belong to $X_a$.
The initial-data parametrization therefore gives
$\operatorname{rank}P_*(x)=\dim\ker_{X_a}\mathbf L_{0,0}(x)=2N$.

The inverse-symbol bounds and the Fourier tails in \eqref{eq:GzD}, as in \cite[Proposition~5.3]{CCB26}, give
$$
\|\mathscr R_{\varepsilon,t}^{(l)}(z)-\mathscr R_{0,0}(z)\|_{s;\mathcal L(\mathscr B_{\varrho})}\leqslant Cs,\qquad
\|\mathscr R_{\varepsilon,t}^{(l)}(z)-\mathscr R_{\varepsilon,t}^{(\infty)}(z)\|_{s;\mathcal L(\mathscr B_{\varrho})}\leqslant Cs e^{-\mu l}.
$$

Using \eqref{eq:projcoordD} on constant sections $\Psi(x)=u\in X_{\varrho}$, the analytic coefficient estimates give
$$
\begin{aligned}
|[(P_{\varepsilon,t}^{(l)}(x)-P_{*}(x))u]_{m}|&\leqslant Cs e^{\delta|m|}\|u\|_{\varrho},\\
|[(P_{\varepsilon,t}^{(l)}(x)-P_{\varepsilon,t}^{(\infty)}(x))u]_{m}|&\leqslant Cs e^{-\mu l}e^{\delta|m|}\|u\|_{\varrho},
\end{aligned}
$$
uniformly on $\mathbb{T}_{h}\times\mathbb{D}_{s}^{2}$. Multiplication by $e^{-a|m|}$ and summation over $m\in\mathbb Z$ prove \eqref{eq:projbdD}, because $\delta<a$. 

In particular, after decreasing $s_*$, $\|P_{\varepsilon,t}^{(l)}(x)-P_*(x)\|_{\mathcal L(X_{\varrho})}<1$.
The projection identities extend holomorphically in $x$, so rank
stability gives $\operatorname{rank}P_{\varepsilon,t}^{(l)}(x)=2N$
throughout $\mathbb T_h\times\mathbb D_{s_*}^2$. This completes the proof.
\end{proof}

\medskip
\noindent\textbf{Step 3. Normalizing the initial data.}
We now choose a holomorphic frame whose initial-data matrix is $I_{2N}$.
Recall the initial-data map $\mathcal E:X_a\to\mathbb C^{2N}$,
\begin{equation}\label{eq:initD}
\mathcal Eu=(u_{N-1},\ldots,u_0,u_{-1},\ldots,u_{-N})^{\mathsf T}.
\end{equation}

At $(\varepsilon,t)=(0,0)$, the equation becomes
$$
\lambda(u_{m+N}+u_{m-N})+(2\cos(2\pi(x+m\alpha))-{E_*})u_m=0.
$$
Every initial vector $b\in\mathbb C^{2N}$ determines a unique solution. The previously defined map $\widetilde U_*(x):\mathbb C^{2N}\to X_a$ therefore satisfies
$
{\mathbf L_{0,0}^{(\infty)}}(x){\widetilde U_*}(x)=0, \mathcal E{\widetilde U_*}(x)=I_{2N},
$
and ${\widetilde U_*}$ is analytic in $x$.
For $l\in[N,\infty]$, define
$
K_{\varepsilon, t}^{(l)}(x)=\mathcal EP_{\varepsilon, t}^{(l)}(x){\widetilde U_*}(x).
$

\begin{lemma}\label{lem:initD}
After decreasing $s_*$, the matrices $K_{\varepsilon, t}^{(l)}(x)$ are invertible throughout $\mathbb T_h\times\mathbb D_{s_*}^2$, uniformly in $l\in[N,\infty]$. Moreover,
$
\widetilde U_{\varepsilon, t}^{(l)}(x)
=P_{\varepsilon, t}^{(l)}(x){\widetilde U_*}(x)\big(K_{\varepsilon, t}^{(l)}(x)\big)^{-1}
$
belongs to $\mathscr H_{ h,s}(\mathcal L(\mathbb C^{2N},X_a))$ for every $0<s\leqslant s_*$ and is a frame of $\mathcal{M}_{\varepsilon, t}^{(l)}(x)$ satisfying
\begin{equation}\label{eq:rankD}
\mathcal E\widetilde U_{\varepsilon, t}^{(l)}(x)=I_{2N},\qquad
\operatorname{Ran}\widetilde U_{\varepsilon, t}^{(l)}(x)=\mathcal{M}_{\varepsilon, t}^{(l)}(x).
\end{equation}
Let $
A_{\varepsilon, t}^{(l)}(x)=\mathcal ES^{{q}}\widetilde U_{\varepsilon, t}^{(l)}(x)
$. Then \begin{equation}\label{eq:intertwineD}
S^{{q}}\widetilde U_{\varepsilon, t}^{(l)}(x)
=\widetilde U_{\varepsilon, t}^{(l)}(x+{q}\alpha)A_{\varepsilon, t}^{(l)}(x).
\end{equation} Furthermore,
$$
\|\widetilde U_{\varepsilon, t}^{(l)}-{\widetilde U_*}\|_{ h,s}
+\|A_{\varepsilon, t}^{(l)}-{A_{*,q}}\|_{ h,s}\leqslant Cs,
\|\widetilde U_{\varepsilon, t}^{(l)}-\widetilde U_{\varepsilon, t}^{(\infty)}\|_{ h,s}
+\|A_{\varepsilon, t}^{(l)}-A_{\varepsilon, t}^{(\infty)}\|_{ h,s}
\leqslant Cs e^{-\mu l}.
$$
\end{lemma}

\begin{proof}
The projection estimates and ${P_*}{\widetilde U_*}={\widetilde U_*}$ imply
$
P_{\varepsilon, t}^{(l)}{\widetilde U_*}
={\widetilde U_*}+O_{ h,s}(s).
$
Applying the bounded map $\mathcal E$ and using $\mathcal E{\widetilde U_*}=I_{2N}$, we obtain
$
K_{\varepsilon, t}^{(l)}=I_{2N}+O_{ h,s}(s).
$
Thus, after decreasing $s_*$, the Neumann series for $\big(K_{\varepsilon, t}^{(l)}\big)^{-1}$ converges in $\mathscr H_{ h,s}(\mathrm{M}_{2N}(\mathbb C))$, uniformly for $0<s\leqslant s_*$, and
$
\big(K_{\varepsilon, t}^{(l)}\big)^{-1}
=I_{2N}+O_{ h,s}(s).
$
In particular, $\mathcal E\widetilde U_{\varepsilon, t}^{(l)}=I_{2N}$. Since $K_{\varepsilon, t}^{(l)}$ is invertible, $P_{\varepsilon, t}^{(l)}{\widetilde U_*}$ is injective. Its range is contained in the solution space and has dimension $2N$; by the preceding rank result, it therefore equals $\mathcal{M}_{\varepsilon, t}^{(l)}(x)$. This proves \eqref{eq:rankD}.

The invariance of the solution spaces under $S^{{q}}$ gives an invertible matrix $A_{\varepsilon, t}^{(l)}(x)$, which satisfies \eqref{eq:intertwineD}. Applying $\mathcal E$ and using \eqref{eq:rankD} yields $A_{\varepsilon, t}^{(l)}(x)=\mathcal ES^{{q}}\widetilde U_{\varepsilon, t}^{(l)}(x)$. At ${(\varepsilon,t)=(0,0)}$, the recurrence and the initial-data coordinates \eqref{eq:initD} give
$$
\mathcal ES^N\widetilde U_*(x)=A_{*,N}(x).
$$
{For the fixed $q$, $\mathcal ES^q\widetilde U_*(x)=A_{*,q}(x)$.}
The estimate for $\widetilde U_{\varepsilon, t}^{(l)}-{\widetilde U_*}$ follows from the estimates for the projection and $\big(K_{\varepsilon, t}^{(l)}\big)^{-1}$. Since $\mathcal ES^{{q}}:X_a\to\mathbb C^{2N}$ is bounded, 
$$
\|A_{\varepsilon, t}^{(l)}-{A_{*,q}}\|_{ h,s;\mathrm{M}_{2N}(\mathbb C)}
\leqslant\|\mathcal ES^{{q}}\|_{\mathcal L(X_a,\mathbb C^{2N})}
\|\widetilde U_{\varepsilon, t}^{(l)}-{\widetilde U_*}\|_{ h,s;\mathcal L(\mathbb C^{2N},X_a)}.
$$
This proves both the matrix estimate.

Finally, the convergence of the Riesz projections, together with
$$
(K_{\varepsilon, t}^{(l)})^{-1}-(K_{\varepsilon, t}^{(\infty)})^{-1}=(K_{\varepsilon, t}^{(l)})^{-1}(K_{\varepsilon, t}^{(\infty)}-K_{\varepsilon, t}^{(l)})(K_{\varepsilon, t}^{(\infty)})^{-1},
$$
 gives the estimates for $\widetilde U_{\varepsilon, t}^{(l)}-\widetilde U_{\varepsilon, t}^{(\infty)}$ and $A_{\varepsilon, t}^{(l)}-A_{\varepsilon, t}^{(\infty)}$.
\end{proof}

\end{proof}

\section{Symplectification, monotonicity and dynamical invariants}
We first symplectify the frames constructed in the preceding section, then establish monotonicity through the energy and potential derivatives. Finally, we identify the slow Lyapunov exponents and relate the rotation number to the IDS.
\subsection{Symplectification}

\begin{proposition}\label{prop:slow}
Let $\alpha\in\mathbb{R}\backslash\mathbb{Q}$ and fix an integer $q\geqslant1$.  There are $s_{*},h,b,C>0$, uniform in $E_{*}\in\Sigma_{0}$, such that, for $l\in[N,\infty]$ and $0<s\leqslant s_{*}$, the spaces $\mathcal{M}_{\varepsilon,t}^{(l)}(x)$ have dimension $2N$ and admit jointly analytic frames $U_{\varepsilon,t}^{(l)}(x)$ satisfying
	$$
	S^{q}U_{\varepsilon,t}^{(l)}(x)=U_{\varepsilon,t}^{(l)}(x+q\alpha)C_{\varepsilon,t}^{(l)}(x),\qquad
	C_{\varepsilon,t}^{(l)}(x)=A_{*,q}(x)e^{F_{\varepsilon,t}^{(l)}(x)}.
	$$
	For real $x,\varepsilon,t$, one has $C_{\varepsilon,t}^{(l)}\in\mathrm{HSp}(J_{2N})$ and $F_{\varepsilon,t}^{(l)}\in\mathfrak{hsp}(J_{2N})$. Moreover,
	\begin{equation*}\label{slowbd}
		\|U_{\varepsilon,t}^{(l)}-\widetilde U_{*}\|_{h,s}+\|F_{\varepsilon,t}^{(l)}\|_{h,s}\leqslant Cs,
	\end{equation*}
	and, for $|\varepsilon|,|t|<s_{*}$,
	\begin{equation*}\label{slowtail}
		\begin{aligned}
			&\|U_{\varepsilon,t}^{(l)}-U_{\varepsilon,t}^{(\infty)}\|_{h}
			+\|C_{\varepsilon,t}^{(l)}-C_{\varepsilon,t}^{(\infty)}\|_{h}
			+\|F_{\varepsilon,t}^{(l)}-F_{\varepsilon,t}^{(\infty)}\|_{h}
			\leqslant C|\varepsilon|e^{-bl}.
		\end{aligned}
	\end{equation*}
\end{proposition}

Recall $t_k$ from \eqref{equ:part}.
For $l\in[N,\infty]$ and real $\varepsilon$, define the symplectic form on $X_a$ by
\begin{equation}\label{green}
\mathfrak J_\varepsilon^{(l)}(u,v)
=\frac1\lambda\sum_{k=1}^l\sum_{j=1}^k
\big(t_k\overline{u_{-j}}v_{k-j}-\overline{t_k}\overline{u_{k-j}}v_{-j}\big).
\end{equation}
Since $a<2\pi r$, the series in \eqref{green} converges absolutely, and
$$
|\mathfrak J_\varepsilon^{(l)}(u,v)|
\leqslant\frac2\lambda\|u\|_a\|v\|_a\sum_{k=1}^l k|t_k|e^{ak}
\leqslant C\|u\|_a\|v\|_a,
$$
uniformly in $l$ and small real $\varepsilon$.
The form is skew-Hermitian and independent of $x,t$; at $\varepsilon=0$, it is independent of $l$ and denoted by $\mathfrak J_0$.
For frames, we use the same notation columnwise:
$$
\big(\mathfrak J_\varepsilon^{(l)}(U,V)\big)_{ij}
=\mathfrak J_\varepsilon^{(l)}(U\mathbf e_i,V\mathbf e_j).
$$

For finite $l$, we always take $t_l\ne0$, with $l=N$ when $\varepsilon=0$. Put $v=v_0+\varepsilon W^{(l)}$ and use $C_v,J_v$ from Section~\ref{sec:aubry-duality}. In the coordinates $\mathcal E_lu$, the normalized symplectic form satisfies
\begin{equation}\label{eq:Jld}
\mathfrak J_\varepsilon^{(l)}(u,w)
=\lambda^{-1}(\mathcal E_lu)^*J_v(\mathcal E_lw).
\end{equation}
The one-step companion cocycle $\mathcal A_{\varepsilon,t}^{(l)}=\mathcal A_{E_*+t}^{v}$ preserves $J_v$:
\begin{equation}\label{eq:finiteJD}
\mathcal A_{\varepsilon,t}^{(l)}(x)^*J_v\mathcal A_{\varepsilon,t}^{(l)}(x)=J_v,
\qquad x,\varepsilon,t\in\mathbb R.
\end{equation}

Let
$$
\widetilde\Phi_l(x)=\mathcal E_l\widetilde U_{\varepsilon, t}^{(l)}(x):\mathbb C^{2N}\to\mathbb C^{2l}.
$$
This map is injective, because $l\geqslant N$ and $\mathcal E\widetilde U_{\varepsilon, t}^{(l)}=I_{2N}$. The companion relation and \eqref{eq:intertwineD} imply
\begin{equation}\label{eq:bridgeD}
\big(\mathcal A_{\varepsilon, t}^{(l)}\big)^{({q})}(x)\widetilde\Phi_l(x)=\widetilde\Phi_l(x+{q}\alpha)A_{\varepsilon, t}^{(l)}(x).
\end{equation}
At $\varepsilon=0$, this identity holds with $l=N$ and $\widetilde\Phi_N=I_{2N}$.
For $l\in[N,\infty]$ and real $x,\varepsilon,t$, define
\begin{equation}\label{eq:GlD}
G_{\varepsilon,t}^{(l)}(x)
=\mathfrak J_{\varepsilon}^{(l)}\big(\widetilde U_{\varepsilon,t}^{(l)}(x),\widetilde U_{\varepsilon,t}^{(l)}(x)\big).
\end{equation}
For finite $l$, \eqref{eq:Jld} gives
$$
G_{\varepsilon,t}^{(l)}(x)=\lambda^{-1}\widetilde\Phi_{l}(x)^{*}J_v\widetilde\Phi_{l}(x).
$$
By \eqref{eq:finiteJD} and \eqref{eq:bridgeD},
\begin{equation}\label{eq:GinvD}
A_{\varepsilon,t}^{(l)}(x)^{*}G_{\varepsilon,t}^{(l)}(x+q\alpha)A_{\varepsilon,t}^{(l)}(x)=G_{\varepsilon,t}^{(l)}(x).
\end{equation}
At $\varepsilon=0$, the slow cocycles and the pulled-back forms are independent of $l$, so \eqref{eq:GinvD} holds for every $l\geqslant N$ as well.

\begin{lemma}\label{lem:formD}
The matrices $G_{\varepsilon, t}^{(l)}$ admit holomorphic extensions to $\mathbb T_h\times\mathbb D_{s_*}^2$ and converge in $\mathscr H_{ h,s}(\mathrm{M}_{2N}(\mathbb C))$ to $G_{\varepsilon, t}^{(\infty)}$. For $l\in[N,\infty]$ and every $0<s\leqslant s_*$,
\begin{equation}\label{eq:GbdD}
\|G_{\varepsilon, t}^{(l)}-J_{2N}\|_{ h,s}\leqslant Cs,\qquad \|G_{\varepsilon, t}^{(l)}-G_{\varepsilon, t}^{(\infty)}\|_{ h,s}\leqslant Cs e^{-\mu l}.
\end{equation}
Moreover, for real $x,\varepsilon, t$,
\begin{equation}\label{eq:infiniteJD}
A_{\varepsilon, t}^{(\infty)}(x)^*G_{\varepsilon, t}^{(\infty)}(x+{q}\alpha)A_{\varepsilon, t}^{(\infty)}(x)=G_{\varepsilon, t}^{(\infty)}(x).
\end{equation}
\end{lemma}

\begin{proof}
	Write $\widetilde U_l=\widetilde U_{\varepsilon, t}^{(l)}$ and let $\widetilde U_{l,j}\in\mathbb C^{1\times2N}$ be its $j$th coordinate row and $\widetilde U_{l,j}^{\sharp}(x,\varepsilon, t):=\widetilde U_{l,j}(\overline x,\overline\varepsilon,\overline t)^*.$ Directly from \eqref{eq:GlD} and \eqref{green}, using $\mathcal E\widetilde U_l=I_{2N}$, we obtain
	\begin{equation}\label{eq:GseriesD}
		\begin{aligned}
			G_{\varepsilon, t}^{(l)}=J_{2N}+\frac{\varepsilon}{\lambda}\sum_{k=1}^{l}\sum_{j=0}^{k-1}\big(&-\overline{w_k}\,\widetilde U_{l,j}^{\sharp}\widetilde U_{l,j-k}
			+w_k\,\widetilde U_{l,j-k}^{\sharp}\widetilde U_{l,j}\big).
		\end{aligned}
	\end{equation}
	 This formula agrees with the finite-range definition for real $x,\varepsilon, t$.
	And \eqref{eq:GseriesD} gives
	$$
	\|G_{\varepsilon, t}^{(l)}- {J_{2N}}\|_{ h,s}
	\leqslant\frac{2s}{\lambda}\|\widetilde U_l\|_{ h,s}^{2}
	\sum_{k=1}^{l}k|w_k|e^{ak}\leqslant Cs.
	$$
Therefore the series
	$$
	\begin{aligned}
		G_{\varepsilon, t}^{(\infty)}=J_{2N}+\frac{\varepsilon}{\lambda}\sum_{k\geqslant1}\sum_{j=0}^{k-1}\big(&-\overline{w_k}\,\widetilde U_{\infty,j}^{\sharp}\widetilde U_{\infty,j-k}
		+w_k\,\widetilde U_{\infty,j-k}^{\sharp}\widetilde U_{\infty,j}\big)
	\end{aligned}
	$$
	converges absolutely in $\mathscr H_{ h,s}(\mathrm{M}_{2N}(\mathbb C))$. 
	Using \eqref{eq:GseriesD}, the corresponding infinite series, and Lemma~\ref{lem:initD}, we obtain
	$$
	\begin{aligned}
		\|G_{\varepsilon, t}^{(l)}-G_{\varepsilon, t}^{(\infty)}\|_{ h,s}\leqslant Cs\big(\|\widetilde U_l-\widetilde U_\infty\|_{ h,s}+e^{-\mu l}\big)\leqslant Cs e^{-\mu l}.
	\end{aligned}
	$$ This proves \eqref{eq:GbdD}.
	
	Finally, by Lemma~\ref{lem:initD}, $A_{\varepsilon, t}^{(l)}\to A_{\varepsilon, t}^{(\infty)}$ in the same norm. Passing to the limit in \eqref{eq:GinvD} for real variables gives \eqref{eq:infiniteJD}.
\end{proof}

\begin{lemma}\label{lem:sympD}
After decreasing $s_*$, for $l\in[N,\infty]$ set
\begin{equation*}\label{eq:QD}
Q_{\varepsilon, t}^{(l)}(x)=\big(-J_{2N}G_{\varepsilon, t}^{(l)}(x)\big)^{-1/2},\qquad U_{\varepsilon, t}^{(l)}(x)=\widetilde U_{\varepsilon, t}^{(l)}(x)Q_{\varepsilon, t}^{(l)}(x).
\end{equation*}
and
\begin{equation*}\label{eq:CD}
C_{\varepsilon, t}^{(l)}(x)
=Q_{\varepsilon, t}^{(l)}(x+{q}\alpha)^{-1}
A_{\varepsilon, t}^{(l)}(x)Q_{\varepsilon, t}^{(l)}(x).
\end{equation*}
	Then $U_{\varepsilon,t}^{(l)}$ and $C_{\varepsilon,t}^{(l)}$ satisfy Proposition~\ref{prop:slow}.  Moreover,
\begin{equation}\label{eq:normalbdD}
\begin{aligned}
&\|U_{\varepsilon, t}^{(l)}-\widetilde U_*\|_{ h,s}+\|Q_{\varepsilon, t}^{(l)}-I_{2N}\|_{ h,s}
+\|C_{\varepsilon, t}^{(l)}-A_{*,q}\|_{ h,s}
+\|F_{\varepsilon, t}^{(l)}\|_{ h,s}\leqslant Cs,\\
&\|U_{\varepsilon, t}^{(l)}-U_{\varepsilon, t}^{(\infty)}\|_{ h,s}+\|Q_{\varepsilon, t}^{(l)}-Q_{\varepsilon, t}^{(\infty)}\|_{ h,s}
+\|C_{\varepsilon, t}^{(l)}-C_{\varepsilon, t}^{(\infty)}\|_{ h,s}
+\|F_{\varepsilon, t}^{(l)}-F_{\varepsilon, t}^{(\infty)}\|_{ h,s}
\leqslant Cs e^{-\mu l}.
\end{aligned}
\end{equation}
\end{lemma}

\begin{proof}
Let $T_l=-J_{2N}G_{\varepsilon, t}^{(l)}-I_{2N}=-J_{2N}(G_{\varepsilon, t}^{(l)}-J_{2N})$. By Lemma~\ref{lem:formD},
$$
\|T_l\|_{ h,s}\leqslant Cs,
\qquad
\|T_l-T_\infty\|_{ h,s}\leqslant Cs e^{-\mu l}.
$$
Decreasing $s_*$ so that $Cs_*<1/4$, define
$$
Q_{\varepsilon, t}^{(l)}
=\sum_{n\geqslant0}\binom{-1/2}{n}T_l^n,
\qquad
\big(Q_{\varepsilon, t}^{(l)}\big)^{-1}
=\sum_{n\geqslant0}\binom{1/2}{n}T_l^n.
$$
Both series converge absolutely in the Banach algebra $\mathscr H_{ h,s}(\mathrm{M}_{2N}(\mathbb C))$, uniformly in $l$ and $0<s\leqslant s_*$. In particular,
$$
\|Q_{\varepsilon, t}^{(l)}-I_{2N}\|_{ h,s}
+\|\big(Q_{\varepsilon, t}^{(l)}\big)^{-1}-I_{2N}\|_{ h,s}
\leqslant Cs.
$$
For $n\geqslant1$, the identity
$
T_l^n-T_\infty^n
=\sum_{j=0}^{n-1}T_l^j(T_l-T_\infty)T_\infty^{n-1-j}
$
 give
$$
\|T_l^n-T_\infty^n\|_{ h,s}
\leqslant n4^{-(n-1)}\|T_l-T_\infty\|_{ h,s}.
$$
Summing the binomial series yields
\begin{equation*}\label{eq:QdiffD}
\|Q_{\varepsilon, t}^{(l)}-Q_{\varepsilon, t}^{(\infty)}\|_{ h,s}
+\|\big(Q_{\varepsilon, t}^{(l)}\big)^{-1}
-\big(Q_{\varepsilon, t}^{(\infty)}\big)^{-1}\|_{ h,s}
\leqslant Cs e^{-\mu l}.
\end{equation*}

On the real axis, $G_{\varepsilon, t}^{(l)}$ and $J_{2N}$ are skew-Hermitian and
$
\big(-J_{2N}G_{\varepsilon, t}^{(l)}\big)^{*}J_{2N}=J_{2N}\big(-J_{2N}G_{\varepsilon, t}^{(l)}\big).
$
The binomial series has real coefficients, so $Q_{\varepsilon, t}^{(l)}{}^{*}J_{2N}=J_{2N}Q_{\varepsilon, t}^{(l)}$. Since $Q_{\varepsilon, t}^{(l)}$ commutes with $-J_{2N}G_{\varepsilon, t}^{(l)}$, we conclude that
\begin{equation}\label{equ:sym} 
Q_{\varepsilon, t}^{(l)}{}^{*}G_{\varepsilon, t}^{(l)}Q_{\varepsilon, t}^{(l)}=Q_{\varepsilon, t}^{(l)}{}^*J_{2N}\big(-J_{2N}G_{\varepsilon, t}^{(l)}\big)Q_{\varepsilon, t}^{(l)}=J_{2N}\big(Q_{\varepsilon, t}^{(l)}\big)^2\big(-J_{2N}G_{\varepsilon, t}^{(l)}\big)=J_{2N}.
\end{equation}
	Combining \eqref{eq:GinvD} with \eqref{equ:sym} proves that $C_{\varepsilon,t}^{(l)}$ preserves $J_{2N}$. Its transfer relation follows from \eqref{eq:intertwineD}. Since $G_{0,t}^{(l)}=J_{2N}$, we have $Q_{0,t}^{(l)}=I_{2N}$ and $C_{0,0}^{(l)}=A_{*,q}$.

We now prove the uniform estimates for the cocycle. The preceding estimates for $Q_{\varepsilon,t}^{(l)}$ give
$$
\|Q_{\varepsilon,t}^{(l)}-I_{2N}\|_{h,s}\leqslant Cs,\qquad
\|Q_{\varepsilon,t}^{(l)}-Q_{\varepsilon,t}^{(\infty)}\|_{h,s}\leqslant Cs e^{-\mu l}.
$$
Since $U_{\varepsilon,t}^{(l)}=\widetilde U_{l}Q_{\varepsilon,t}^{(l)}$, these bounds and Lemma~\ref{lem:initD} also give
$$
\|U_{\varepsilon,t}^{(l)}-\widetilde U_{*}\|_{h,s}\leqslant Cs,\qquad
\|U_{\varepsilon,t}^{(l)}-U_{\varepsilon,t}^{(\infty)}\|_{h,s}\leqslant Cs e^{-\mu l}.
$$
The corresponding estimates for the inverses of $Q_{l}=Q_{\varepsilon,t}^{(l)}$ follow from
$$
Q_{l}^{-1}-I_{2N}=Q_{l}^{-1}(I_{2N}-Q_{l}),\qquad
Q_{l}^{-1}-Q_{\infty}^{-1}=Q_{l}^{-1}(Q_{\infty}-Q_{l})Q_{\infty}^{-1}.
$$
Since $A_{\varepsilon,t}^{(l)}=A_{*,q}+O_{h,s}(s)$, the formula
$$
C_{\varepsilon,t}^{(l)}(x)=Q_{l}(x+q\alpha)^{-1}A_{\varepsilon,t}^{(l)}(x)Q_{l}(x)
$$
gives $\|C_{\varepsilon,t}^{(l)}-A_{*,q}\|_{h,s}\leqslant Cs$. Comparing this formula for $l$ and $\infty$, and using the truncation estimate in Lemma~\ref{lem:initD}, gives
$
\|C_{\varepsilon,t}^{(l)}-C_{\varepsilon,t}^{(\infty)}\|_{h,s}\leqslant Cs e^{-\mu l}.
$

Set $Y_{l}(x)=A_{*,q}(x)^{-1}C_{\varepsilon,t}^{(l)}(x)-I_{2N}$. Since $A_{*,q}^{-1}$ is uniformly bounded, we have
$$
\|Y_{l}\|_{h,s}\leqslant Cs,\qquad \|Y_{l}-Y_{\infty}\|_{h,s}\leqslant Cs e^{-\mu l}.
$$
After decreasing $s_{*}$, the logarithm series
$$
F_{\varepsilon,t}^{(l)}=\sum_{n\geqslant1}\frac{(-1)^{n+1}}{n}Y_{l}^{n}
$$
converges absolutely in the same norm. Applying again
$$
Y_{l}^{n}-Y_{\infty}^{n}
=\sum_{j=0}^{n-1}Y_{l}^{j}(Y_{l}-Y_{\infty})Y_{\infty}^{n-1-j}
$$
gives
$$
\|F_{\varepsilon,t}^{(l)}\|_{h,s}\leqslant Cs,\qquad
\|F_{\varepsilon,t}^{(l)}-F_{\varepsilon,t}^{(\infty)}\|_{h,s}\leqslant Cs e^{-\mu l}.
$$
Since both $A_{*,q}(x)$ and $C_{\varepsilon,t}^{(l)}(x)$ preserve $J_{2N}$ for real variables, the local logarithm satisfies
$
(F_{\varepsilon,t}^{(l)})^{*}J_{2N}+J_{2N}F_{\varepsilon,t}^{(l)}=0
$
for real variables. This completes \eqref{eq:normalbdD}. In particular, $F_{\varepsilon,t}^{(l)}$ is jointly holomorphic in $(x,\varepsilon,t)$ on a common domain, and $F_{0,0}^{(l)}=0$.
\end{proof}
\subsection{Infinite-range discrete
Hellmann--Feynman identity}

The Hellmann--Feynman identity expresses the derivative of an
eigenvalue along a differentiable normalized eigenpair
$H_su_s=E_su_s$ as $E_s'=\langle u_s,H_s'u_s\rangle$.
Thus the first-order response to a perturbation can be read from
the unperturbed state. For a finite-range difference operator,
Li--Xu--Zhou \cite[Lemma~5.9]{LXZ26} obtain a discrete analogue:
the energy variation of an iterated transfer matrix, paired with
the symplectic form, is a sum of squared solution amplitudes.
Its positivity gives monotonicity in energy after finitely many
steps.

We derive the corresponding identity on the slow solution bundle
of the infinite-range Aubry dual, allowing both the energy and the
potential to vary. We omit the superscript $(\infty)$. Fix $W$ and
a sufficiently small real parameter point $(\bar\varepsilon,\bar t)$,
and write $U=U_{\bar\varepsilon,\bar t}$ and
$C=C_{\bar\varepsilon,\bar t}$. For
$\mathbf v=(v_\varepsilon,v_t)\in\mathbb R^2$, let
$$
 \delta_{\mathbf v}
 =v_\varepsilon\partial_\varepsilon+v_t\partial_t
 \quad\text{at }(\bar\varepsilon,\bar t).
$$
These derivatives exist by Proposition~\ref{prop:slow}. The directions
$(0,1)$ and $(1,0)$ vary the energy and the potential, respectively.
At the origin, $(1,-\xi)$ is tangent to the line $t=-\varepsilon\xi$.

Differentiating $C^*J_{2N}C=J_{2N}$ shows that
$C^*J_{2N}\delta_{\mathbf v}C$ is Hermitian. To account for the
variation of the frame, set
$$
 K_{\mathbf v}(x)
 =\mathfrak J_{\bar\varepsilon}
    \bigl(U(x),\delta_{\mathbf v}U(x)\bigr).
$$
Then we have the following infinite-range discrete
Hellmann--Feynman identity:

\begin{proposition}\label{energyvar}
One has
\begin{equation}\label{varid}
		\begin{aligned}
			C(x)^{*}J_{2N}\delta_{\mathbf{v}}C(x)
			&=\frac{1}{\lambda}\sum_{m=0}^{q-1}U_{m}(x)^{*}\bigg(v_{t}U_{m}(x)-v_{\varepsilon}\sum_{n\in\mathbb{Z}}w_{n}U_{m+n}(x)\bigg)\\
			&\quad+K_{\mathbf{v}}(x)-C(x)^{*}K_{\mathbf{v}}(x+q\alpha)C(x).
		\end{aligned}
\end{equation}

\end{proposition}

\begin{proof}
Since $t=E-E_{*}$, the operator satisfies
$$
\mathbf{L}_{\varepsilon,t}(x)
=\mathbf{L}_{0,0}(x)+\varepsilon\sum_{n\in\mathbb{Z}}w_{n}S^{n}-tI.
$$
Thus, along $(\bar\varepsilon+sv_{\varepsilon},\bar t+sv_{t})$,
$$
\frac{\mathrm{d}}{\mathrm{d}s}\mathbf{L}_{\bar\varepsilon+sv_{\varepsilon},\bar t+sv_{t}}(x)\bigg|_{s=0}
=v_{\varepsilon}\sum_{n\in\mathbb{Z}}w_{n}S^{n}-v_{t}I.
$$
By Proposition~\ref{prop:slow}, the frame is jointly analytic in the parameters and satisfies
$$
\mathbf{L}_{\bar\varepsilon+sv_{\varepsilon},\bar t+sv_{t}}(x)
U_{\bar\varepsilon+sv_{\varepsilon},\bar t+sv_{t}}(x)=0.
$$
Differentiating this identity at $s=0$ gives
$$
0=\bigg(v_{\varepsilon}\sum_{n\in\mathbb{Z}}w_{n}S^{n}-v_{t}I\bigg)U(x)
+\mathbf{L}_{\bar\varepsilon,\bar t}(x)\delta_{\mathbf{v}}U(x).
$$
Hence
$$
\mathbf{L}_{\bar\varepsilon,\bar t}(x)\delta_{\mathbf{v}}U(x)
=v_{t}U(x)-v_{\varepsilon}\sum_{n\in\mathbb{Z}}w_{n}S^{n}U(x).
$$

At $(\bar\varepsilon,\bar t)$, a direct calculation gives
$$
\mathfrak J_{\bar\varepsilon}(Su,Sv)-\mathfrak J_{\bar\varepsilon}(u,v)
=\frac{1}{\lambda}\big(\overline{u_{0}}(\mathbf{L}_{\bar\varepsilon,\bar t}(x)v)_{0}
-\overline{(\mathbf{L}_{\bar\varepsilon,\bar t}(x)u)_{0}}v_{0}\big).
$$
Apply this identity to $U$ and $\delta_{\mathbf{v}}U$ at successive lattice sites. Since $\mathbf{L}_{\bar\varepsilon,\bar t}(x)U(x)=0$, summation over $m=0,\ldots,q-1$ gives
$$
\mathfrak J_{\bar\varepsilon}(S^{q}U(x),S^{q}\delta_{\mathbf{v}}U(x))-K_{\mathbf{v}}(x)=\frac{1}{\lambda}\sum_{m=0}^{q-1}U_{m}(x)^{*}\bigg(v_{t}U_{m}(x)-v_{\varepsilon}\sum_{n\in\mathbb{Z}}w_{n}U_{m+n}(x)\bigg).
$$
The first variation of the transfer relation is
$$
S^{q}\delta_{\mathbf{v}}U(x)
=(\delta_{\mathbf{v}}U)(x+q\alpha)C(x)+U(x+q\alpha)\delta_{\mathbf{v}}C(x).
$$
Since
$$
S^{q}U(x)=U(x+q\alpha)C(x),\qquad
\mathfrak J_{\bar\varepsilon}(U(x+q\alpha),U(x+q\alpha))=J_{2N},
$$
the terminal pairing becomes
$$
\begin{aligned}
\mathfrak J_{\bar\varepsilon}(S^{q}U(x),S^{q}\delta_{\mathbf{v}}U(x))
&=C(x)^{*}K_{\mathbf{v}}(x+q\alpha)C(x)
+C(x)^{*}J_{2N}\delta_{\mathbf{v}}C(x).
\end{aligned}
$$
Substitution proves \eqref{varid}. 
\end{proof}

\begin{corollary}\label{cor:modevar}
Suppose that $C(x)=C_{*}$ is constant.
For $z\neq0$ with $C_{*}z=\zeta z$ and $|\zeta|=1$,\footnote{This is the observation used in \cite[\S4.2]{JLS26}.} one has
$$
z^{*}C_{*}^{*}J_{2N}\langle\delta_{\mathbf{v}}C\rangle z
=\frac{1}{\lambda}\sum_{m=0}^{q-1}
\int_{\mathbb{T}}\overline{U_{m}(x)z}
\bigg(v_{t}U_{m}(x)z
-v_{\varepsilon}\sum_{n\in\mathbb{Z}}w_{n}U_{m+n}(x)z\bigg)\,\mathrm{d}x.
$$
For the energy direction, set $G=C_{*}^{-1}\delta_{(0,1)}C$. Then
\begin{equation}\label{modepos}
z^{*}J_{2N}\langle G\rangle z
=\frac{1}{\lambda}\sum_{m=0}^{q-1}
\int_{\mathbb{T}}|U_{m}(x)z|^{2}\,\mathrm{d}x>0.
\end{equation}
\end{corollary}

\begin{proof}
Since $C_{*}z=\zeta z$ and $|\zeta|=1$, translation invariance gives
$$
\int_{\mathbb{T}}z^{*}
\big(K_{\mathbf{v}}(x)-C_{*}^{*}K_{\mathbf{v}}(x+q\alpha)C_{*}\big)z\,\mathrm{d}x=0.
$$
Taking the phase average of \eqref{varid} against $z$ proves the first identity. For $\mathbf{v}=(0,1)$, the relation $C_{*}^{*}J_{2N}=J_{2N}C_{*}^{-1}$ gives the equality in \eqref{modepos}.

If its right-hand side vanished, continuity would give $U_{m}(x)z=0$ for every $x$ and $0\leqslant m<q$. The relation $S^{q}U(x)z=\zeta U(x+q\alpha)z$ then gives $U_{m}(x)z=0$ for every $m\in\mathbb{Z}$, contradicting injectivity of the frame.
\end{proof}

\begin{corollary}\label{parmono}
Let $L=\lceil 2N/q\rceil$. For sufficiently small real $\varepsilon,t$ and every $x\in\mathbb T$, one has
	\begin{equation}\label{enmono}
		\big((C_{\varepsilon,t}^{(l)})^{(L)}(x)\big)^{*}J_{2N}\partial_{t}(C_{\varepsilon,t}^{(l)})^{(L)}(x)
		\geqslant c_{0}I_{2N},
	\end{equation}
where $c_0>0$ is independent of $E_*\in\Sigma_0$ and $l\in[N,\infty]$.
\end{corollary}

\begin{proof}
	At $\varepsilon=0$, one has $\mathfrak J_{0}(u,v)=(\mathcal{E}u)^{*}J_{2N}\mathcal{E}v$. Since $\mathcal{E}\widetilde U_{0,t}^{(l)}=I_{2N}$, we have $G_{0,t}^{(l)}=J_{2N}$ and $Q_{0,t}^{(l)}=I_{2N}$. Consequently, $\mathcal{E}U_{0,t}^{(l)}=I_{2N}$ and $\mathcal{E}\partial_{t}U_{0,t}^{(l)}=0$, so
	$
	\mathfrak J_{0}(U_{0,t}^{(l)},\partial_{t}U_{0,t}^{(l)})=0.
	$
	The energy variation formula in Proposition~\ref{energyvar} applies to every truncation with the form $\mathfrak J_{\varepsilon}^{(l)}$, since $\partial_{t}\mathbf{L}_{\varepsilon,t}^{(l)}=-I$. Put $p=qL\geqslant2N$. Apply that formula to the $p$-step transfer relation at $(\bar\varepsilon,\bar t)=(0,0)$. The two terms involving $\mathfrak J_{0}(U_{0,t}^{(l)},\partial_{t}U_{0,t}^{(l)})$ vanish. 
This gives
	$$
	\left.
	\big((C_{\varepsilon,t}^{(l)})^{(L)}(x)\big)^{*}J_{2N}
	\partial_{t}(C_{\varepsilon,t}^{(l)})^{(L)}(x)
	\right|_{(\varepsilon,t)=(0,0)}
	=\frac{1}{\lambda}\sum_{m=0}^{p-1}U_{*,m}(x)^{*}U_{*,m}(x),
	$$
	where $U_*(x)=U_{0,0}^{(l)}(x)$ and $(C_{\varepsilon,t}^{(l)})^{(L)}$ is the $L$-step iterate over $q\alpha$.
	
	Let $R_{*}(x)$ be the matrix with rows $U_{*,0}(x),\ldots,U_{*,2N-1}(x)$. It differs from $\mathcal{E}S^{N}U_{*}(x)=A_{*,N}(x)$ by a fixed row permutation. The reference matrices $A_{*,N}(x)$ are invertible and depend continuously on $(E_{*},x)$. Therefore
	$
	d_{*}:=\bigg(\max_{E_{*}\in\Sigma_{0},\,x\in\mathbb{T}}\|R_{*}(x)^{-1}\|\bigg)^{-1}>0.
	$
	Since $p\geqslant2N$, it follows that
	$$
	\frac{1}{\lambda}\sum_{m=0}^{p-1}U_{*,m}(x)^{*}U_{*,m}(x)
	\geqslant\frac{1}{\lambda}R_{*}(x)^{*}R_{*}(x)
	\geqslant\frac{d_{*}^{2}}{\lambda}I_{2N}.
	$$
	The holomorphic estimates in Lemma~\ref{lem:sympD} are uniform in $E_{*}$ and $l$ on a fixed parameter domain. For the fixed iterate $L$, Cauchy estimates on a smaller domain bound its first and second parameter derivatives uniformly. Thus, for real $\varepsilon,t$,
	$$
	\bigg\|
	\big((C_{\varepsilon,t}^{(l)})^{(L)}\big)^{*}J_{2N}\partial_{t}(C_{\varepsilon,t}^{(l)})^{(L)}
	-\frac{1}{\lambda}\sum_{m=0}^{p-1}U_{*,m}^{*}U_{*,m}
	\bigg\|_{C^{0}(\mathbb{T})}\leqslant K(|\varepsilon|+|t|).
	$$
	The matrix on the left of \eqref{enmono} is Hermitian by differentiation of the preservation identity for $J_{2N}$. After decreasing the parameter domain so that $K(|\varepsilon|+|t|)\leqslant d_{*}^{2}/(2\lambda)$, the preceding bounds prove \eqref{enmono} with $c_{0}=d_{*}^{2}/(2\lambda)$.
\end{proof}

\subsection{Lyapunov exponents and rotation number}

\subsubsection{Dominated splitting and acceleration}

\begin{lemma}\label{lem:par}
	For sufficiently small real $\varepsilon,t$ with $\varepsilon\neq0$
and finite $l$,  the companion cocycle $\mathcal A_{\varepsilon,t}^{(l)}$ admits a dominated splitting
	\begin{equation*}
		\mathbb C^{2l}=E_l^u(x)\oplus E_l^c(x)\oplus E_l^s(x),\qquad \dim E_l^c=2N,\quad \dim E_l^u=\dim E_l^s=l-N,
	\end{equation*}
    and the bundles $E_{l}^{u}$, $E_{l}^{c}$, and $E_{l}^{s}$ depend continuously on $(x,\varepsilon,t)$.
	Moreover,
$		E_l^c(x)=\operatorname{Ran}\big(\mathcal E_l U_{\varepsilon,t}^{(l)}(x)\big), x\in\mathbb T. $
    And the Lyapunov exponents restricted to $E_l^c(x)$ have absolute value less than $a/2$, and all remaining exponents have absolute value not smaller than $a$.

As a consequence, let $L_{1}(q\alpha,C_{\varepsilon,t}^{(\infty)})\geqslant\cdots\geqslant L_{2N}(q\alpha,C_{\varepsilon,t}^{(\infty)})$ be the Lyapunov exponents, then
\begin{equation}\label{eq:accD}
\#\{j:L_{j}(q\alpha,C_{\varepsilon,t}^{(\infty)})=0\}=2\omega_{\varepsilon,W}(E_*+t).
\end{equation}
\end{lemma}

\begin{proof}
	\eqref{eq:reference-growth} gives
	$$
	\sup_{E_{*}\in\Sigma_{0},\,x\in\mathbb{T}}\|A_{*,q}^{(n)}(x)\|\leqslant C e^{q\chi|n|},\qquad n\in\mathbb{Z}.
	$$
	 Choose $\chi<\chi_{1}<\chi_{2}<a/2$ and fix an integer $m\geqslant1$ with $Ce^{q\chi m}\leqslant e^{q\chi_{1}m}$. By Proposition~\ref{prop:slow}, the cocycle $C_{\varepsilon,t}^{(l)}$ and its inverse are uniformly close to $A_{*,q}$ and its inverse. Hence, after decreasing the parameter domain, both iterates of length $m$ satisfy
	$$
	\sup_{E_{*}\in\Sigma_{0},\,x\in\mathbb{T},\,l\in[N,\infty]}
	\|C_{\varepsilon,t}^{(l),(\pm m)}(x)\|\leqslant e^{q\chi_{2}m}.
	$$
	Splitting each positive or negative iterate into blocks of length $m$ and a bounded remainder gives
	\begin{equation}\label{eq:slow}
		\sup_{x\in\mathbb{T}}
		\big\|C_{\varepsilon,t}^{(l),(n)}(x)\big\|
		\leqslant K e^{q\chi_{2}|n|},
		\qquad l\in[N,\infty],
	\end{equation}
	where $K>0$ is independent of $E_{*}$ and $l$. Thus every Lyapunov exponent of
	$(q\alpha,C_{\varepsilon,t}^{(l)})$ belongs to
	$[-q\chi_{2},q\chi_{2}]$.
	
	Fix a finite $l$ and write
	$$
	\mathcal{C}_{l}(x)
	:=\operatorname{Ran}\big(\mathcal{E}_{l}U_{\varepsilon,t}^{(l)}(x)\big)
	=\mathcal{E}_{l}\ker_{X_{a}}\mathbf{L}_{\varepsilon,t}^{(l)}(x).
	$$
	By Proposition~\ref{prop:slow} and covariance under $S$, this is a continuous invariant subbundle of dimension $2N$. The frame relation gives, for every
	$n\in\mathbb{Z}$,  the perturbed companion cocycle satisfies
	$$
	\mathcal{A}_{\varepsilon,t}^{(l),(qn)}(x)
	\mathcal{E}_{l}U_{\varepsilon,t}^{(l)}(x)=
	\mathcal{E}_{l}U_{\varepsilon,t}^{(l)}(x+qn\alpha)
	C_{\varepsilon,t}^{(l),(n)}(x).
	$$
	For fixed $l,\varepsilon,t$, the frame $\mathcal{E}_{l}U_{\varepsilon,t}^{(l)}$
	and its inverse on its range are uniformly bounded in $x$. Together
	with \eqref{eq:slow}, this gives
	$$
	\big\|\mathcal{A}_{\varepsilon,t}^{(l),(qn)}(x)v\big\|
	\leqslant K_{l}e^{q\chi_{2}|n|}\|v\|,
	\qquad v\in\mathcal{C}_{l}(x),\quad n\in\mathbb{Z}.
	$$
	The finitely many iterates between successive returns can be absorbed
	into the constant. Consequently,
	$$
	\big\|\mathcal{A}_{\varepsilon,t}^{(l),(m)}(x)v\big\|
	\leqslant K_{l}e^{\chi_{2}|m|}\|v\|,
	\qquad v\in\mathcal{C}_{l}(x),\quad m\in\mathbb{Z}.
	$$
	
	At an Oseledets-regular phase $x$, let
	$\mathcal{E}_{\chi}^{(l)}(x)$ be the Oseledets subspace of the
	companion cocycle corresponding to the exponent $\chi$. Since
	$\chi_{2}<a/2$, the preceding two-sided bound gives
	$
	\mathcal{C}_{l}(x)\subseteq
	\bigoplus_{|\chi|\leqslant a/2}\mathcal{E}_{\chi}^{(l)}(x).
	$
	Conversely, let $v\in\mathcal{E}_{\chi}^{(l)}(x)$ with $|\chi|<a$,
	and let $u$ be the bilateral solution satisfying $\mathcal{E}_{l}u=v$.
	For every $0<\eta<a-|\chi|$, the two-sided Oseledets estimates imply
	$$
	|u_{m}|\leqslant C_{\eta}(x,v)e^{(|\chi|+\eta)|m|},
	\qquad m\in\mathbb{Z}.
	$$
	Consequently,
	$$
	\|u\|_{a}
	\leqslant C_{\eta}(x,v)
	\sum_{m\in\mathbb{Z}}e^{-(a-|\chi|-\eta)|m|}
	<\infty.
	$$
	Thus $u\in\ker_{X_{a}}\mathbf{L}_{\varepsilon,t}^{(l)}(x)$ and
	$v\in\mathcal{C}_{l}(x)$. In particular,
	$
	\bigoplus_{|\chi|\leqslant a/2}\mathcal{E}_{\chi}^{(l)}(x)
	\subseteq\mathcal{C}_{l}(x).
	$
	We conclude that, for almost every $x$,
	$
	\mathcal{C}_{l}(x)
	=\bigoplus_{|\chi|\leqslant a/2}\mathcal{E}_{\chi}^{(l)}(x).
	$
	The exponents on $\mathcal{C}_{l}$ have absolute value at most
	$\chi_{2}<a/2$. The same argument with $|\chi|<a$ shows that every
	remaining exponent has absolute value at least $a$.
	
	By \cite[Theorem~1.4]{GJYZ26}, the spectral gaps at the boundaries of
	the slow exponents give a dominated splitting
	$
	\mathbb{C}^{2l}
	=E_{l}^{u}(x)\oplus E_{l}^{c}(x)\oplus E_{l}^{s}(x).
	$
	The pairing of the Lyapunov exponents gives
	$$
	\dim E_{l}^{u}=\dim E_{l}^{s}=l-N,
	\qquad \dim E_{l}^{c}=2N.
	$$
	Its middle bundle agrees almost everywhere with $\mathcal{C}_{l}$.
	Since both bundles are continuous, they agree for every $x$:
	$$
	E_{l}^{c}(x)
	=\mathcal{C}_{l}(x)
	=\operatorname{Ran}\big(\mathcal{E}_{l}U_{\varepsilon,t}^{(l)}(x)\big).
	$$
	 Persistence and uniqueness of dominated
	splittings with these dimensions give the asserted joint continuity
	in $(x,\varepsilon,t)$, with $\varepsilon\neq0$
	\cite{AJS14}.
	
	Finally, put $E=E_{*}+t$ and let
	$
	0\leqslant\widehat L_{1}^{(l)}(E)
	\leqslant\cdots\leqslant\widehat L_{l}^{(l)}(E)
	$
	be the nonnegative Lyapunov exponents of the one-step companion
	cocycle. The preceding identification and the $q$-step frame relation
	give
	$$
	\widehat L_{j}^{(l)}(E)
	=\frac{1}{q}L_{N-j+1}(q\alpha,C_{\varepsilon,t}^{(l)}),
	\qquad 1\leqslant j\leqslant N.
	$$
	Moreover, if $l>N$, then
	$
	\widehat L_{j}^{(l)}(E)\geqslant a,
	\qquad N+1\leqslant j\leqslant l.
	$
	The convergence of the slow cocycles and continuity of Lyapunov
	exponents for analytic one-frequency cocycles \cite{AJS14} therefore
	give
	$$
	\widehat L_{j}^{(l)}(E)\longrightarrow
	\widehat L_{j}(E)
	=\frac{1}{q}L_{N-j+1}(q\alpha,C_{\varepsilon,t}^{(\infty)}),
	\qquad 1\leqslant j\leqslant N.
	$$
	The remaining limiting dual exponents, whenever present, are at least
	$a$. Quantitative global theory \cite[Corollary~1.1]{GJYZ26} gives
	$$
	\omega_{\varepsilon,W}(E)
	=\#\{1\leqslant j\leqslant N:\widehat L_{j}(E)=0\}.
	$$
	At $\varepsilon=0$ the operator has range $N$, so the same identification follows directly from the $q$-step frame relation and \cite[Corollary~1.1]{GJYZ26}.
	Since $C_{\varepsilon,t}^{(\infty)}$ preserves $J_{2N}$, its Lyapunov
	exponents occur in opposite pairs. Hence
	$$
	\#\{j:L_{j}(q\alpha,C_{\varepsilon,t}^{(\infty)})=0\}
	=2\#\{1\leqslant j\leqslant N:\widehat L_{j}(E)=0\}
	=2\omega_{\varepsilon,W}(E),
	$$
	which proves \eqref{eq:accD}.
\end{proof}

\subsubsection{Rotation number and IDS}

Let $\mathcal N_0$ be the IDS of $v_0(x)=2\lambda\cos(2\pi Nx)$.
We use the rotation number from Section~\ref{sec:rotation-prelim}, with continuous real lifts for the slow cocycles chosen to agree for all truncations at $(\varepsilon,t)=(0,0)$.

\begin{lemma}\label{lem:localids}
By a suitable continuous real lift, for sufficiently small $\varepsilon$ and $t$, one has\footnote{Equation~\eqref{ids:local-absolute} uses the special form of $v_0$; \eqref{ids:local} holds more generally.}
\begin{equation}\label{ids:local-absolute}
1-\mathcal N_{\varepsilon,W}(E_*+t)
=\frac1q\rho(q\alpha,C_{\varepsilon,t}^{(\infty)}),
\end{equation}
\begin{equation}\label{ids:local}
\mathcal{N}_{\varepsilon,W}(E_{*}+t)-\mathcal{N}_{0}(E_{*})
=-\frac{1}{q}\big(
\rho(q\alpha,C_{\varepsilon,t}^{(\infty)})
-\rho(q\alpha,C_{0,0}^{(\infty)})
\big).
\end{equation}
\end{lemma}

\begin{proof}
	Fix a finite truncation of actual degree $l$ and take $0<|\varepsilon|<\varepsilon_0$. Put $v=v_0+\varepsilon W^{(l)}$ and $E=E_*+t$, and use the explicit choice of $T_v$ in Section~\ref{sec:aubry-duality}.

	Using the normalization in \eqref{eq:Jld}, define
	$$
	\Phi_l=\lambda^{-1/2}T_v^{-1}\mathcal E_lU_{\varepsilon,t}^{(l)}=(\Phi_{l,1},\Phi_{l,2}).
	$$
	By Lemma~\ref{lem:par}, $\Phi_l$ is a frame of $T_v^{-1}E_l^c$. The normalization of the slow frame and its transfer relation give
	\begin{equation}\label{idsframe}
	\begin{aligned}
	\Phi_l^*J_{2l}\Phi_l
	&=\lambda^{-1}(\mathcal E_lU_{\varepsilon,t}^{(l)})^*J_v(\mathcal E_lU_{\varepsilon,t}^{(l)})=J_{2N},\\
	(\widehat{\mathcal A}_E^v)^{(q)}(x)\Phi_l(x)
	&=\Phi_l(x+q\alpha)C_{\varepsilon,t}^{(l)}(x).
	\end{aligned}
	\end{equation}

	By Lemma~\ref{lem:par}, the bundles $E_l^s,E_l^c,E_l^u$ are jointly continuous in $(x,\varepsilon,t)$ on $\mathbb T\times\{0<|\varepsilon|<\varepsilon_0,\ |t|<t_0\}$. Each connected component is homotopy equivalent to $\mathbb T$, so its complex vector bundles are trivial \cite{Forstneric2017}. Choose continuous periodic frames $X_l^s$ and $X_l^u$ of $T_v^{-1}E_l^s$ and $T_v^{-1}E_l^u$, respectively.

	Fix the parameters and choose $\chi_2<\gamma<a$, where $\chi_2$ is as in the proof of Lemma~\ref{lem:par}. The bounds on the fast exponents, together with the uniform subadditive bound, give $K_l>0$ such that, for every $n\geqslant0$,
	$$
	\begin{aligned}
	\big\|(\widehat{\mathcal A}_E^v)^{(n)}(x)v_s\big\|
	&\leqslant K_le^{-\gamma n}\|v_s\|,
	&&v_s\in T_v^{-1}E_l^s(x),\\
	\big\|(\widehat{\mathcal A}_E^v)^{(-n)}(x)v_u\big\|
	&\leqslant K_le^{-\gamma n}\|v_u\|,
	&&v_u\in T_v^{-1}E_l^u(x).
	\end{aligned}
	$$
	The two-sided estimate \eqref{eq:slow} and the frame relation give
	$$
	\big\|(\widehat{\mathcal A}_E^v)^{(\pm n)}(x)v_c\big\|
	\leqslant K_le^{\chi_2 n}\|v_c\|,\qquad v_c\in T_v^{-1}E_l^c(x),
	$$
	after increasing $K_l$. These estimates are uniform in $x$ for fixed $l,\varepsilon,t$. Since $\widehat{\mathcal A}_E^v$ preserves $J_{2l}$ and $\gamma>\chi_2$, applying forward iterates to the preserved pairings shows that $T_v^{-1}E_l^s$ is isotropic and $J_{2l}$-orthogonal to $T_v^{-1}E_l^c$. Backward iterates give the same assertions for $T_v^{-1}E_l^u$.

	The pairing between $T_v^{-1}E_l^s$ and $T_v^{-1}E_l^u$ is nondegenerate. We may therefore choose $X_l^u$ as the symplectic dual frame to $X_l^s$, normalized by $(X_l^s)^*J_{2l}X_l^u=-I_{l-N}$. Define
	$$
	P_l=(X_l^s,\Phi_{l,1},X_l^u,\Phi_{l,2}).
	$$
	The preceding pairing identities and \eqref{idsframe} give $P_l^*J_{2l}P_l=J_{2l}$. Let $D_l(x)$ be determined by
	$$
	(\widehat{\mathcal A}_E^v)^{(q)}(x)X_l^s(x)=X_l^s(x+q\alpha)D_l(x).
	$$
	Write $C_{\varepsilon,t}^{(l)}(x)=(C_{ij}(x))_{i,j=1}^2$. Invariance of the three bundles, preservation of their pairings, and \eqref{idsframe} imply
	\begin{equation}\label{idsblocks}
	P_l(x+q\alpha)^{-1}(\widehat{\mathcal A}_E^v)^{(q)}(x)P_l(x)
	=\begin{pmatrix}
	D_l(x)&0&0&0\\
	0&C_{11}(x)&0&C_{12}(x)\\
	0&0&(D_l(x)^*)^{-1}&0\\
	0&C_{21}(x)&0&C_{22}(x)
	\end{pmatrix}.
	\end{equation}

	With compatible real lifts, Lemma~\ref{lem:idsrot} and Remark~\ref{rem:rotation-properties} give
	$$
	\begin{aligned}
	q\big(1-\mathcal N_{\varepsilon,W^{(l)}}(E_*+t)\big)
	&=\rho\big(q\alpha,(\mathcal A_E^v)^{(q)}\big)=\rho\big(q\alpha,(\widehat{\mathcal A}_E^v)^{(q)}\big).
	\end{aligned}
	$$
	After a fixed coordinate permutation, the right-hand side of \eqref{idsblocks} is $\operatorname{diag}(D_l,(D_l^*)^{-1})\oplus C_{\varepsilon,t}^{(l)}$, preserving $J_{2(l-N)}\oplus J_{2N}$. The first summand preserves a constant Lagrangian subspace, so its rotation number is zero modulo $\mathbb Z$. Applying Remark~\ref{rem:rotation-properties} and Lemma~\ref{lem:rotation-conjugacy} directly to \eqref{idsblocks} gives
	\begin{equation}\label{idsdegree}
	q\big(1-\mathcal N_{\varepsilon,W^{(l)}}(E_*+t)\big)
	-\rho(q\alpha,C_{\varepsilon,t}^{(l)})
	-q\alpha\deg P_l\in\mathbb Z.
	\end{equation}
	By continuity, both $\deg P_l$ and the integer in \eqref{idsdegree} are constant on each connected component. Subtract \eqref{idsdegree} at $(\varepsilon,t)$ and $(\varepsilon',0)$ in the same component and let $\varepsilon'\to0$. Since $C_{0,0}^{(l)}=C_{0,0}^{(\infty)}=A_{*,q}$, Proposition~\ref{prop:slow} and the continuity of the IDS \cite[Lemma~8.1]{CCB26} and rotation number (Remark~\ref{rem:rotation-properties}) give
	\begin{equation}\label{idsfin}
	\mathcal N_{\varepsilon,W^{(l)}}(E_*+t)-\mathcal N_0(E_*)
	=-\frac1q\big(\rho(q\alpha,C_{\varepsilon,t}^{(l)})-\rho(q\alpha,C_{0,0}^{(\infty)})\big).
	\end{equation}
	Continuity extends \eqref{idsfin} to $\varepsilon=0$.

	Finally, pass to the full potential $W$. Uniform convergence of its Fourier truncations gives convergence of the IDS at every energy \cite[Lemma~8.1]{CCB26}. Proposition~\ref{prop:slow} gives uniform convergence of the slow cocycles, and their rotation numbers converge in the real lifts fixed at $(0,0)$ by Remark~\ref{rem:rotation-properties}. Passing to the limit in \eqref{idsfin} proves \eqref{ids:local}.

By Lemma~\ref{lem:idsrot}, choose the reference lift so that $1-\mathcal N_0(E_*)=\rho(\alpha,A_{*,1})$. Since $C_{0,0}^{(\infty)}=A_{*,q}$, the iteration property and \eqref{ids:local} give \eqref{ids:local-absolute}.
\end{proof}

\subsection*{Proof of Theorem~\ref{thm:slow}}
\pushQED{\qed}
Apply Propositions~\ref{prop:slow1} and~\ref{prop:slow} with
$l=\infty$, $q=N$, and $t=E-E_*$. The first gives the rank-$2N$
analytic solution bundle and its global frame; the second normalizes
the symplectic form to $J_{2N}$ and gives the stated translation
relation. Neither construction requires a Diophantine condition:
the reference growth bound \eqref{eq:reference-growth} holds for
every irrational $\alpha$. \qedhere
\popQED

\section{Acceleration Stratification}\label{sec:acceleration-strata}

We apply Proposition~\ref{prop:slow} at the upper spectral edge of the unperturbed operator. The parabolic conjugacy determines the reference frame, and Corollary~\ref{cor:modevar} gives its first-order response to energy and potential.

In the rest paper, we will always assume $\alpha \in DC$.
Fix the parameters of Theorem~\ref{thm:strata}, and set $E_0=\sup\Sigma_0$. 
Fix $0<a<2\pi r$. We write $E=E_0-\varepsilon\xi$.

\subsection{Positive edge state}

\begin{lemma}\label{lem:edge}
 There exist a one-periodic real analytic function $p>0$ and a one-periodic real analytic map $Z:\mathbb{T}\to\mathrm{SL}(2,\mathbb{R})$, with first column $(p(x),p(x-N\alpha))^{\mathsf T}$, such that
\begin{equation}\label{eq:edge}
\lambda\bigl[p(x+N\alpha)+p(x-N\alpha)\bigr]+2\cos(2\pi x)p(x)=E_{0}p(x)
\end{equation}
and
\begin{equation}\label{eq:para}
Z(x+N\alpha)^{-1}S_{\lambda^{-1},E_{0}/\lambda}(x)Z(x)=\begin{pmatrix}1&c\\0&1\end{pmatrix},\qquad
c=-\int_{\mathbb{T}}\frac{\mathrm{d}x}{p(x)p(x+N\alpha)}<0.
\end{equation}
\end{lemma}

\begin{proof}
Since $\lambda>1$, $(N\alpha,S_{\lambda^{-1},E_{0}/\lambda})$ is a subcritical almost Mathieu cocycle. As $1-\mathcal N_0(E_0)=0$ and $N\alpha\in DC$,  \cite{Eliasson92}
and \cite[Theorem~4.1(3)]{AJ10} give
$$
B(x+N\alpha)^{-1}S_{\lambda^{-1},E_{0}/\lambda}(x)B(x)
=\begin{pmatrix}\sigma&c\\0&\sigma\end{pmatrix},\qquad
\sigma\in\{\pm1\},\quad c\in\mathbb{R},
$$
where $B:\mathbb R\to\mathrm{SL}(2,\mathbb R)$ is analytic and
$B(x+1)=\pm B(x)$. By Remark~\ref{rem:rotation-conjugacy},
the degree $k$ of the induced $\mathrm{PSL}(2,\mathbb R)$ map
satisfies $kN\alpha\in\mathbb Z$. Thus $k=0$, so $B$ is
one-periodic. It preserves the usual fibred rotation number,
whose value at $E_0$ is zero, and hence $\sigma=1$.
Writing its first column as
$(\varphi(x),\varphi(x-N\alpha))^{\mathsf T}$, we obtain a
nonzero real analytic solution $\varphi$ of \eqref{eq:edge}.

We now prove positivity. On $L^2(\mathbb T)$,
consider
$$
 (\mathcal Lf)(x)
 =\lambda\bigl(f(x+N\alpha)+f(x-N\alpha)\bigr)
   +2\cos(2\pi x)f(x).
$$
Fourier transformation identifies $\mathcal L$ with the almost
Mathieu operator of coupling $\lambda$ and frequency $N\alpha$
at phase zero. Thus $\mathcal L\leqslant E_0I$, and
$$
\begin{aligned}
 0\geqslant
 \langle|\varphi|,(\mathcal L-E_0)|\varphi|\rangle=2\lambda\int_{\mathbb T}
 \bigl(|\varphi(x)\varphi(x+N\alpha)|
       -\varphi(x)\varphi(x+N\alpha)\bigr)\,dx
 \geqslant0.
\end{aligned}
$$
Equality implies $\mathcal L|\varphi|=E_0|\varphi|$, pointwise
by continuity. Any zero would propagate along the dense
$N\alpha$-orbit, forcing $\varphi\equiv0$. Hence $\varphi$ has
a fixed sign, and $p=|\varphi|>0$ is real analytic and satisfies
\eqref{eq:edge}.
Set
$$
 Q(x)=\begin{pmatrix}
 p(x)&0\\ p(x-N\alpha)&p(x)^{-1}
 \end{pmatrix},
 \qquad b(x)=-\frac{1}{p(x)p(x+N\alpha)}.
$$
A direct computation gives
$$
Q(x+N\alpha)^{-1}S_{\lambda^{-1},E_{0}/\lambda}(x)Q(x)
=\begin{pmatrix}1&b(x)\\0&1\end{pmatrix}.
$$
Let $c=\int_{\mathbb T}b(x)\,dx<0$. Since $N\alpha$ is
Diophantine, there is a one-periodic real analytic solution of
$ f(x+N\alpha)-f(x)=b(x)-c.$
Then
$ Z(x)=Q(x)\begin{pmatrix}1&f(x)\\0&1\end{pmatrix}$
has the prescribed first column and satisfies \eqref{eq:para}.
\end{proof}

Let $p>0$ and $c<0$ be as in Lemma~\ref{lem:edge}. For
$W(x)=\sum_{n\in\mathbb Z}w_ne^{2\pi inx}$, define
$$
 A=\int_{\mathbb T}p(x)^2\,dx,\qquad
 I_n=\int_{\mathbb T}p(x)p(x+n\alpha)\,dx,\qquad n\in\mathbb Z,
$$
and
$ M_W=\sum_{n\in\mathbb Z}w_nI_n\mathsf S^n,$
where $(\mathsf Sb)_j=b_{j-1}$ with indices modulo $N$.
Since $I_{-n}=I_n\in\mathbb R$ and
$w_{-n}=\overline{w_n}$, $M_W$ is Hermitian.

\subsection{Monotonicity in the multidimensional Moser--P\"oschel unfolding}
Fix a compact interval $J_0\subset\mathbb R$.
\begin{proposition}\label{prop:center}
For every $W\in C_{r}^{\omega}(\mathbb{T},\mathbb{R})$ and each integer $p_{0}\geqslant1$, there are $h,\varepsilon_{0},K_{p_{0}}>0$ with the following property. The rank-$2N$ slow bundle at $E_{0}-\varepsilon\xi$ admits an analytic frame $U_{\varepsilon,\xi}$ near $\{0\}\times J_{0}$. For $0<\varepsilon<\varepsilon_{0}$, its $N$-step cocycle $C_{\varepsilon,\xi}$ admits the normal form
\begin{equation}\label{mpred}
B_{0}(x+N\alpha,\xi)^{-1}C_{\varepsilon,\xi}(x)B_{0}(x,\xi)
=C_{0}(\xi)e^{F_{0}(x,\xi)},
\end{equation}
where $B_{0},C_{0},F_{0}$ are analytic, $B_{0}$ is homotopic to the identity, and, for real variables,
$
C_{\varepsilon,\xi},B_{0},C_{0}\in\mathrm{HSp}(J_{2N}),
F_{0}\in\mathfrak{hsp}(J_{2N}).
$
Let
$$
K_{0}(\xi)=\lambda^{-1}(M_{W}+A\xi I_{N}),\qquad
L(\xi)=\begin{pmatrix}0&cI_{N}\\K_{0}(\xi)&0\end{pmatrix}.
$$
Then
\begin{equation}\label{mpstart}
A_{0}:=\varepsilon^{-\frac{1}{2}}\log C_{0}
=L+O_{C^{1}(J_{0})}(\sqrt{\varepsilon}),
\sup_{\xi\in J_{0}}\big(\|F_{0}(\cdot,\xi)\|_{h}
+\|\partial_{\xi}F_{0}(\cdot,\xi)\|_{h}\big)
\leqslant K_{p_{0}}\varepsilon^{p_{0}+\frac{1}{2}}
\end{equation}
With the continuous real lifts tending to zero as $\varepsilon\to0$, one has
\begin{equation}\label{idsrho}
1-\mathcal{N}_{\varepsilon,W}(E_{0}-\varepsilon\xi)
=\frac{1}{N}\rho(N\alpha,C_{\varepsilon,\xi})
=\frac{1}{N}\rho(N\alpha,C_{0}e^{F_{0}}).
\end{equation}
For fixed $p_{0}$, the constants and estimates are locally uniform in $W$.
\end{proposition}

\begin{proof}
We first construct the cocycle before averaging. In the coordinates \eqref{eq:initD}, define $H_{*}(x)$ by applying $Z(x+(N-1-j)\alpha)$ to each coordinate pair $(u_{N-1-j},u_{-1-j})$, $0\leqslant j<N$. By \eqref{equ:direct} and \eqref{eq:para},
$$
H_{*}(x)^{*}J_{2N}H_{*}(x)=J_{2N},\qquad
A_{*,N}(x)H_{*}(x)=H_{*}(x+N\alpha)C_{*},
\qquad
C_{*}=\begin{pmatrix}I_{N}&cI_{N}\\0&I_{N}\end{pmatrix}.
$$
Since $Z_{11}=p>0$, $Z$ is homotopic to a constant. The same holds for $H_{*}$, so its degree is zero.

Apply Proposition~\ref{prop:slow} with $q=N$ and $E_{*}=E_{0}$. We take the frame $U_{\varepsilon,-\varepsilon\xi}^{(\infty)}H_{*}$ and denote it and its associated cocycle by $U_{\varepsilon,\xi}$ and $C_{\varepsilon,\xi}$. Then
\begin{equation}\label{eq:icc}
S^{N}U_{\varepsilon,\xi}(x)
=U_{\varepsilon,\xi}(x+N\alpha)C_{\varepsilon,\xi}(x),\qquad
C_{\varepsilon,\xi}=C_{*}e^{F_{\varepsilon,\xi}}.
\end{equation}
These families are analytic near $\{0\}\times J_{0}$ in $(\varepsilon,\xi)$. On a fixed complex neighborhood of $J_{0}$ and a phase strip of width $h_{1}>0$,
$
\|F_{\varepsilon,\xi}\|_{h_{1}}\leqslant K|\varepsilon|,
F_{0,\xi}=0.
$
The bounds are locally uniform in $W$. For real variables, $C_{\varepsilon,\xi}$ preserves $J_{2N}$ and $F_{\varepsilon,\xi}\in\mathfrak{hsp}(J_{2N})$.

\begin{lemma}\label{lem:average}
For $\varepsilon\neq0$, set
$K_{\varepsilon}(\xi)=\varepsilon^{-1}
\langle(F_{\varepsilon,\xi})_{21}\rangle$.
It extends analytically to $\varepsilon=0$, with value
$K_{0}(\xi)=\lambda^{-1}(M_{W}+A\xi I_{N})$, and
$
\langle(F_{\varepsilon,\xi})_{21}\rangle
=\varepsilon K_{0}(\xi)+O_{C^{1}(J_{0})}(\varepsilon^{2}).
$
For sufficiently small real $\varepsilon$, $K_{\varepsilon}$
is Hermitian and
$\partial_{\xi}K_{\varepsilon}\geqslant
\frac{A}{2\lambda}I_{N}$ on $J_{0}$.
Every ordered eigenvalue of $K_{\varepsilon}$ is therefore
strictly increasing in $\xi$.
\end{lemma}

\begin{proof}
For $z=(b,0)^{\mathsf T}$, set $u_m(x)=(U_{0,\xi}(x)z)_m$ and extend
$b$ periodically. Since $Z(y)e_1=(p(y),p(y-N\alpha))^{\mathsf T}$,
the construction of $H_*$ and the edge equation give
$u_m(x)=p(x+m\alpha)b_{N-1-m}$ for all $m\in\mathbb Z$.
Since $H_{*}$ is independent of $\varepsilon,t$ and preserves $J_{2N}$, multiplying \eqref{varid} on the left by $H_{*}(x)^{*}$ and on the right by $H_{*}(x)$ gives the same identity for this frame. 
We apply Corollary~\ref{cor:modevar} at $(0,0)$ with $q=N$
and direction $(1,-\xi)$, using $C_{*}z=z$.
Since
$$
z^{*}J_{2N}
\partial_{\varepsilon}F_{\varepsilon,\xi}\big|_{\varepsilon=0}z
=-b^{*}
\big(\partial_{\varepsilon}F_{\varepsilon,\xi}
\big|_{\varepsilon=0}\big)_{21}b,
$$
we obtain
$$
\lambda b^{*}
\big\langle\big(\partial_{\varepsilon}F_{\varepsilon,\xi}
\big|_{\varepsilon=0}\big)_{21}\big\rangle b
=\sum_{m=0}^{N-1}\int_{\mathbb{T}}\overline{u_{m}(x)}
\bigg(\xi u_{m}(x)+\sum_{n\in\mathbb{Z}}w_{n}u_{m+n}(x)\bigg)
\,\mathrm{d}x=b^{*}(A\xi I_{N}+M_{W})b.
$$
 Polarization and Taylor
expansion give $\big\langle(F_{\varepsilon,\xi})_{21}\big\rangle
 =\varepsilon K_0(\xi)+O_{C^1(J_0)}(\varepsilon^2)$, with the $C^1$ bound by Cauchy's
estimate. 
Since $F_{\varepsilon,\xi}\in\mathfrak{hsp}(J_{2N})$ for
real variables, its lower-left block is Hermitian.
Thus
$\partial_{\xi}K_{\varepsilon}
=\frac{A}{\lambda}I_{N}+O(|\varepsilon|)$
gives the stated lower bound.
For $\xi<\xi'$ in $J_{0}$,
$
K_{\varepsilon}(\xi')-K_{\varepsilon}(\xi)
\geqslant\frac{A}{2\lambda}(\xi'-\xi)I_{N},
$
so the min--max principle gives strict increase of every
ordered eigenvalue.
\end{proof}
Fix $p_{0}\geqslant1$. We apply finite-order averaging to the cocycle \eqref{eq:icc}, always solving the homological equation at $C_{*}$.

At order $j$, let $F_{j}$ be the term homogeneous of degree $j$ in $\varepsilon$ in the current remainder. Solve
$$
C_{*}^{-1}X_{j}(x+N\alpha,\xi)C_{*}-X_{j}(x,\xi)
=F_{j}(x,\xi)-\langle F_{j}\rangle,\qquad
\langle X_{j}\rangle=0.
$$
At this stage the current constant is $C_{*}+O(\varepsilon)$. Replacing it by $C_{*}$ in the homological equation introduces an error of order $\varepsilon^{j+1}$. Conjugation by $e^{X_{j}}$ therefore removes the nonconstant term of order $j$, and its average is absorbed into the constant. Repeating this procedure $p_{0}$ times gives, on a smaller strip of width $h>0$,
\begin{equation}\label{aprepe}
Q_{p_{0}}(x+N\alpha,\xi)^{-1}C_{\varepsilon,\xi}(x)Q_{p_{0}}(x,\xi)
=\widetilde C_{p_{0}}(\xi)e^{\widetilde F_{p_{0}}(x,\xi)},
\end{equation}
where
$$
\begin{aligned}
\sup_{\xi\in J_{0}}\big(\|Q_{p_{0}}(\cdot,\xi)-I_{2N}\|_{h}
+\|\partial_{\xi}Q_{p_{0}}(\cdot,\xi)\|_{h}\big)
\leqslant K_{p_{0}}\varepsilon,
\widetilde C_{p_{0}}
&=C_{*}e^{\langle F_{\varepsilon,\xi}\rangle}
+O_{C^{1}(J_{0})}(\varepsilon^{2}),\\
\sup_{\xi\in J_{0}}\big(\|\widetilde F_{p_{0}}(\cdot,\xi)\|_{h}
+\|\partial_{\xi}\widetilde F_{p_{0}}(\cdot,\xi)\|_{h}\big)
&\leqslant K_{p_{0}}\varepsilon^{p_{0}+1}.
\end{aligned}
$$
The expansion of $\widetilde C_{p_{0}}$ follows from the first averaging step, since all later constant corrections have order at least $\varepsilon^{2}$. The derivative bounds follow by carrying out the construction on a slightly larger complex neighborhood of $J_{0}$ and applying Cauchy's estimate.

Each averaging transformation preserves $J_{2N}$. Hence $Q_{p_{0}}$ and $\widetilde C_{p_{0}}$ preserve $J_{2N}$, and $\widetilde F_{p_{0}}\in\mathfrak{hsp}(J_{2N})$ for real variables.

By Lemma~\ref{lem:average},
$$
\widetilde C_{p_{0}}
=\begin{pmatrix}
I_{N}+O(\varepsilon)&cI_{N}+O(\varepsilon)\\
\varepsilon K_{0}+O(\varepsilon^{2})&I_{N}+O(\varepsilon)
\end{pmatrix},
$$
where all errors are estimated in $C^{1}(J_{0})$. Set
$$
D_{\varepsilon}
=\operatorname{diag}\big(\varepsilon^{-\frac{1}{4}}I_{N},
\varepsilon^{\frac{1}{4}}I_{N}\big),\qquad
B_{0}=Q_{p_{0}}D_{\varepsilon},
$$
and define
$$
C_{0}=D_{\varepsilon}^{-1}\widetilde C_{p_{0}}D_{\varepsilon},
\qquad
F_{0}=D_{\varepsilon}^{-1}\widetilde F_{p_{0}}D_{\varepsilon}.
$$
Then \eqref{aprepe} gives \eqref{mpred}, and direct block multiplication yields
\begin{equation*}\label{mpscale}
C_{0}=I_{2N}+\sqrt{\varepsilon}\,L
+O_{C^{1}(J_{0})}(\varepsilon).
\end{equation*}
After decreasing $\varepsilon_{0}$, the local logarithm is defined and satisfies
$$
\log C_{0}
=\sqrt{\varepsilon}\,L+O_{C^{1}(J_{0})}(\varepsilon).
$$
This proves the estimate for $A_{0}$ in \eqref{mpstart}. Moreover,
$
\|D_{\varepsilon}^{-1}\|\|D_{\varepsilon}\|
=\varepsilon^{-\frac{1}{2}},
$
so the estimate for $\widetilde F_{p_{0}}$ gives the stated bound for $F_{0}$.

The matrix $D_{\varepsilon}$ preserves $J_{2N}$ and is joined to the identity by $D_{\varepsilon}^{s}$, $0\leqslant s\leqslant1$. Since $Q_{p_{0}}$ is a product of exponentials in the same group, $B_{0}$ is homotopic to the identity. The matrices $C_{0}$ and $F_{0}$ have the asserted group and Lie algebra properties.

It remains to verify the IDS formula. Since $\mathcal{N}_{0}(E_{0})=1$ and $H_{*}$ has degree zero, Lemma~\ref{lem:localids} gives
$$
1-\mathcal{N}_{\varepsilon,W}(E_{0}-\varepsilon\xi)
=\frac{1}{N}\rho(N\alpha,C_{\varepsilon,\xi}),
$$
with the lift normalized by $\rho(N\alpha,C_{*})=0$. The conjugation by $B_{0}$ is homotopic to the identity, so \eqref{mpred} gives
$$
\rho(N\alpha,C_{\varepsilon,\xi})
=\rho(N\alpha,C_{0}e^{F_{0}})
$$
for the compatible real lifts. These lifts tend to zero as $\varepsilon\to0$, since $C_{\varepsilon,\xi}\to C_{*}$ and $C_{0}e^{F_{0}}\to I_{2N}$. This proves \eqref{idsrho}.

\end{proof}

We first determine when the eigenvalues of $M_W$
are distinct. Let
\begin{equation*}\label{eq:generic}
 \mathcal G=\bigl\{W\in C_r^\omega(\mathbb T,\mathbb R):
                    M_W\text{ has simple spectrum}\bigr\}.
\end{equation*}

\begin{lemma}\label{lem:generic}\label{prop:generic}
The complement of $\mathcal G$ is a finite union of proper closed
real linear hyperplanes in $C_r^\omega(\mathbb T,\mathbb R)$.
It is Gauss-null, and $\mathcal G$ is open, dense, and prevalent.
\end{lemma}

\begin{proof}
Let $\zeta=e^{2\pi i/N}$ and let $\mathcal{F}$ be the unitary
Fourier matrix with entries
$\mathcal{F}_{i,j}=N^{-\frac{1}{2}}\zeta^{(N-1-i)j}$,
$0\leqslant i,j<N$.
Then
\begin{equation}\label{eq:nu}
\mathcal{F}^{*}M_{W}\mathcal{F}
=\operatorname{diag}\big(\widetilde\nu_{0}(W),\ldots,
\widetilde\nu_{N-1}(W)\big),\qquad
\widetilde\nu_{a}(W)=\sum_{n\in\mathbb{Z}}w_{n}I_{n}\zeta^{na}.
\end{equation}
These functionals are real-valued and real linear, since
$w_{-n}=\overline{w_{n}}$ and $I_{-n}=I_{n}\in\mathbb{R}$.
The bounds $|I_{n}|\leqslant A$ and
$|w_{n}|\leqslant\|W\|_{r}e^{-2\pi r|n|}$ give
$|\widetilde\nu_{a}(W)|\leqslant K\|W\|_{r}$.
Thus they are continuous.

For $0\leqslant a<b<N$, put
$L_{a,b}=\widetilde\nu_{a}-\widetilde\nu_{b}$.
For $W_{s,t}(x)=s\cos(2\pi x)+t\sin(2\pi x)$,
 $$
L_{a,b}(W_{s,t})
=I_{1}\operatorname{Re}\big((s-it)(\zeta^{a}-\zeta^{b})\big).
$$
Taking $s=\operatorname{Re}(\zeta^{a}-\zeta^{b})$ and
$t=\operatorname{Im}(\zeta^{a}-\zeta^{b})$ gives
$L_{a,b}(W_{s,t})=I_{1}|\zeta^{a}-\zeta^{b}|^{2}>0$.
Thus $L_{a,b}$ is nonzero.
Consequently,
$$
C_{r}^{\omega}(\mathbb{T},\mathbb{R})\setminus\mathcal{G}
=\bigcup_{0\leqslant a<b<N}\ker L_{a,b}
$$
is a finite union of proper closed real linear hyperplanes.
In particular, $\mathcal{G}$ is open and dense.

Let $\mathfrak{m}$ be a nondegenerate Gaussian probability
measure on $C_{r}^{\omega}(\mathbb{T},\mathbb{R})$.
Each $L_{a,b}$ has a Gaussian distribution with positive
variance, so $\mathfrak{m}(\ker L_{a,b})=0$.
This proves the Gauss-null assertion. In particular,
$\mathcal{G}$ is prevalent.
\end{proof}

For $W\in\mathcal G$, order the eigenvalues of $M_W$ as
$\nu_1(W)>\cdots>\nu_N(W)$ and put
\begin{equation}\label{eq:transition-values}
 \xi_j(W)=-\nu_j(W)/A.
\end{equation}
Fix $W_{*}\in\mathcal{G}$ and a nondegenerate compact interval
$I\subset\mathbb{R}\setminus
\{\xi_{1}(W_{*}),\ldots,\xi_{N}(W_{*})\}$.
There are an open ball $\mathcal{V}\subset\mathcal{G}$ about
$W_{*}$ and constants $R_{*},\Delta_{*},\sigma_{*}>0$ such that
$$
\|W\|_{r}\leqslant R_{*},\qquad
\nu_{j}(W)-\nu_{j+1}(W)\geqslant\Delta_{*},\quad1\leqslant j< N,\quad
|\nu_{j}(W)+A\xi|\geqslant\sigma_{*},\quad
1\leqslant j\leqslant N,
$$
for every $W\in\mathcal{V}$ and $\xi\in I$.
Hence $k=\#\{j:\xi_{j}(W)<\xi\}$ is constant on
$\mathcal{V}\times I$.

For $W\in\mathcal V$, the eigenvalues $(\nu_j(W)+A\xi)/\lambda$ of
$K_0$ move with the common slope $A/\lambda$. Their distinct zeros
$\xi_j(W)$ are the leading thresholds.

Fix $p_0\geqslant1$ and let $C_0,A_0,L$ be given by
Proposition~\ref{prop:center}, with $J_0$ containing $I$ in its interior.

\begin{lemma}\label{lem:mp-transversality}
There is $c_*>0$ such that, for sufficiently small $\varepsilon>0$,
the following hold uniformly for $W\in\mathcal V$.
Write $\ell_a$ for the eigenvalues of $L$. The matrix $C_0$ admits a
diagonalization $V_0^{-1}C_0V_0=\operatorname{diag}(e^{2\pi i\mu_a})$
with $V_0^{\pm1}=O_{C^1(I)}(1)$ and
\begin{equation}\label{mpfreq}
 \max_a\|2\pi i\mu_a-\sqrt\varepsilon\,\ell_a\|_{C^1(I)}
 \leqslant K_{p_0}\varepsilon .
\end{equation}
All the multipliers are simple. Exactly $2k$ lie on the unit circle,
in $k$ pairs whose phases move in opposite directions as $\xi$ increases.
One also has
$$
\|C_{0}-I_{2N}\|_{C^1(I)}\leqslant K_{p_{0}}\sqrt{\varepsilon},\qquad \inf_{\xi\in I,\,a\neq b}|e^{2\pi i\mu_{a}(\xi)}-e^{2\pi i\mu_{b}(\xi)}|\geqslant c_{*}\sqrt{\varepsilon}.
$$
Moreover, the scalar divisors
$$
g_{ab}(n,\xi)=e^{2\pi inN\alpha}-e^{2\pi i(\mu_{a}(\xi)-\mu_{b}(\xi))},\qquad n\in\mathbb{Z},\quad a\neq b,
$$
satisfy
\begin{equation*}\label{eq:mp-transversality}
	\big|\partial_{\xi}\operatorname{Re}g_{ab}(n,\xi)\big|\geqslant c_*\varepsilon ,\qquad \xi\in I,\quad n\in\mathbb{Z},\quad a\neq b.
\end{equation*}
\end{lemma}

\begin{proof}
After permuting the Fourier columns in \eqref{eq:nu}
according to $\nu_{1}(W)>\cdots>\nu_{N}(W)$, the same
permutation works for every $W\in\mathcal{V}$.
Applying this basis to both coordinate blocks reduces $L$
to the blocks
$\begin{pmatrix}0&c\\s_{j}&0\end{pmatrix}$,
where $s_{j}=\lambda^{-1}(\nu_{j}(W)+A\xi)$.
 Since $c<0$, their roots are
\begin{equation}\label{aroots}
 \ell_j=-\ell_{N+j}=
 \begin{cases}
   i r_j,\quad r_j=\sqrt{-c\,s_j}, & j\leqslant k,\\
   g_j,\quad g_j=\sqrt{c\,s_j},    & j>k .
 \end{cases}
\end{equation}
The bounds on $\mathcal V\times I$ give
\begin{equation}\label{equ:gap}
 \min_a|\ell_a|\geqslant c_*,\qquad
 \min_{a\ne b}|\ell_a-\ell_b|\geqslant c_*,\qquad
 \max_a\|\ell_a\|_{C^1(I)}\leqslant K
\end{equation}
for some $c_*>0$. The block eigenvectors $(c,\ell_a)^{\mathsf T}$
give a uniformly $C^1$-bounded diagonalizer of $L$ with bounded
$C^1$ inverse. By \eqref{mpstart}, the Riesz projections of $A_0$
differ from those of $L$ by $O_{C^1(I)}(\sqrt\varepsilon)$.
Applying them to these eigenvectors gives $V_0$ and simple eigenvalues
$\beta_a=\ell_a+O_{C^1(I)}(\sqrt\varepsilon)$. Taking
$2\pi i\mu_a=\sqrt\varepsilon\,\beta_a$ proves \eqref{mpfreq}.

The spectrum of $A_0\in\mathfrak{hsp}(J_{2N})$ is invariant under
$\beta\mapsto-\overline\beta$. Uniqueness of the branch near each
$\pm i r_j$ keeps it on the imaginary axis, while the branches near
$\pm g_j$ have real parts bounded away from zero. The exponential is
injective on the small disk containing $\sqrt\varepsilon\,\operatorname{spec}A_0$.
Thus exactly $2k$ simple multipliers are on the circle. Since
$r_j'=-cA/(2\lambda r_j)$ is uniformly positive, \eqref{mpfreq}
gives opposite signs for the two phase derivatives in each elliptic pair.

The bounds for $C_0-I_{2N}$ and the separation of its eigenvalues
follow from \eqref{mpstart}, \eqref{equ:gap}, and \eqref{mpfreq}.
We next verify transversality for $g_{ab}$.
By \eqref{aroots},
$$
\ell_{a}'=\frac{cA}{2\lambda\ell_{a}},\qquad \ell_{a}'-\ell_{b}'=-\frac{cA}{2\lambda}\frac{\ell_{a}-\ell_{b}}{\ell_{a}\ell_{b}}.
$$
Thus $|\ell_{a}'-\ell_{b}'|\geqslant c_{*}$. Put $d_{ab}=\beta_{a}-\beta_{b}$,
$$
\partial_{\xi}\operatorname{Re}g_{ab}(n,\xi)=-\sqrt{\varepsilon}\operatorname{Re}\big(e^{\sqrt{\varepsilon}d_{ab}}d_{ab}'\big).
$$
If the two branches are not both elliptic, then $|\operatorname{Re}(\ell_{a}'-\ell_{b}')|\geqslant c_{*}$: for two hyperbolic branches the difference is real, and for a mixed pair its real part is the derivative of the hyperbolic branch. By \eqref{mpfreq},
$$
\big|\partial_{\xi}\operatorname{Re}g_{ab}(n,\xi)\big|\geqslant c_{*}\sqrt{\varepsilon}-K_{p_{0}}\varepsilon\geqslant c_{1}\sqrt{\varepsilon}.
$$
If both branches are elliptic, then $d_{ab}=i\omega_{ab}$ with $\omega_{ab}$ real. The $C^{1}$ eigenvalue estimate gives $|\omega_{ab}|,|\omega_{ab}'|\geqslant c_{*}$ and $|\omega_{ab}|\leqslant K$. Hence
$$
\big|\partial_{\xi}\operatorname{Re}g_{ab}(n,\xi)\big|=\sqrt{\varepsilon}\big|\sin(\sqrt{\varepsilon}\omega_{ab})\omega_{ab}'\big|\geqslant c_{1}\varepsilon.
$$
\end{proof}

\subsection{Positive-measure reducibility}

For a compact parameter set $\Lambda\subset\mathbb R$, set
\begin{equation*}\label{eq:parameter-norm}
 \|F\|_{h,\Lambda}
 =\sup_{\xi\in\Lambda}\|F(\cdot,\xi)\|_h
 +\sup_{\substack{\xi,\xi'\in\Lambda\\\xi\ne\xi'}}
 \frac{\|F(\cdot,\xi)-F(\cdot,\xi')\|_h}{|\xi-\xi'|}.
\end{equation*}
For functions independent of $x$, we omit $h$.

\begin{proposition}\label{unfold}
There are $h,\varepsilon_{0},K>0$ such that, for $W\in\mathcal{V}$ and $0<\varepsilon<\varepsilon_{0}$, there is a compact set $T_{\varepsilon,W}\subset I$ satisfying
\begin{equation}\label{unfmeas}
\operatorname{Leb}(I\setminus T_{\varepsilon,W})\leqslant K\varepsilon^{\frac52}.
\end{equation}
On $T_{\varepsilon,W}$, there is a $B\in C^\omega_{\frac{h}{2}}(\mathbb T,\mathrm{HSp}(J_{2N}))$, homotopic to the identity, such that
\begin{equation}\label{unfred}
B(x+N\alpha,\xi)^{-1}C_{\varepsilon,\xi}(x)B(x,\xi)=C_{\infty}(\xi).
\end{equation}
The constant $C_{\infty}$ belongs to $\mathrm{HSp}(J_{2N})$ and satisfies
\begin{equation}\label{unflog}
\bigg\|\varepsilon^{-\frac12}\log C_{\infty}-\begin{pmatrix}0&cI_{N}\\\lambda^{-1}(M_{W}+A\xi I_{N})&0\end{pmatrix}\bigg\|_{T_{\varepsilon,W}}\leqslant K\sqrt{\varepsilon}.
\end{equation}
It has exactly $2k$ simple eigenvalues on the unit circle; the remaining eigenvalues are off the unit circle.
\end{proposition}

\begin{proof}

We use the following discrete simple-spectrum version of the
parameter iteration in \cite{HeYou06}, with explicit dependence on
the spectral separation and transversality constants.

Let $C:I\to\mathrm{HSp}(J_{2N})$ and $V:I\to\mathrm{GL}(2N,\mathbb{C})$ be $C^{1}$ families such that
$$
V^{-1}CV=\operatorname{diag}(e^{2\pi i\mu_{1}},\ldots,e^{2\pi i\mu_{2N}}),\qquad \mu_{a}\in C^{1}(I,\mathbb{C}).
$$

Assume
\begin{equation}\label{data}
	\|C-I_{2N}\|_{I}\leqslant\frac{1}{2},\qquad \|V\|_{I}+\|V^{-1}\|_{I}+\max_{a}\|\mu_{a}\|_{I}\leqslant M_{0},
\end{equation}
where $M_{0}\geqslant1$ is fixed, and
\begin{equation}\label{sep}
	\inf_{\xi\in I}|e^{2\pi i\mu_{a}(\xi)}-e^{2\pi i\mu_{b}(\xi)}|
	\geqslant\rho>0,\qquad a\ne b.
\end{equation}
For this family, define $g_{ab}$ by the same formula as above.
Assume also that
\begin{equation}\label{nd1}
	\big|\partial_{\xi}\operatorname{Re}g_{ab}(n,\xi)\big|\geqslant\chi>0,\qquad \xi\in I,\quad n\in\mathbb{Z},\quad a\neq b.
\end{equation}

Fix $\tau_{1}>\max\{\tau,1\}$ and choose
\begin{equation*}\label{hypar}
	\nu\geqslant1,\qquad b_{\mathrm{K}}>16N^{2}\nu+1,\qquad
	s_{\mathrm{K}}=2\big(4N^{2}(\tau_{1}+1)+2N\tau\big)+1,\qquad
	M=M_{0}+\rho^{-1}+\chi^{-1}.
\end{equation*}

\begin{theorem}\cite{HeYou06}\label{kam}
	Suppose \eqref{data}, \eqref{sep}, and \eqref{nd1} hold. Let $F$ be analytic on $|\operatorname{Im}x|<h$, Lipschitz in $\xi\in I$, and $\mathfrak{hsp}(J_{2N})$-valued on real phases. Put $\eta=\|F\|_{h,I}$. There are $c_{\mathrm{K}}>0$, $0<\theta_{0}<1$, and $K>0$, depending only on $2N,\gamma,\tau,\tau_{1},\nu,b_{\mathrm{K}},h,I$ and $M_{0}$, such that
	\begin{equation}\label{small}
		\eta\leqslant\min\big\{c_{\mathrm{K}}h^{b_{\mathrm{K}}(2s_{\mathrm{K}}+1)}M^{-b_{\mathrm{K}}(32N^{4}+1)},\theta_{0}^{b_{\mathrm{K}}}e^{-\frac{6}{5}}\big\}
	\end{equation}
	implies the existence of a compact set $T\subset I$ and matrices $B_{\infty},C_{\infty}$ satisfying
	\begin{equation*}
		B_{\infty}(x+N\alpha,\xi)^{-1}C(\xi)e^{F(x,\xi)}B_{\infty}(x,\xi)=C_{\infty}(\xi),\qquad \xi\in T.
	\end{equation*}
	Here $B_{\infty}$ is analytic on $|\operatorname{Im}x|<h/2$, both matrices preserve $J_{2N}$ on real phases, and
	\begin{equation}\label{kamb}
		\operatorname{Leb}(I\setminus T)\leqslant C_{\mathrm{K}}|I|(10\eta)^{\nu/b_{\mathrm{K}}},\qquad C_{\mathrm{K}}\leqslant K(1+\chi^{-1}).
	\end{equation}
	Moreover,
	\begin{equation}\label{const}
		\|C_{\infty}-C\|_{T}+\|\log C_{\infty}-\log C\|_{T}\leqslant K\eta,\quad
\sup_{\xi\in T}\|B_{\infty}(\cdot,\xi)-I_{2N}\|_{h/2}<\frac{1}{2}.
	\end{equation}
\end{theorem}

\begin{remark}
We follow the KAM scheme of \cite{HeYou06}, with the notation
and iteration parameters adapted to the discrete setting.
Although their result state the result for $gl(m,\C)$ linear system, the Hermitian-symplectic assertion follows from the standard
structure-preserving KAM iteration, as the homological equation preserve the   Hermitian-symplectic structure, consult \cite{Chavaudret13,Eliasson92} and   Lemma \ref{cl:KAM} for details.
\end{remark}

Apply Proposition~\ref{prop:center} with $p_0\geqslant1$,
to be chosen below.
We apply Theorem~\ref{kam} with $C=C_{0}$, $F=F_{0}$, and $V=V_{0}$.
By \eqref{mpstart}, $\eta:=\|F_0\|_{h,I}\leqslant K_{p_0}\varepsilon^{p_0+1/2}$.
Lemma~\ref{lem:mp-transversality} gives \eqref{data}, \eqref{sep}, and \eqref{nd1}, with $M_{0}$ independent of $\varepsilon$, $\rho=c_{*}\sqrt{\varepsilon}$, and $\chi=c_{*}\varepsilon$.

Choose $p_{0}$ so that
$p_{0}+\frac{1}{2}>Q:=b_{\mathrm{K}}(32N^{4}+1)$.
Then
$
\big(p_{0}+\frac{1}{2}\big)\frac{\nu}{b_{\mathrm{K}}}
>\nu(32N^{4}+1)>\frac{7}{2}.
$
The preceding bounds give
$$
M\leqslant K\varepsilon^{-1},\qquad
c_{\mathrm{K}}h^{b_{\mathrm{K}}(2s_{\mathrm{K}}+1)}
M^{-b_{\mathrm{K}}(32N^{4}+1)}
\geqslant c_{1}\varepsilon^{Q}.
$$
Thus \eqref{mpstart} implies \eqref{small} for sufficiently
small $\varepsilon$, uniformly for $W\in\mathcal{V}$.
By \eqref{kamb},
$$
\operatorname{Leb}(I\setminus T_{\varepsilon,W})
\leqslant K\varepsilon^{-1}
\big(10K_{p_{0}}\varepsilon^{p_{0}+\frac{1}{2}}\big)^{\nu/b_{\mathrm{K}}}
\leqslant K\varepsilon^{\frac{5}{2}}.
$$
This proves \eqref{unfmeas}.
Set $B=Q_{p_0}D_\varepsilon B_\infty$. Then \eqref{unfred} holds.
The matrices $B_\infty$ and $Q_{p_0}$ are close to the identity, while $D_\varepsilon$ is positive diagonal. Their logarithms lie in $\mathfrak{hsp}(J_{2N})$, so their product $B=Q_{p_0}D_\varepsilon B_\infty$ is homotopic to the identity.

Finally, \eqref{const} and \eqref{mpstart} give
$$
A_{\infty}=\varepsilon^{-\frac{1}{2}}\log C_{\infty}\in\mathfrak{hsp}(J_{2N}),\qquad \|A_{\infty}-L\|_{T_{\varepsilon,W}}\leqslant K_{p_{0}}\sqrt{\varepsilon}+K\varepsilon^{-\frac{1}{2}}\eta\leqslant K\sqrt{\varepsilon}.
$$
This proves \eqref{unflog}. 
The simple-spectrum perturbation argument used for $A_{0}$
applies to $A_{\infty}$.
Its spectrum is invariant under
$z\mapsto-\overline{z}$.
Each neighborhood of $\pm ir_{j}$ contains one eigenvalue,
which is therefore imaginary, while the remaining eigenvalues
have real parts bounded away from zero.
Exponentiation gives exactly $2k$ simple eigenvalues of
$C_{\infty}$ on the unit circle; all others are off it.

\end{proof}

\subsection{Proofs of the main results}
Fix $W_*\in\mathcal G$ and an interval $I$ as in Theorem~\ref{thm:landscape}. Choose a compact interval $I^+$ containing $I$ in its interior and avoiding every $\xi_j(W_*)$. Apply Proposition~\ref{unfold} on $I^+$, with a common neighborhood $\mathcal V$ of $W_*$. The integer $k=\#\{j:\xi_j(W)<\xi\}$ is constant for $W\in\mathcal V$ and $\xi\in I^+$. For every sufficiently small $\varepsilon>0$, the proposition gives a compact set $T_{\varepsilon,W}^+\subset I^+$ such that
\begin{equation}\label{measure}
\operatorname{Leb}(I^+\setminus T_{\varepsilon,W}^+)\leqslant K\varepsilon^{5/2},
\end{equation}
with \eqref{unfred}--\eqref{unflog} on this set. The bounds are uniform in $W$, and the conjugacy has degree zero.

\subsection*{Proof of Theorem~\ref{thm:landscape}}
\pushQED{\qed}
Set $\kappa=N^{-1}\sqrt{-cA/\lambda}>0$.
(i)
Assume $I\subset(\xi_1(W_*),\infty)$, so that $k\geqslant1$, and let $T_{\varepsilon,W}=I\cap T_{\varepsilon,W}^{+}$. Equation~\eqref{measure} gives the measure estimate. We suppress $W$ where no confusion is possible.
	
	On $T_{\varepsilon,W}^{+}$, Proposition~\ref{unfold} gives exactly $2k$ zero slow   exponents. Lemma~\ref{lem:par} therefore gives $\omega_{\varepsilon,W}(E_{0}-\varepsilon\xi)=k$. Since the Schr\"odinger cocycle has zero acceleration outside the spectrum, $k\geqslant1$ also gives $E_0-\varepsilon\xi\in\Sigma_{\varepsilon,W}$. This proves the spectral assertion and \eqref{eq:intro-acceleration-count}.
	
	For $1\leqslant j\leqslant k$, write $r_{j}=N\kappa\sqrt{\xi-\xi_{j}(W)}$. Write the elliptic eigenvalues of $C_\infty$ as $e^{i\vartheta_j}$ and $e^{i\vartheta_{N+j}}$, choosing the phases near zero. By \eqref{unflog},
	\begin{equation}\label{angles}
		\vartheta_{j}=\sqrt{\varepsilon}r_{j}+O(\varepsilon),\qquad \vartheta_{N+j}=-\sqrt{\varepsilon}r_{j}+O(\varepsilon),\qquad 1\leqslant j\leqslant k.
	\end{equation}
	Set $k_{\varepsilon,j}=\vartheta_{j}/N$, for small $\varepsilon$,
	$$
	k_{\varepsilon,j}=\kappa\sqrt{\varepsilon}\sqrt{\xi-\xi_{j}(W)}+O(\varepsilon)>0,\qquad
	\kappa^{-2}k_{\varepsilon,j}^{2}=\varepsilon(\xi-\xi_{j}(W))+O(\varepsilon^{\frac{3}{2}}).
	$$
	This proves \eqref{eq:intro-dispersion}.
	
	Choose $z_j\neq0$ with $C_\infty z_j=e^{iNk_{\varepsilon,j}}z_j$ and set
	$
	Y_j(x)=U_{\varepsilon,\xi}(x)B(x,\xi)z_j.
	$
	Equations~\eqref{eq:icc} and \eqref{unfred} give $S^{N}Y_{j}(x)=e^{iNk_{\varepsilon,j}}Y_{j}(x+N\alpha)$. Thus $F_{j,s}(x)=(Y_{j}(x))_{s}$ gives
	$$
	(Y_{j}(x))_{Nm+s}=e^{iNmk_{\varepsilon,j}}F_{j,s}(x+mN\alpha),
	$$
	which is the required nonzero analytic block Bloch solution \eqref{eq:intro-block-Bloch}.

(ii)
We return to the general interval $I$, allowing $k=0$. The elliptic eigenlines for $\vartheta_j$
and $\vartheta_{N+j}$ have negative and positive $iJ_{2N}$-sign,
respectively, by $c<0$ and \eqref{unflog}; the hyperbolic part preserves
a constant Lagrangian subspace. Lemma~\ref{lem:unitary-rotation} and
Remark~\ref{rem:rotation-properties} therefore give
	$
	\rho(N\alpha,C_{\infty})=\frac{1}{2\pi}\sum_{j=1}^{k}(\vartheta_{j}-\vartheta_{N+j})
	$ in the real lift tending to zero as $\varepsilon\downarrow0$.
Equations~\eqref{idsrho} and \eqref{angles} now give
	$$
	1-\mathcal{N}_{\varepsilon,W}(E_{0}-\varepsilon\xi)=\sqrt{\varepsilon}\Phi_{W}(\xi)+O(\varepsilon),\qquad \xi\in T_{\varepsilon,W}^{+},
	$$
	where
	$
	\Phi_{W}(\xi)=\frac{\kappa}{\pi}\sum_{j=1}^{N}\sqrt{(\xi-\xi_{j}(W))_{+}}.
	$
	The derivative of $\Phi_{W}$ is uniformly bounded on $I^{+}$.

	Set $d_{\varepsilon}=2K\varepsilon^{\frac{5}{2}}$, with $K$ bounding the right-hand side of \eqref{measure}. Since $I\subset\operatorname{int}I^{+}$, for small $\varepsilon$ both intervals $[\xi-d_{\varepsilon},\xi]$ and $[\xi,\xi+d_{\varepsilon}]$ lie in $I^{+}$ and meet $T_{\varepsilon,W}^{+}$. Choose points $\xi_{-}\leqslant\xi\leqslant\xi_{+}$ in that set with $|\xi_{\pm}-\xi|\leqslant d_{\varepsilon}$. Monotonicity of the IDS gives
	$$
		1-\mathcal{N}_{\varepsilon,W}(E_{0}-\varepsilon\xi)
		\leqslant1-\mathcal{N}_{\varepsilon,W}(E_{0}-\varepsilon\xi_{+})\leqslant\sqrt{\varepsilon}\Phi_{W}(\xi)+K\sqrt{\varepsilon}d_{\varepsilon}+K\varepsilon.
	$$
	The lower bound follows using $\xi_{-}$. Since $\sqrt{\varepsilon}d_{\varepsilon}=O(\varepsilon^{3})$, we obtain
	\begin{equation}\label{idsprof}
		\pi\big[1-\mathcal{N}_{\varepsilon,W}(E_{0}-\varepsilon\xi)\big]=\kappa\sqrt{\varepsilon}\sum_{j=1}^{N}\sqrt{(\xi-\xi_{j}(W))_{+}}+O(\varepsilon),\qquad \xi\in I.
	\end{equation}

We next establish the Lyapunov exponent asymptotics on the whole interval $I$.  Use the normal form $C_{0}e^{F_{0}}$ and the diagonalizing matrix $V_{0}$ given by Proposition~\ref{prop:center} and Lemma~\ref{lem:mp-transversality}, respectively. By \eqref{mpstart} and the bounds preceding \eqref{mpfreq}, we have
$$\widehat{C}(x)=V_{0}^{-1}C_{0}e^{F_{0}(x)}V_{0}=\mathcal{D}(I_{2N}+\mathcal{R}(x)),
	\mathcal{D}=\operatorname{diag}\big(e^{2\pi i\mu_{1}},\ldots,e^{2\pi i\mu_{2N}}\big),$$
where $
	\mathcal{R}=e^{V_{0}^{-1}F_{0}V_{0}}-I_{2N}, \|\mathcal{R}\|_{\infty}\leqslant K\eta\leqslant K\varepsilon^{\frac{3}{2}}.$ 
Let $b_{a}=\operatorname{Re}(2\pi i\mu_{a})$. \eqref{aroots}, \eqref{equ:gap} and \eqref{mpfreq} give $2k$ zero values and separate all distinct $b_{a}$ by at least $c_{*}\sqrt{\varepsilon}$. For any distinct value $b\in\{b_a\}_{a=1}^{2N}$, set $t_{\pm}=b\pm\varepsilon$. For small $\varepsilon$, the diagonal system $e^{-t_{\pm}}\mathcal{D}$ is uniformly hyperbolic, with contraction constant $1$ and rate $\varepsilon$. By quantitative stability of uniform hyperbolicity \cite[Theorem~5.2]{BFP21}, the perturbed system $e^{-t_{\pm}}\widehat C$ is uniformly hyperbolic with the same stable dimension, since
$
K\eta\frac{1+e^{-\varepsilon}}{1-e^{-\varepsilon}}\leqslant K\sqrt{\varepsilon}<1.
$
 Let $\Lambda_{1}\geqslant\cdots\geqslant\Lambda_{2N}$ be the Lyapunov exponents of $\widehat{C}$. The weighted systems $e^{-t_{\pm}}\widehat C$ have exponents $\Lambda_{j}-t_{\pm}$, so preservation of the stable dimension gives
$$
\#\{j:\Lambda_{j}<t_{\pm}\}=\#\{a:b_{a}<t_{\pm}\}.
$$ Since $b$ is the only value among the $b_{a}$ in $(b-\varepsilon,b+\varepsilon)$, subtraction gives
	$
	\#\{j:b-\varepsilon<\Lambda_{j}<b+\varepsilon\}=\#\{a:b_{a}=b\}.
	$
The nonnegative slow exponents per lattice site are
$$
\gamma_j(W,\xi)=\frac1N L_{N-j+1}(N\alpha,C_{\varepsilon,\xi})
=\frac1N\Lambda_{N-j+1},\qquad 1\leqslant j\leqslant N.
$$
Combining the preceding count with \eqref{aroots} and \eqref{mpfreq} gives
\begin{equation}\label{legamma}
	\bigg|\gamma_{j}(W,\xi)-\kappa\sqrt{\varepsilon}\sqrt{\max\{\xi_{j}(W)-\xi,0\}}\bigg|\leqslant K\varepsilon,\qquad 1\leqslant j\leqslant N.
\end{equation}
	By Lemma~\ref{lem:par}, these are the first $N$ nonnegative limiting dual exponents, and all remaining exponents are at least $a$. Fix $y_{0}>0$ with $2\pi y_{0}<a$. For small $\varepsilon$, \eqref{legamma} gives $\gamma_{j}<2\pi y_{0}$ for every $j$. The multiplicative Jensen formula \cite[Theorem~1.2]{GJYZ26} yields, with $E=E_{0}-\varepsilon\xi$,
	\begin{equation}\label{lejensen}
		L_{\varepsilon,W}(E;y)=L_{\varepsilon,W}(E;y_{0})-2\pi Ny_{0}+\sum_{j=1}^{N}\max\{\gamma_{j}(W,\xi),2\pi|y|\},\qquad 2\pi|y|<a.
	\end{equation}
	At the fixed height $y_{0}$, the unperturbed formula $L_{0}(E_{0};y)=\log\lambda+2\pi N|y|$ is positive and affine. The cocycle is therefore uniformly hyperbolic, and smooth dependence of the Lyapunov exponent gives \cite{Avila15}
	$
	L_{\varepsilon,W}(E_{0}-\varepsilon\xi;y_{0})=\log\lambda+2\pi Ny_{0}+O(\varepsilon)
	$
	uniformly in $W\in\mathcal{V}$ and $\xi\in I$. Substitute this and \eqref{legamma} into \eqref{lejensen}, using $|\max\{u,t\}-\max\{v,t\}|\leqslant|u-v|$, to obtain
	\begin{equation}\label{leprof}
			L_{\varepsilon,W}(E_{0}-\varepsilon\xi;y)=\log\lambda+\sum_{j=1}^{N}\max\big\{\kappa\sqrt{\varepsilon}\sqrt{(\xi_{j}(W)-\xi)_{+}},2\pi|y|\big\}+O(\varepsilon), \quad\xi\in I, 2\pi|y|<a.
	\end{equation}
	Combining \eqref{idsprof} and \eqref{leprof}, and setting $y=\sqrt{\varepsilon}\eta$, proves \eqref{eq:intro-uniform-landscape}. \qedhere
\popQED

\subsection*{Proof of Theorem~\ref{thm:strata}}
\pushQED{\qed}
Set $\xi_{N+1}(W_{*})=+\infty$ and choose compact intervals of positive length
	$
	I_{k}\Subset(\xi_{k}(W_{*}),\xi_{k+1}(W_{*})),\quad 1\leqslant k\leqslant N.
	$
	Apply Theorem~\ref{thm:landscape} to these finitely many intervals. Take the minimum of their neighborhood radii and thresholds and the maximum of their constants. Decrease the threshold so that
	$K\varepsilon^{\frac{5}{2}}\leqslant\frac{1}{2}\min_{k}|I_{k}|.$
	The map $\xi\mapsto E_0-\varepsilon\xi$ multiplies Lebesgue
measure by $\varepsilon$. Thus \eqref{eq:intro-acceleration-count} gives
$$
 \operatorname{Leb}\{E\in\Sigma_{\varepsilon,W}:\omega_{\varepsilon,W}(E)=k\}
 \geqslant\varepsilon(|I_k|-K\varepsilon^{5/2})
 \geqslant\tfrac12|I_k|\varepsilon.
$$
This is \eqref{eq:positive-measure} with $c_k=|I_k|/2$.
Lemma~\ref{lem:generic} gives the assertions about $\mathcal G$. \qedhere
\popQED

\subsection*{Proof of Corollary~\ref{cor:non-typeI-open}}
\pushQED{\qed}
Apply Theorem~\ref{thm:strata} with $N=2$ and $\lambda=2$, and fix $W_*\in\mathcal G$. Let $\mathcal V$, $\varepsilon_0$, and $c_2$ be the resulting neighborhood and constants. Choose $\delta>0$ so that the ball of radius $\delta$ about $W_*$ lies in $\mathcal V$. Fix $0<\varepsilon<\varepsilon_0$ and set
$$
\mathcal U=\{v\in C_r^\omega(\mathbb T,\mathbb R):\|v-v_{\varepsilon,W_*}\|_r<\varepsilon\delta\}.
$$
Every $v\in\mathcal U$ has the form $v_{\varepsilon,W}$ with $\|W-W_*\|_r<\delta$. Thus \eqref{eq:positive-measure} gives the asserted bound with $c=c_2\varepsilon$.
On the energy windows used above, \eqref{leprof} gives
$L_v(E;0)\geqslant\frac12\log2>0$ after decreasing $\varepsilon_0$.
Hence these acceleration-two energies are not Type~I.
Thus $\mathcal{U}$ contains no potential whose operator at
the fixed frequency $\alpha$ is Type~I.
Potentials with the Type~I property throughout a neighborhood
uniformly over all irrational frequencies also have this
property at the fixed $\alpha$.
Hence both assertions of
Conjecture~\ref{conj:whole-spectrum-typeI} are false.\qedhere
\popQED

\section{Collapsed gaps}
The potential $v_{0}(x)=2\lambda\cos(2\pi Nx)$ has the IDS of the almost Mathieu operator at frequency $N\alpha$. By \cite{AYZ23}, every label in $N\mathbb{Z}\setminus\{0\}$ corresponds to an open gap. The gap-labelling theorem excludes open gaps with labels outside $N\mathbb{Z}$. Hence, for each $m_{0}\in\mathbb{Z}\setminus N\mathbb{Z}$, there is a unique energy $E_{*}$ satisfying
$\mathcal{N}_{0}(E_{*})=\{m_{0}\alpha\}.$
Fix such a label and its energy, and fix $W\in C_{r}^{\omega}(\mathbb{T},\mathbb{R})$. We study the slow cocycle near $E_{*}$, writing $E=E_{*}+t$.
In this section we take $q=1$ in Proposition~\ref{prop:slow}.
The single-channel transfer matrix still
advances $N$ sites, whereas the slow cocycle advances one site.

\begin{lemma}\label{cor:ell}
There exist a diagonal unitary matrix $D_{*}$, an analytic reference frame $U_{*}$, and constants $s_{*},h,C>0$ such that the following hold for $0<s\leqslant s_{*}$:
The spaces $\mathcal{M}_{\varepsilon,t}^{(\infty)}(x)$ have dimension $2N$ and admit jointly analytic frames satisfying
\begin{equation*}\label{cl:slow}
SU_{\varepsilon,t}(x)=U_{\varepsilon,t}(x+\alpha)C_{\varepsilon,t}(x),\qquad C_{\varepsilon,t}=D_{*}e^{F_{\varepsilon,t}},
\end{equation*}
where for real $x,\varepsilon,t$,
$C_{\varepsilon,t}\in\mathrm{HSp}(\widehat J_{2N}), F_{\varepsilon,t}\in\mathfrak{hsp}(\widehat J_{2N}).
$
Moreover,
$$
\|U_{\varepsilon,t}-U_{*}\|_{h,s}+\|F_{\varepsilon,t}\|_{h,s}\leqslant Cs.
$$
For real $x$ and sufficiently small real $\varepsilon,t$,
\begin{equation}\label{cl:ids}
\mathcal{N}_{\varepsilon,W}(E_{*}+t)-\{m_{0}\alpha\}=-\rho(\alpha,C_{\varepsilon,t})+\rho(\alpha,D_{*}).
\end{equation}
\end{lemma}

\begin{proof}
The almost Mathieu cocycle $(N\alpha,S_{\lambda^{-1},E_{*}/\lambda})$ has coupling $\lambda^{-1}<1$. By Lemma \ref{lem:idsrot}, its rotation number is $1-\{m_{0}\alpha\}$. The Diophantine condition on $\alpha$ gives
$$
\|1-\{m_{0}\alpha\}-kN\alpha\|_{\mathbb{T}}=\|(m_{0}+Nk)\alpha\|_{\mathbb{T}}\geqslant c(1+|k|)^{-\tau},\qquad k\in\mathbb{Z}.
$$
Thus 
\cite[Theorem~4.1]{AJ10} gives a one-periodic analytic $\mathrm{SL}(2,\mathbb R)$ conjugacy to a constant rotation through angle $2\pi\varrho_*$.
Applying Lemmas~\ref{lem:unitary-rotation} and~\ref{lem:rotation-conjugacy} to the real conjugacy gives
$$
2\varrho_{*}\equiv1-\{m_{0}\alpha\}\pmod{N\alpha\mathbb{Z}+\mathbb{Z}}.
$$
By the Cayley isomorphism $\mathrm{SL}(2,\mathbb R)\simeq\mathrm{SU}(1,1)$,
we obtain $Z_*\in C^\omega(\mathbb T,\mathrm{SL}(2,\mathbb C))$ satisfying
\begin{equation}\label{refred}
\begin{aligned}
Z_{*}(x+N\alpha)^{-1}S_{\lambda^{-1},E_{*}/\lambda}(x)Z_{*}(x)
&=\operatorname{diag}(e^{-2\pi i\varrho_{*}},e^{2\pi i\varrho_{*}}),\\
Z_{*}(\bar{x})^{*}J_{2}Z_{*}(x)&=\widehat J_{2}.
\end{aligned}
\end{equation}
To distinguish the one-step phases, set
$\theta_{\sigma,j}=(\sigma\varrho_*+j)/N$ for $\sigma=\pm$ and $0\le j<N$.
Since $e^{2\pi iN\theta_{\sigma,j}}=e^{2\pi i\sigma\varrho_*}$,
\eqref{refred} defines scalar solutions by
$$
 \binom{u_{\sigma,j,m}(x)}{u_{\sigma,j,m-N}(x)}
 =e^{2\pi im\theta_{\sigma,j}}Z_*(x+m\alpha)\mathbf e_\sigma,
 \qquad m\in\mathbb Z,
$$
where $\mathbf e_-=(1,0)^{\mathsf T}$ and $\mathbf e_+=(0,1)^{\mathsf T}$.
They satisfy $Su_{\sigma,j}(x)=e^{2\pi i\theta_{\sigma,j}}u_{\sigma,j}(x+\alpha)$,
and hence
$
 A_{*,1}(x)\mathcal E u_{\sigma,j}(x)
 =e^{2\pi i\theta_{\sigma,j}}\mathcal E u_{\sigma,j}(x+\alpha).
$
Moreover, \eqref{refred} gives
$$
 (\mathcal E u_{\sigma,j}(\bar x))^*J_{2N}\mathcal E u_{\sigma',k}(x)
 =(\widehat J_2)_{\sigma\sigma'}
 \sum_{s=0}^{N-1}e^{2\pi is(\theta_{\sigma',k}-\theta_{\sigma,j})}
 =N(\widehat J_2)_{\sigma\sigma'}\delta_{jk}.
$$
Thus $H_*(x)=N^{-1/2}\mathcal E\bigl(u_{-,0}(x),\ldots,u_{-,N-1}(x)
\,\bigm|\,u_{+,0}(x),\ldots,u_{+,N-1}(x)\bigr)$ satisfies
$$
 H_*(x+\alpha)^{-1}A_{*,1}(x)H_*(x)
 =D_*:=\operatorname{diag}_{\sigma,j}(e^{2\pi i\theta_{\sigma,j}}),\qquad
 H_*(\bar x)^*J_{2N}H_*(x)=\widehat J_{2N}.
$$
The solution frame $U_*=\widetilde U_*H_*$ has columns
$N^{-1/2}u_{\sigma,j}$.

\begin{lemma}\label{lem:seed}
The phases $\theta_{\sigma,j}$ are non-resonant. That is, there exists $c>0$ such that
\begin{equation}\label{cl:seed}
\|\theta_{\sigma,j}-\theta_{\sigma',k}-n\alpha\|_{\mathbb{T}}\geqslant c(1+|n|)^{-\tau},\qquad (\sigma,j)\neq(\sigma',k),\quad n\in\mathbb{Z}.
\end{equation}
\end{lemma}
\begin{proof}
For $\sigma'=\sigma$ and $n\ne0$,
$N(\theta_{\sigma,j}-\theta_{\sigma,k}-n\alpha)=j-k-Nn\alpha$,
so the Diophantine condition and $\|x\|_{\mathbb T}\geqslant N^{-1}\|Nx\|_{\mathbb T}$ give \eqref{cl:seed}.
When $n=0$, $j\ne k$ gives
$\|\theta_{\sigma,j}-\theta_{\sigma,k}\|_{\mathbb T}
=\|(j-k)/N\|_{\mathbb T}\geqslant N^{-1}$.
For $\sigma'=-\sigma$, write
$2\varrho_*\equiv(dN-m_0)\alpha\pmod{\mathbb Z}$, $d\in\mathbb Z$.
Then $N(\theta_{\sigma,j}-\theta_{-\sigma,k}-n\alpha)
\equiv-\sigma[m_0+N(\sigma n-d)]\alpha\pmod{\mathbb Z}$.
Since $m_0+N(\sigma n-d)\equiv m_0\not\equiv0\pmod N$,
the same estimate gives \eqref{cl:seed}.
\end{proof}
Applying Proposition~\ref{prop:slow}, we take the frame $U_{\varepsilon,t}^{(\infty)}H_{*}$ and denote it and its associated cocycle by $U_{\varepsilon,t}$ and $C_{\varepsilon,t}$, respectively. They satisfy the stated estimates. Since $H_{*}$ is independent of $(\varepsilon,t)$, \eqref{ids:local} and \eqref{eq:rotation-conjugacy-difference} give \eqref{cl:ids}. \qedhere

\end{proof}

\subsection{Preliminary averaging}

Let
$
G=\partial_{t}F_{\varepsilon,t}\big|_{(0,0)},F_{W}=\partial_{\varepsilon}F_{\varepsilon,t}\big|_{(0,0)}.$
Thus
\begin{equation}\label{cl:first-variation-expansion}
C_{\varepsilon,t}=D_{*}\exp\big(tG+\varepsilon F_{W}+O((|t|+|\varepsilon|)^{2})\big).
\end{equation}
We first show that the rotation number can be prescribed by changing the energy.

\begin{lemma}\label{cl:prep}
For $0<h_{1}<h$, $0<s\leqslant s_{*}$, there are analytic matrices $B,D_W,R$ satisfying
\begin{equation}\label{cl:nf}
B(\varepsilon,E,x+\alpha)^{-1}C_{\varepsilon,E-E_{*}}(x)B(\varepsilon,E,x)=D_W(\varepsilon,E)e^{R(\varepsilon,E,x)}.
\end{equation}
Here $B(0,E_{*},x)=I_{2N}$, $B\in\mathrm{HSp}(\widehat J_{2N})$ for real variables, $D_W$ is diagonal unitary. Moreover,
$$
D_W(\varepsilon,E)=D_{*}\exp\big((E-E_{*})\langle G\rangle_{\mathrm{diag}}+\varepsilon\langle F_{W}\rangle_{\mathrm{diag}}+O((|\varepsilon|+|E-E_{*}|)^{2})\big),
$$
and
$$
\|B-I_{2N}\|_{h_{1},s}\leqslant Ks,\quad
\|R\|_{h_{1},s}\leqslant Ks^{12}.
$$
There is a unique local real analytic function $E_W^{(0)}(\varepsilon,\vartheta)$ near $(0,\{m_{0}\alpha\})$ such that
$$
E_W^{(0)}(0,\{m_{0}\alpha\})=E_{*},\qquad \rho\big(\alpha,D_W(\varepsilon,E_W^{(0)}(\varepsilon,\vartheta))\big)-\rho(\alpha,D_{*})=\{m_{0}\alpha\}-\vartheta,
$$
and $|E_W^{(0)}(\varepsilon,\vartheta)-E_{*}|\leqslant K(|\varepsilon|+|\vartheta-\{m_{0}\alpha\}|)$.
\end{lemma}

\begin{proof}
At step $k=1,\ldots,11$, write the current cocycle as
$$
 D_*e^A\exp\bigl(F_k(x)+O_{h',s}(s^{k+1})\bigr),
$$
where $A=O_s(s)$ is diagonal, with $A\in\mathfrak{hsp}(\widehat J_{2N})$
for real variables, and $F_k$ is homogeneous of degree $k$ in
$(\varepsilon,t)$. Initially, $A=0$ and $F_1=tG+\varepsilon F_W$.
Solve 
$$
 D_*^{-1}X_k(x+\alpha)D_*-X_k(x)
 =F_k(x)-\langle F_k\rangle_{\mathrm{diag}},\qquad
 \langle X_k\rangle_{\mathrm{diag}}=0.
$$
The Diophantine condition on $\alpha$ and \eqref{cl:seed} give a unique
solution, homogeneous of degree $k$, with
$$
 \|X_k\|_{h'',s}\leqslant C(h'-h'')^{-(\tau+1)}\|F_k\|_{h',s},
 \qquad 0<h''<h'\leqslant h.
$$
Since $A=O_s(s)$ and $X_k=O_{h'',s}(s^k)$,
$$
 e^{-A}D_*^{-1}X_k(x+\alpha)D_*e^A
 =D_*^{-1}X_k(x+\alpha)D_*+O_{h'',s}(s^{k+1}).
$$
Thus we solve the same homological equation with $D_*$ at each
step; the resulting error is of order $s^{k+1}$. The Baker--Campbell--Hausdorff
formula shows that conjugation by $e^{X_k}$ gives
$$
 D_*e^{A+\langle F_k\rangle_{\mathrm{diag}}}
 \exp\bigl(O_{h'',s}(s^{k+1})\bigr).
$$
Accumulating a total strip loss of $h-h_1$ over these eleven steps gives
\eqref{cl:nf} and the stated bounds. The first step gives the stated
expansion of $D_W$. Each $X_k$ lies in
$\mathfrak{hsp}(\widehat J_{2N})$ for real variables, by uniqueness of the
homological solution; the diagonal averages are skew-Hermitian.
Hence $B$ preserves $\widehat J_{2N}$ and $D_W$ is diagonal unitary.

By \eqref{modepos},
$\operatorname{tr}(\widehat J_{2N}\langle G\rangle_{\mathrm{diag}})>0$.
Lemma~\ref{lem:unitary-rotation} and the expansion of $D_W$ therefore give
$$
 \left.\partial_E\rho(\alpha,D_W(\varepsilon,E))\right|_{(0,E_*)}
 =-\frac1{2\pi}\operatorname{tr}
   (\widehat J_{2N}\langle G\rangle_{\mathrm{diag}})<0.
$$
The rotation number of $D_W$ is real analytic for the chosen local lift,
so the implicit function theorem gives $E_W^{(0)}$ and the stated energy estimate. \qedhere

\end{proof}

\subsection{Transversality and typical perturbations}
We index the modes by $\mathbf{a},\mathbf{b},\mathbf{c},\ldots$, writing $u_{\mathbf{a},m}=u_{\sigma,j,m}$
for $\mathbf{a}=(\sigma,j)$. Since $(u_{\mathbf{a},0},u_{\mathbf{a},-N})^{\mathsf T}=Z_*\mathbf e_\sigma$
is independent of $j$, define $\psi_\sigma$ by
$$
\frac{u_{\mathbf{a},0}(x)}{\|u_{\mathbf{a},0}\|_{L^2}}
=\sum_{n\in\mathbb Z}\psi_\sigma(n)e^{2\pi inx},\qquad
\sum_{n\in\mathbb Z}|\psi_\sigma(n)|^2=1.
$$
Substituting $u_{\mathbf{a},m}(x)=e^{2\pi i m\theta_{\sigma,j}}u_{\mathbf{a},0}(x+m\alpha)$
into the dual equation at site zero and comparing Fourier coefficients gives
$$
\psi_\sigma(n+1)+\psi_\sigma(n-1)
+2\lambda\cos\bigl(2\pi N(\theta_{\sigma,j}+n\alpha)\bigr)\psi_\sigma(n)
=E_*\psi_\sigma(n).
$$
Thus $\psi_\sigma$ is a normalized eigenfunction of $H_{v_0,\alpha,\theta_{\sigma,j}}$.
Define\footnote{Which comes from the famous Hellmann–Feynman formula.}
\begin{equation*}\label{cl:energy-velocities}
v_{\sigma,j}=\sum_{n\in\mathbb Z}|\psi_\sigma(n)|^2v_0'(\theta_{\sigma,j}+n\alpha).
\end{equation*}
Applying \eqref{refred}
to the initial pair and using Parseval's identity gives
$$
\begin{aligned}
-\sigma
&=-2\operatorname{Im}\int_{\mathbb T}\overline{u_{\mathbf{a},-N}(x)}u_{\mathbf{a},0}(x)\,\mathrm dx\\
&=-2\operatorname{Im}\bigg(e^{2\pi iN\theta_{\sigma,j}}
\int_{\mathbb T}\overline{u_{\mathbf{a},0}(x)}u_{\mathbf{a},0}(x+N\alpha)\,\mathrm dx\bigg)\\
&=-2\|u_{\mathbf{a},0}\|_{L^2}^2\sum_{n\in\mathbb Z}|\psi_\sigma(n)|^2
\sin\bigl(2\pi N(\theta_{\sigma,j}+n\alpha)\bigr)=\frac{\|u_{\mathbf{a},0}\|_{L^2}^2}{2\pi\lambda N}v_{\sigma,j}.
\end{aligned}
$$
Hence
\begin{equation}\label{cl:profile-norm}
\operatorname{sgn}v_{\sigma,j}=-\sigma,\qquad
\frac{\|u_{\mathbf{a},0}\|_{L^2}^2}{N}=\frac{2\pi\lambda}{|v_{\sigma,j}|}.
\end{equation}

Abbreviate $E_W(\varepsilon)=E_W^{(0)}(\varepsilon,\{m_{0}\alpha\})$ and choose local real analytic phases by
$$
D_W(\varepsilon,E)_{\mathbf{a}\mathbf{a}}=e^{2\pi i\theta_{\mathbf{a},W}(\varepsilon,E)},\qquad \theta_{\mathbf{a},W}(0,E_{*})=\theta_{\mathbf{a}}.
$$
Write
$$
\partial_{\varepsilon}E_W=\frac{\mathrm{d}}{\mathrm{d}\varepsilon}E_W(\varepsilon)\bigg|_{\varepsilon=0},\ 
\partial_{E}\theta_{\mathbf{a}}=\frac{\partial}{\partial E}\theta_{\mathbf{a},W}(0,E)\bigg|_{E=E_{*}},\ 
\partial_{\varepsilon}\theta_{\mathbf{a},W}=\frac{\mathrm{d}}{\mathrm{d}\varepsilon}\theta_{\mathbf{a},W}(\varepsilon,E_W(\varepsilon))\bigg|_{\varepsilon=0}.
$$
By Lemma~\ref{cl:prep} and the definition of $E_W$, we have
\begin{equation}\label{cl:fixed-rotation}
\rho\bigl(\alpha,D_W(\varepsilon,E_W(\varepsilon))\bigr)=\rho(\alpha,D_*).
\end{equation}
Thus the last derivative is taken along the energy curve $E_W(\varepsilon)$, on which the rotation number of the prepared constant remains fixed. Define
$$
w_{\mathbf{a}}(W)=\sum_{n\in\mathbb{Z}}|\psi_{\sigma}(n)|^{2}W(\theta_{\mathbf{a}}+n\alpha).
$$

\begin{lemma}[Transversality]\label{cl:variation}
One has
\begin{equation}\label{cl:phase-jets}
\partial_{\varepsilon}\theta_{\mathbf{a},W}=\frac{\partial_{\varepsilon}E_W-w_{\mathbf{a}}(W)}{v_{\mathbf{a}}},\qquad \partial_{\varepsilon}E_W=\frac{\sum_{\mathbf{a}}|v_{\mathbf{a}}|^{-1}w_{\mathbf{a}}(W)}{\sum_{\mathbf{a}}|v_{\mathbf{a}}|^{-1}}
\end{equation}
and 
\begin{equation}\label{cl:FH}
\partial_{E}\theta_{\mathbf{a}}=v_{\mathbf{a}}^{-1},\qquad\partial_{E}\mathcal{N}_{0}(E_{*})=\sum_{\mathbf{a}}|v_{\mathbf{a}}|^{-1}>0.
\end{equation}
Moreover, for $\mathbf{a}\neq \mathbf{b}$, the continuous real linear functional
$W\mapsto\partial_{\varepsilon}\theta_{\mathbf{a},W}
-\partial_{\varepsilon}\theta_{\mathbf{b},W}$
on $C_{r}^{\omega}(\mathbb{T},\mathbb{R})$ is nonzero.
\end{lemma}

\begin{proof}
Since $H_{*}$ is independent of $\varepsilon,t$ and satisfies $H_{*}(x)^{*}J_{2N}H_{*}(x)=\widehat J_{2N}$, multiplying \eqref{varid} on the left by $H_{*}(x)^{*}$ and on the right by $H_{*}(x)$ gives the same identity in this frame, with $J_{2N}$ replaced by $\widehat J_{2N}$. The averaging argument in Corollary~\ref{cor:modevar} therefore applies. Apply Corollary~\ref{cor:modevar} to
$C_{\varepsilon,t}$ in \eqref{cl:first-variation-expansion} at $(0,0)$,
with $q=1$ and $z=\mathbf e_{\mathbf{a}}$. By Lemma~\ref{cor:ell}, the column indexed by $\mathbf{a}$ of $U_{0,0}$ is $N^{-\frac{1}{2}}u_{\mathbf{a}}$. Hence
$$
\big\langle(D_{*}^{*}\widehat J_{2N}\delta_{\mathbf{v}}C)_{\mathbf{a}\mathbf{a}}\big\rangle
=\frac{1}{\lambda N}\bigg\langle\overline{u_{\mathbf{a},0}}
\bigg(v_{t}u_{\mathbf{a},0}-v_{\varepsilon}\sum_{r\in\mathbb{Z}}w_{r}u_{\mathbf{a},r}\bigg)\bigg\rangle.
$$

By \eqref{cl:profile-norm},
$
\frac1{\lambda N}\big\langle|u_{\mathbf{a},0}|^{2}\big\rangle=\frac{2\pi}{|v_{\mathbf{a}}|}.
$
Moreover, Parseval's identity gives
$$
\int_{\mathbb{T}}\overline{u_{\mathbf{a},0}(x)}u_{\mathbf{a},0}(x+r\alpha)\,\mathrm{d}x
=\|u_{\mathbf{a},0}\|_{L^2}^2\sum_{n\in\mathbb{Z}}|\psi_{\sigma}(n)|^{2}e^{2\pi inr\alpha}.
$$
It follows that
$$
\frac1{\lambda N}\bigg\langle\overline{u_{\mathbf{a},0}}\sum_{r\in\mathbb{Z}}w_{r}u_{\mathbf{a},r}\bigg\rangle
=\frac{2\pi}{|v_{\mathbf{a}}|}\sum_{r,n\in\mathbb{Z}}w_{r}|\psi_{\sigma}(n)|^{2}e^{2\pi ir(\theta_{\mathbf{a}}+n\alpha)}=\frac{2\pi}{|v_{\mathbf{a}}|}w_{\mathbf{a}}(W).
$$
The sums are absolutely convergent because $\sum_{r}|w_{r}|<\infty$ and $\sum_{n}|\psi_{\sigma}(n)|^{2}=1$.

Using $(\widehat J_{2N})_{\mathbf{a}\mathbf{a}}=i\sigma$
and \eqref{cl:profile-norm}, we obtain
\begin{equation}\label{cl:variation-form}
\frac{1}{2\pi i}\big\langle(D_{*}^{-1}\delta_{\mathbf{v}}C)_{\mathbf{a}\mathbf{a}}\big\rangle
=\frac{1}{2\pi i}[v_t\langle G_{\mathbf{a}\mathbf{a}}\rangle+v_\varepsilon\langle(F_W)_{\mathbf{a}\mathbf{a}}\rangle]
=\frac{v_{t}-v_{\varepsilon}w_{\mathbf{a}}(W)}{v_{\mathbf{a}}}.
\end{equation}
	By \eqref{cl:variation-form} and the first-order expansion of $D_W$ in Lemma~\ref{cl:prep},
	$$
	\theta_{\mathbf{a},W}(\varepsilon,E)
	=\theta_{\mathbf{a}}+\frac{E-E_{*}-\varepsilon w_{\mathbf{a}}(W)}{v_{\mathbf{a}}}
	+O((|\varepsilon|+|E-E_{*}|)^{2}).
	$$
	In particular, $\partial_{E}\theta_{\mathbf{a}}=v_{\mathbf{a}}^{-1}$.

Set $\varepsilon=0$. Lemmas~\ref{lem:rotation-conjugacy}
and~\ref{lem:unitary-rotation}, applied to \eqref{cl:nf} with
$\deg B=0$ and $R=O(t^2)$, together with \eqref{cl:ids} and
$\partial_E\theta_{\mathbf{a}}=v_{\mathbf{a}}^{-1}$, give
$$
\mathcal{N}_{0}(E_{*}+t)=\{m_{0}\alpha\}+t\sum_{\mathbf{a}}|v_{\mathbf{a}}|^{-1}+O(t^{2}).
$$
This proves the formula for the derivative of the IDS in \eqref{cl:FH}.

Along the energy curve (with $\mathbf v=(1,\partial_\varepsilon E_W)$) from Lemma~\ref{cl:prep}, one has
$$
\partial_{\varepsilon}\theta_{\mathbf{a},W}
=\frac{\partial_{\varepsilon}E_W-w_{\mathbf{a}}(W)}{v_{\mathbf{a}}}.
$$
Differentiating \eqref{cl:fixed-rotation}
at $\varepsilon=0$ and using Lemma~\ref{lem:unitary-rotation} gives
\begin{equation}\label{equ:phase}
    \sum_{j=0}^{N-1}\bigl(\partial_{\varepsilon}\theta_{+,j,W}
-\partial_{\varepsilon}\theta_{-,j,W}\bigr)=0.
\end{equation}
Substitution and \eqref{cl:profile-norm} give
$
\partial_{\varepsilon}E_W
=\frac{\sum_{\mathbf{a}}|v_{\mathbf{a}}|^{-1}w_{\mathbf{a}}(W)}{\sum_{\mathbf{a}}|v_{\mathbf{a}}|^{-1}}.
$
This gives the stated formula for $\partial_{\varepsilon}E_W$ and completes \eqref{cl:phase-jets}.

It remains to prove that each pairwise phase difference is not identically zero. Define the probability measures
$$
\mu_{\mathbf{a}}=\sum_{n\in\mathbb{Z}}|\psi_{\sigma}(n)|^{2}\delta_{\theta_{\mathbf{a}}+n\alpha},
\qquad w_{\mathbf{a}}(W)=\int_{\mathbb{T}}W\,\mathrm{d}\mu_{\mathbf{a}}.
$$
By \eqref{cl:seed}, the orbits
$
\mathcal{O}_{\mathbf{a}}=\{\theta_{\mathbf{a}}+n\alpha:n\in\mathbb{Z}\}
$
are pairwise disjoint. Each $\mu_{\mathbf{a}}$ is concentrated on $\mathcal{O}_{\mathbf{a}}$ and has mass one. Evaluating a relation $\sum_{\mathbf{a}}c_{\mathbf{a}}\mu_{\mathbf{a}}=0$ on each $\mathcal{O}_{\mathbf{a}}$ gives $c_{\mathbf{a}}=0$. Thus the measures are linearly independent.

The map $W\mapsto(w_{\mathbf{a}}(W))_{\mathbf{a}}$ is onto $\mathbb{R}^{2N}$. Indeed, otherwise there would be a nonzero real vector $(c_{\mathbf{a}})_{\mathbf{a}}$ such that
$
\int_{\mathbb{T}}W\,\mathrm{d}\bigg(\sum_{\mathbf{a}}c_{\mathbf{a}}\mu_{\mathbf{a}}\bigg)=0
$
for every real analytic $W$. Taking all real trigonometric polynomials would imply $\sum_{\mathbf{a}}c_{\mathbf{a}}\mu_{\mathbf{a}}=0$, contradicting linear independence.

By \eqref{equ:phase}, the phase-variation vectors belong to
$$
\mathcal{H}=\bigg\{y\in\mathbb{R}^{2N}:\sum_{j=0}^{N-1}y_{+,j}=\sum_{j=0}^{N-1}y_{-,j}\bigg\}.
$$
Conversely, fix $y\in\mathcal{H}$. The surjectivity of $W\mapsto(w_{\mathbf{a}}(W))_{\mathbf{a}}$ allows us to choose $W$ such that
$w_{\mathbf{a}}(W)=-v_{\mathbf{a}}y_{\mathbf{a}}$ for every $\mathbf{a}$.
Since $y\in\mathcal{H}$, we have $\sum_{\mathbf{a}}\operatorname{sgn}(v_{\mathbf{a}})y_{\mathbf{a}}=0$. Equations~\eqref{cl:profile-norm} and~\eqref{cl:phase-jets} therefore give
$
\partial_{\varepsilon}E_{W}=0,
\partial_{\varepsilon}\theta_{\mathbf{a},W}
=y_{\mathbf{a}}.
$
Thus the phase-variation vectors fill $\mathcal{H}$.

Fix distinct indices $\mathbf{a},\mathbf{b}$. Since $N\geqslant2$, choose $\mathbf{c}\neq \mathbf{a}$ in the same sign group as $\mathbf{a}$. Then $y=\mathbf{e}_{\mathbf{a}}-\mathbf{e}_{\mathbf{c}}$ belongs to $\mathcal{H}$. Moreover, $y_{\mathbf{a}}=1$ and $y_{\mathbf{b}}\in\{0,-1\}$, so $y_{\mathbf{a}}-y_{\mathbf{b}}\neq0$. Consequently, no coordinate difference vanishes identically on $\mathcal{H}$. Since every vector in $\mathcal{H}$ is realized by a potential $W$, the functional
$
W\mapsto\partial_{\varepsilon}\theta_{\mathbf{a},W}-\partial_{\varepsilon}\theta_{\mathbf{b},W}
$
is nonzero for every $\mathbf{a}\neq \mathbf{b}$.
Finally, $|w_{\mathbf{a}}(W)|\leqslant\|W\|_{r}$, and the preceding formulas express $\partial_{\varepsilon}E_W$ and $\partial_{\varepsilon}\theta_{\mathbf{a},W}$ as fixed real linear combinations of the $w_{\mathbf{a}}(W)$. Their linearity and continuity follow.
\end{proof}

Define $\mathcal{G}_{m_{0}}\subset C_{r}^{\omega}(\mathbb{T},\mathbb{R})$ by
\begin{equation}\label{cl:generic}
\operatorname{Re}w_{1}\neq0,\qquad \partial_{\varepsilon}\theta_{\mathbf{a},W}\neq\partial_{\varepsilon}\theta_{\mathbf{b},W}\quad(\mathbf{a}\neq \mathbf{b}),
\end{equation}
where $w_{1}$ is the first Fourier coefficient of $W$. The first condition ensures that the perturbed potential has fundamental period one. Lemma~\ref{cl:variation} shows that the complement of $\mathcal{G}_{m_{0}}$ is a finite union of proper closed real linear hyperplanes.

\subsection{Positive-measure reducibility at a prescribed IDS value}
Henceforth, we assume $W\in\mathcal{G}_{m_0}$ and suppress the dependence on $W$ in $D$ and $E^{(0)}$. Fix $\tau_{1}>\max\{\tau+1,2\}$.

\begin{lemma}\label{cl:KAM}
There are $h_{0},\delta_{0},c,C>0$ such that, for $0<\delta<\delta_{0}$ and $|\vartheta-\{m_{0}\alpha\}|\leqslant\delta$, a compact set $\mathcal{E}_{\delta}(\vartheta)\subset[\delta,2\delta]$ satisfies
\begin{equation}\label{cl:measure}
\operatorname{Leb}([\delta,2\delta]\setminus\mathcal{E}_{\delta}(\vartheta))\leqslant C\delta^{3}.
\end{equation}
For every $\varepsilon\in\mathcal{E}_{\delta}(\vartheta)$, there are $E_{\infty}\in\mathbb{R}$, $B_{\infty}\in C^\omega(\mathbb T, {\mathrm{HSp}}(\widehat J_{2N}))$, and a diagonal unitary matrix $D_{\infty}$ such that
$$
B_{\infty}(x+\alpha)^{-1}C_{\varepsilon,E_{\infty}-E_{*}}(x)B_{\infty}(x)=D_{\infty},\qquad \mathcal{N}_{\varepsilon,W}(E_{\infty})=\vartheta,
$$
and
$$
|E_{\infty}-E_{*}|+\|B_{\infty}-I_{2N}\|_{h_{0}/2}\leqslant C\delta.
$$
Writing $D_{\infty}=\operatorname{diag}(e^{2\pi i\theta_{\mathbf{a}}^{(\infty)}})$, one has
\begin{equation}\label{cl:final-NR}
\|\theta_{\mathbf{a}}^{(\infty)}-\theta_{\mathbf{b}}^{(\infty)}-n\alpha\|_{\mathbb{T}}\geqslant c\delta^{2}(1+|n|)^{-\tau_{1}-1},\qquad \mathbf{a}\neq \mathbf{b},\quad n\in\mathbb{Z}.
\end{equation}
\end{lemma}

\begin{proof}

Fix $\vartheta$ as in the statement. Choose a strip $|\operatorname{Im}x|<h_0\le1$
and a bidisc $|\varepsilon|,|E-E_*|<s\ll1$ within the domain of
Lemma~\ref{cl:prep}. Then take $0<\delta\ll s$ and set
$h_k=\frac{h_0}{2}(1+2^{-k})$ and $\mathcal E_0=[\delta,2\delta]$.
For a compact parameter set $\mathcal E_k\subset\mathcal E_0$, set
\begin{equation}\label{eq:collapse-parameter-norm}
\begin{aligned}
 \|H\|_k
 =\sup_{\substack{\varepsilon\in\mathcal E_k\\|E-E_*|<s}}
 \|H(\cdot,\varepsilon,E)\|_{h_k}+\delta\sup_{\substack{\varepsilon,\varepsilon'\in\mathcal E_k,\ \varepsilon\ne\varepsilon'\\|E-E_*|<s}}
 \frac{\|H(\cdot,\varepsilon,E)-H(\cdot,\varepsilon',E)\|_{h_k}}
 {|\varepsilon-\varepsilon'|}.
\end{aligned}
\end{equation}
For functions independent of $E$, we omit $E$ from the suprema.
All estimates below
are uniform in $\vartheta$.

\noindent\emph{1. Initial step.} We suppress the dependence on $\varepsilon$ when no confusion arises.
Take $E^{(0)}=E_W^{(0)}(\varepsilon,\vartheta)$ from
Lemma~\ref{cl:prep}, and put $B_0(x)=B(\varepsilon,E^{(0)},x)$
and $D_0(E)=D_W(\varepsilon,E)$.
We construct families
\begin{equation}\label{cl:revised-family}
 B_k(x+\alpha)^{-1}C_{\varepsilon,E-E_*}(x)B_k(x)
 =D_k(E)e^{f_k(E,x)}.
\end{equation}
Here $B_k$ is independent of $E$, and $f_k$ is defined by the logarithm near
$I_{2N}$. For real variables, $B_k\in\mathrm{HSp}(\widehat J_{2N})$,
$f_k\in\mathfrak{hsp}(\widehat J_{2N})$, and $D_k$ is diagonal unitary.
Choose $E^{(k)}$ by
\begin{equation}\label{cl:revised-rotation}
 \rho\bigl(\alpha,D_k(E^{(k)})\bigr)
 =\rho(\alpha,D_*)+\{m_0\alpha\}-\vartheta.
\end{equation}
Write $r_k=\|f_k(E^{(k)})\|_k$. Initially, Lemma~\ref{cl:prep} gives
$\|E^{(0)}-E_*\|_0+\|B_0-I_{2N}\|_0\le C\delta$
and $r_0\le C\delta^{12}$.
Indeed, $f_0(E^{(0)})=R(\varepsilon,E^{(0)},\cdot)$,
$|\partial_\varepsilon E^{(0)}|\le C$, and $\|R\|_{h_0}\le C\delta^{12}$
on a bidisc of radius $\delta$ about $(\varepsilon,E^{(0)}(\varepsilon))$.
Cauchy's estimate and the chain rule at this point give
$$
 \left\|\partial_\varepsilon R(\varepsilon,E^{(0)})\right\|_{h_0}
 =\left\|\partial_\varepsilon R+\partial_E R\,\partial_\varepsilon E^{(0)}\right\|_{h_0}
 \le C\delta^{11}.
$$
The same calculation applies to $B_0-I_{2N}$.

Write $D_{0}(E^{(0)})_{\mathbf{a}\mathbf{a}}=e^{2\pi i\theta_{\mathbf{a}}^{(0)}}$, choosing real phase lifts near those of $D_{*}$. By Lemma~\ref{cl:prep}, the functions $\theta_{\mathbf{a}}^{(0)}(\varepsilon,\vartheta)$ are jointly real analytic on a fixed neighborhood of $(0,\{m_{0}\alpha\})$. Lemma~\ref{cl:variation} gives
$
\partial_{\varepsilon}\theta_{\mathbf{a}}^{(0)}(0,\{m_{0}\alpha\})=\partial_{\varepsilon}\theta_{\mathbf{a},W}.$
 Therefore 
$
\partial_{\varepsilon}(\theta_{\mathbf{a}}^{(0)}-\theta_{\mathbf{b}}^{(0)})
=\partial_{\varepsilon}\theta_{\mathbf{a},W}-\partial_{\varepsilon}\theta_{\mathbf{b},W}
+O(\delta),
$
uniformly in these parameters. By \eqref{cl:generic}, the leading term is nonzero for every $\mathbf{a}\neq \mathbf{b}$.

\medskip
\noindent\emph{2. KAM iteration.}
Assume $\|B_k-I_{2N}\|_k+\|E^{(k)}-E_*\|_k\le C\delta$
and $r_j\le C\delta^{12}2^{-(4\tau_1+6)j}$ for $0\le j\le k$.
These bounds hold at $k=0$ by Lemma~\ref{cl:prep}.
On $\mathcal E_k$, set $D_{k+1}(E)=D_k(E)\exp\langle f_k(E)\rangle_{\mathrm{diag}}$.

We also maintain $\|D_k(E)-D_*\|_k\le Cs$ on the fixed energy disc, starting from
Lemma~\ref{cl:prep}.
Equation~\eqref{cl:revised-family} gives $\|f_k(E)\|_k\le Cs$
on the fixed disc. Since the diagonal factors commute, the
Baker--Campbell--Hausdorff formula in the norm $\|\cdot\|_k$ gives
$$
\begin{aligned}
 \log(D_*^{-1}D_{k+1}(E))
 &=\log(D_*^{-1}D_k(E))+\langle f_k(E)\rangle_{\mathrm{diag}}\\
 &=\left\langle
 \log\!\left(D_*^{-1}B_k(x+\alpha)^{-1}
 C_{\varepsilon,E-E_*}(x)B_k(x)\right)
 \right\rangle_{\mathrm{diag}}+O_{\|\cdot\|_k}(s^2).
\end{aligned}
$$
Thus $\|D_{k+1}(E)-D_*\|_k\le C\bigl(s+\|D_k(E)-D_*\|_k^2\bigr)$,
which preserves a fixed $Cs$ bound for small $s$.
By $\partial_E B_k=0$ and \eqref{cl:first-variation-expansion},
$$
 \partial_E\log\!\left(D_*^{-1}B_k(x+\alpha)^{-1}
 C_{\varepsilon,E-E_*}(x)B_k(x)\right)
 =G(x)+O_{\|\cdot\|_k}(s).
$$
Restricting to real $E\in(E_*-s/2,E_*+s/2)$, Cauchy's estimate gives,
with energy suprema over this interval,
$\partial_E\log(D_*^{-1}D_{k+1})
=\langle G\rangle_{\mathrm{diag}}+O_{\|\cdot\|_k}(s)$ and
$$
 \|\partial_E^j f_k(E)\|_k+
 \|\partial_E^j\log(D_*^{-1}D_{k+1}(E))\|_k\le Cs^{1-j},
 \qquad j=1,2.
$$
For a diagonal unitary constant $D$, write $D_{\mathbf{a}\mathbf{a}}=e^{2\pi i\theta_{\mathbf{a}}(D)}$.
 \eqref{eq:unitary-rotation}, \eqref{cl:profile-norm} and \eqref{cl:FH}, give
$$
 \partial_E\rho(\alpha,D_{k+1}(E))
 =-\partial_E\mathcal N_0(E_*)+O(s)
 \le-\tfrac12\partial_E\mathcal N_0(E_*)<0.
$$
For a real increment $u$, \eqref{eq:unitary-rotation} and
\eqref{cl:revised-rotation} give
$$
\begin{aligned}
 \rho\bigl(\alpha,&D_{k+1}(\varepsilon,E^{(k)}(\varepsilon)+u)\bigr)
   -\rho(\alpha,D_*)-\{m_0\alpha\}+\vartheta\\
 &=-\frac1{2\pi}\operatorname{tr}
   \bigl(\widehat J_{2N}\langle f_k(\varepsilon,E^{(k)}(\varepsilon))\rangle_{\mathrm{diag}}\bigr)+u\int_0^1\partial_E\rho\bigl(
   \alpha,D_{k+1}(\varepsilon,E^{(k)}(\varepsilon)+tu)\bigr)\,dt.
\end{aligned}
$$
Since $\|\langle f_k(E^{(k)})\rangle_{\mathrm{diag}}\|_k\le r_k$,
the implicit function theorem gives a unique nearby $E^{(k+1)}$
satisfying \eqref{cl:revised-rotation} and
$\sup_{\varepsilon\in\mathcal E_k}|E^{(k+1)}(\varepsilon)-E^{(k)}(\varepsilon)|\le Cr_k$.
For fixed $|u|\le Cr_k$, the derivative bounds give
$$
\big\|\int_0^1\partial_E\rho(\alpha,
D_{k+1}(\varepsilon,E^{(k)}(\varepsilon)+tu))\,dt\big\|_k
\le C(1+s^{-1}\|E^{(k)}-E_*\|_k)\le C.
$$
Hence, with the same $u$ at both amplitudes,
$$
 \delta\big|\rho\bigl(\alpha,D_{k+1}(\varepsilon,E^{(k)}(\varepsilon)+u)\bigr)-\rho\bigl(\alpha,D_{k+1}(\varepsilon',E^{(k)}(\varepsilon')+u)\bigr)\big|
 \le C(r_k+|u|)|\varepsilon-\varepsilon'|.
$$
Take $u=E^{(k+1)}(\varepsilon)-E^{(k)}(\varepsilon)$ and apply the
slope bound at $\varepsilon'$ to obtain
$\delta\operatorname{Lip}(E^{(k+1)}-E^{(k)})\le Cr_k$.
The energy derivative bounds also give
\begin{equation}\label{cl:revised-energy-bound}
 \|E^{(k+1)}-E^{(k)}\|_k+\|f_k(E^{(k+1)})\|_k
 +\|D_{k+1}(E^{(k+1)})-D_k(E^{(k)})\|_k\le Cr_k.
\end{equation}

By \eqref{cl:revised-energy-bound}, for
$\varepsilon,\varepsilon'\in\mathcal E_k$, $\varepsilon\ne\varepsilon'$,
$$
 \frac{|(\theta_{\mathbf{a}}^{(k+1)}-\theta_{\mathbf{b}}^{(k+1)})(\varepsilon)
 -(\theta_{\mathbf{a}}^{(k+1)}-\theta_{\mathbf{b}}^{(k+1)})(\varepsilon')|}
 {|\varepsilon-\varepsilon'|}
 \ge c-C\delta^{-1}\sum_{j=0}^{k}r_j\ge\frac c2,
$$
by the inductive bounds on $r_j$.
Restrict to the compact set
$\mathcal E_{k+1}\subset\mathcal E_k$ on which
\begin{equation}\label{cl:revised-exclusion}
 \|\theta_{\mathbf{a}}^{(k+1)}-\theta_{\mathbf{b}}^{(k+1)}-n\alpha\|_{\mathbb T}
 \ge 2^{-k}\delta^2(1+|n|)^{-\tau_1},
 \qquad \mathbf{a}\ne \mathbf{b},\quad n\in\mathbb Z.
\end{equation}

The initial phase estimates and \eqref{cl:revised-energy-bound} give
$
|\theta_{\mathbf{a}}^{(k+1)}-\theta_{\mathbf{a}}|\leqslant C\delta,
$
where $\theta_{\mathbf{a}}$ are the phases of $D_{*}$. By Lemma~\ref{lem:seed}, any resonance excluded in \eqref{cl:revised-exclusion} therefore satisfies
$$
c(1+|n|)^{-\tau}\leqslant C\delta+2^{-k}\delta^{2}(1+|n|)^{-\tau_{1}}\leqslant C\delta.
$$
Hence $1+|n|\geqslant c\delta^{-1/\tau}$. The difference-quotient bound and the bounded range of the phases now give
$$
\begin{aligned}
\operatorname{Leb}(\mathcal{E}_{k}\setminus\mathcal{E}_{k+1})
\leqslant C2^{-k}\delta^{2}\sum_{1+|n|\geqslant c\delta^{-1/\tau}}(1+|n|)^{-\tau_{1}}\leqslant C2^{-k}\delta^{2+(\tau_{1}-1)/\tau}\leqslant C2^{-k}\delta^3.
\end{aligned}
$$

On $\mathcal E_{k+1}$, solve the homological equation with
$\langle X_k\rangle_{\mathrm{diag}}=0$:
$$
 D_{k+1}(E^{(k+1)})^{-1}X_k(x+\alpha)D_{k+1}(E^{(k+1)})-X_k(x)=f_k(E^{(k+1)},x)
       -\langle f_k(E^{(k+1)})\rangle_{\mathrm{diag}}.
$$
For $\mathbf{a}\ne \mathbf{b}$, \eqref{cl:revised-exclusion} gives
$|e^{2\pi i(n\alpha+\theta_{\mathbf{b}}^{(k+1)}-\theta_{\mathbf{a}}^{(k+1)})}-1|
\ge c2^{-k}\delta^2(1+|n|)^{-\tau_1}$.
We obtain
$$
 \|(e^{2\pi i(n\alpha+\theta_{\mathbf{b}}^{(k+1)}-\theta_{\mathbf{a}}^{(k+1)})}-1)^{-1}\|_{k+1}
 \le C\delta^{-4}2^{2k}(1+|n|)^{2\tau_1}.
$$
For $\mathbf{a}= \mathbf{b}$, $n\ne0$, the same bound follows from the Diophantine
condition on $\alpha$ and $\tau_1>\tau$. Thus
$$
 \|X_k\|_{k+1}\le C\delta^{-4}2^{2k}r_k
 \sum_{n\in\mathbb Z}(1+|n|)^{2\tau_1}
 e^{-2\pi|n|(h_k-h_{k+1})}\le C\delta^{-4}2^{(2\tau_1+3)k}r_k,
$$
since $h_k-h_{k+1}=h_0\,2^{-k-2}$ and the sum is
$O((h_k-h_{k+1})^{-2\tau_1-1})$.
As in Lemma~\ref{cl:prep}, uniqueness gives
$X_k\in\mathfrak{hsp}(\widehat J_{2N})$ for real variables.

Set $B_{k+1}=B_ke^{X_k}$ and define $f_{k+1}$ by
\eqref{cl:revised-family}.
At energy $E^{(k+1)}$, suppressing the
energy arguments, the updates give the exact product
$$
 e^{f_{k+1}(x)}
 =e^{-D_{k+1}^{-1}X_k(x+\alpha)D_{k+1}}
  e^{-\langle f_k\rangle_{\mathrm{diag}}}e^{f_k(x)}e^{X_k(x)}.
$$
The linear terms cancel by the homological equation, so Baker--Campbell--Hausdorff formula gives
$$
 r_{k+1}\le C(\|X_k\|_{k+1}+r_k)^2
 \le C\delta^{-8}2^{(4\tau_1+6)k}r_k^2.
$$
For small $\delta$, this gives $r_{k+1}\le C\delta^{12}2^{-(4\tau_1+6)(k+1)}$,
while the bound for $X_k$ gives $\|X_k\|_{k+1}\le C\delta^8 2^{-(2\tau_1+3)k}$.
Consequently,
$$
 \|B_{k+1}-I_{2N}\|_{k+1}+\|E^{(k+1)}-E_*\|_{k+1}
 \le C\delta+C\sum_{j=0}^{k}(\|X_j\|_{j+1}+r_j)\le C\delta
$$
and $\delta^{-1}\sum_{j=0}^{k}r_j\le C\delta^{11}$.
Together with $C\delta<s/4$ and $\|D_{k+1}-D_*\|_{k+1}\le Cs$,
these bounds close the induction.

\medskip
\noindent\emph{3. Convergence.}
The estimates above give
$\sum_{k\ge0}\|X_k\|_{k+1}\le C\delta^8$ and
$\sum_{j\ge k}r_j\le C\delta^{12}2^{-(4\tau_1+6)k}$.

The compact intersection
$\mathcal E_\delta(\vartheta)=\bigcap_{k\ge0}\mathcal E_k$ satisfies
\eqref{cl:measure}. 
On this set $E^{(k)}$, $B_k$ and $D_k(E^{(k)})$
converge to $E_\infty$, $B_\infty$ and $D_\infty$, respectively.
The bounds give $|E_\infty-E_*|+\|B_\infty-I_{2N}\|_{h_0/2}\le C\delta$.
Passing to the limit in \eqref{cl:revised-family} at $E=E^{(k)}$
gives the stated conjugacy. The limit $B_\infty$ belongs to
$\mathrm{HSp}(\widehat J_{2N})$ and has degree zero by the preceding bound.
By Lemma~\ref{lem:rotation-conjugacy},
\eqref{cl:ids} and \eqref{cl:revised-rotation},
$\mathcal N_{\varepsilon,W}(E_\infty)
=\{m_0\alpha\}-\rho(\alpha,D_\infty)+\rho(\alpha,D_*)=\vartheta$.

Finally, for each $n$ take $k=\lceil\log_2(1+|n|)\rceil$ in
\eqref{cl:revised-exclusion}. Since
$\max_{\mathbf{a}}|\theta_{\mathbf{a}}^{(\infty)}-\theta_{\mathbf{a}}^{(k+1)}|
\le C\sum_{j\ge k+1}r_j\le C\delta^{12}2^{-(4\tau_1+6)(k+1)}$, we have
$$
 \|\theta_{\mathbf{a}}^{(\infty)}-\theta_{\mathbf{b}}^{(\infty)}-n\alpha\|_{\mathbb T}
 \ge \tfrac12\delta^2(1+|n|)^{-\tau_1-1}
      -C\delta^{12}(1+|n|)^{-(4\tau_1+6)}\ge \tfrac14\delta^2(1+|n|)^{-\tau_1-1}
$$
for small $\delta$, since $4\tau_1+6>\tau_1+1$. This proves
\eqref{cl:final-NR}. 
\end{proof}

\subsection{Derivative of the IDS and gap collapse}

Fix $\varepsilon\in\mathcal{E}_{\delta}(\vartheta)$, and let $E_{\infty}$, $B_{\infty}$ and $D_{\infty}$ be given by Lemma~\ref{cl:KAM}. Let $\phi_{\mathbf{a}}(x)$ be the zeroth coordinate of column $\mathbf{a}$ of the solution frame $U_{\varepsilon,E_{\infty}-E_{*}}(x)B_{\infty}(x)$.

\begin{lemma}\label{cl:collapse}
One has
\begin{equation}\label{cl:exact-IDS-derivative}
\partial_{E}\mathcal{N}_{\varepsilon,W}(E_{\infty})=\frac{1}{2\pi\lambda}\sum_{\mathbf{a}}\|\phi_{\mathbf{a}}\|_{L^{2}(\mathbb{T})}^{2}>0.
\end{equation}
Moreover,
\begin{equation}\label{cl:IDS-derivative}
\mathcal{N}_{\varepsilon,W}(E)=\vartheta+\partial_{E}\mathcal{N}_{\varepsilon,W}(E_{\infty})(E-E_{\infty})+O((E-E_{\infty})^{2}),
\end{equation}
and $\{E:\mathcal{N}_{\varepsilon,W}(E)=\vartheta\}=\{E_{\infty}\}$.
\end{lemma}

\begin{proof}
Set
$
U(E,x)=U_{\varepsilon,E-E_{*}}(x)B_{\infty}(x).
$
Then
$$
SU(E,x)=U(E,x+\alpha)D_{\infty}e^{F_{\infty}(E,x)},\qquad F_{\infty}(E_{\infty},x)=0,
$$
where $F_{\infty}$ is analytic near $E_{\infty}$. Here $\varepsilon$ is fixed and $B_{\infty}$ is independent of $E$, so Corollary~\ref{cor:modevar} also holds in this frame.
Corollary~\ref{cor:modevar}, applied to this
frame at $E=E_\infty$ with $q=1$ and $z=\mathbf e_{\mathbf{a}}$, gives
\begin{equation}\label{cl:energy-positivity}
\bigl(\widehat J_{2N}\langle\partial_E F_\infty(E_\infty)\rangle_{\mathrm{diag}}\bigr)_{\mathbf{a}\mathbf{a}}
=\frac1\lambda\|\phi_{\mathbf{a}}\|_{L^2(\mathbb T)}^2>0.
\end{equation}

The Diophantine condition on $\alpha$ and \eqref{cl:final-NR} allow us to solve
$$
D_{\infty}^{-1}X(x+\alpha)D_{\infty}-X(x)=-\partial_{E}F_{\infty}(E_{\infty},x)+\big\langle\partial_{E}F_{\infty}(E_{\infty})\big\rangle_{\mathrm{diag}},\qquad \langle X\rangle_{\mathrm{diag}}=0.
$$
The solution is analytic on a smaller strip and belongs to $\mathfrak{hsp}(\widehat J_{2N})$ for real $x$. Expanding the conjugacy gives
\begin{equation}\label{cl:energy-normal-form}
\begin{aligned}
&e^{(E-E_{\infty})X(x+\alpha)}D_{\infty}e^{F_{\infty}(E,x)}e^{-(E-E_{\infty})X(x)}\\
&\qquad=D_{\infty}\exp\big((E-E_{\infty})\big\langle\partial_{E}F_{\infty}(E_{\infty})\big\rangle_{\mathrm{diag}}\big)e^{O((E-E_{\infty})^{2})}.
\end{aligned}
\end{equation}
Applying
Lemmas~\ref{lem:rotation-conjugacy} and~\ref{lem:unitary-rotation}
to \eqref{cl:energy-normal-form}, whose conjugacies have degree zero
and whose constant factor is diagonal unitary, gives
$$
\begin{aligned}
\rho(\alpha,C_{\varepsilon,E-E_{*}})-\rho(\alpha,D_{\infty})
&=-\frac{E-E_{\infty}}{2\pi}\operatorname{tr}\big(\widehat J_{2N}\big\langle\partial_{E}F_{\infty}(E_{\infty})\big\rangle_{\mathrm{diag}}\big)+O((E-E_{\infty})^{2})\\
&=-\frac{E-E_{\infty}}{2\pi\lambda}\sum_{\mathbf{a}}\|\phi_{\mathbf{a}}\|_{L^{2}(\mathbb{T})}^{2}+O((E-E_{\infty})^{2}).
\end{aligned}
$$
The IDS identity \eqref{cl:ids} now proves \eqref{cl:IDS-derivative} and the exact expression in \eqref{cl:exact-IDS-derivative}.
Hence \eqref{cl:IDS-derivative} gives $\mathcal{N}_{\varepsilon,W}(E)<\vartheta$ immediately to the left of $E_{\infty}$ and $\mathcal{N}_{\varepsilon,W}(E)>\vartheta$ immediately to its right. Monotonicity of  IDS then makes its entire level set at $\vartheta$ the single point $E_{\infty}$.
\end{proof}

\subsection*{Proof of Theorem~\ref{thm:collapsed-gaps}}
\pushQED{\qed}
Fix $m_0\notin N\mathbb Z$ and set
$\mathcal G_{\mathrm{cl}}=\mathcal G_{m_0}$.
The assertion about $\mathcal G_{\mathrm{cl}}$ follows from
Lemma~\ref{cl:variation}. Fix $W\in\mathcal G_{\mathrm{cl}}$.
For sufficiently small $\delta>0$, we consider
$v_0+u+\varepsilon W$, where $\|u\|_r\leqslant\delta^{12}$ and
$\varepsilon\in[\delta,2\delta]$. We prove that
\begin{equation}\label{cl:translated-cubic}
 \operatorname{Leb}\{\varepsilon\in[\delta,2\delta]:
        v_0+u+\varepsilon W\notin\mathcal S_N\}
 \leqslant C\delta^3,
 \qquad \|u\|_r\leqslant\delta^{12},
\end{equation}
with $C$ independent of $u$ and $\delta$.

The Neumann-series construction in Proposition~\ref{prop:slow}
is analytic in the total perturbation $u+\varepsilon W$.
Thus adding $u$ changes the slow cocycle and, by Cauchy's estimates,
its first derivatives in $\varepsilon$ and $E$ by $O(\|u\|_r)$
on a smaller fixed domain. Use the same $B_0,D_0,E^{(0)}$ as in
the proof of Lemma~\ref{cl:KAM}. The initial remainder satisfies
$$
 r_0\leqslant C(\delta^{12}+\|u\|_r)\leqslant2C\delta^{12}.
$$ The initial phase
derivatives at fixed rotation are unchanged, and the energy
derivative estimates remain valid for small $\delta$.
The same iteration therefore applies, uniformly in $u$ and
$|\vartheta-\{m_0\alpha\}|\leqslant\delta$. Write
$\mathcal E_\delta(\vartheta)\subset[\delta,2\delta]$ for the
resulting compact set, suppressing its dependence on $u$.
Then $\operatorname{Leb}([\delta,2\delta]\setminus
                   \mathcal E_\delta(\vartheta))
 \leqslant C\delta^3.$
By Lemmas~\ref{lem:localids} and~\ref{cl:collapse}, at every
$\varepsilon\in\mathcal E_\delta(\vartheta)$ the IDS of
$v_0+u+\varepsilon W$ has a singleton level set at $\vartheta$.

For each $a\in\{0,\ldots,N-1\}$, choose distinct positive integers
$m_{a,k}\to\infty$ such that
$$
 m_{a,k}\equiv a\pmod N,\qquad
 \vartheta_{a,k}:=\{m_{a,k}\alpha\},\qquad
 |\vartheta_{a,k}-\{m_0\alpha\}|<\delta.
$$
Such a sequence exists because every rotation by $N\alpha$ is
dense in $\mathbb T$. Define
$$
\begin{aligned}
 \mathcal E_{\delta,a}:=\limsup_{k\to\infty}
          \mathcal E_\delta(\vartheta_{a,k})
 =\bigcap_{K\geqslant1}\bigcup_{k\geqslant K}
          \mathcal E_\delta(\vartheta_{a,k}),\qquad\mathcal E_\delta:=\bigcap_{a=0}^{N-1}\mathcal E_{\delta,a}.
\end{aligned}
$$
For fixed $a$, each tail union
$\bigcup_{k\geqslant K}\mathcal E_\delta(\vartheta_{a,k})$
has measure at least $\delta-C\delta^3$.
These unions decrease with $K$, so continuity of measure from above gives
$$
 \operatorname{Leb}(\mathcal E_{\delta,a})
 =\lim_{K\to\infty}\operatorname{Leb}\bigg(
       \bigcup_{k\geqslant K}\mathcal E_\delta(\vartheta_{a,k})\bigg)
 \geqslant\delta-C\delta^3.
$$
Consequently,
$$
\begin{aligned}
 \operatorname{Leb}([\delta,2\delta]\setminus\mathcal E_{\delta,a})
 \leqslant C\delta^3,\qquad
 \operatorname{Leb}([\delta,2\delta]\setminus\mathcal E_\delta)
 \leqslant\sum_{a=0}^{N-1}
       \operatorname{Leb}([\delta,2\delta]\setminus\mathcal E_{\delta,a})
 \leqslant NC\delta^3.
\end{aligned}
$$
Every $\varepsilon\in\mathcal E_\delta$ belongs to
$\mathcal E_\delta(\vartheta_{a,k})$ for infinitely many $k$,
for each $a$. Thus the same potential has infinitely many
collapsed labels in each residue class modulo $N$. Moreover,
$N\geqslant2$ and \eqref{cl:generic} give
$$
 |\widehat{v_0+u+\varepsilon W}(1)|
 =|\widehat u(1)+\varepsilon w_1|
 \geqslant\delta|w_1|-\delta^{12}>0
$$
for small $\delta$, so the potential has fundamental period one.
This proves \eqref{cl:translated-cubic}, after increasing $C$.
Taking $u=0$ and summing over the intervals
$[2^{-j}\delta,2^{-j+1}\delta]$, $j\geqslant1$, gives
$$
 \operatorname{Leb}\{\varepsilon\in(0,\delta):
              v_{\varepsilon,W}\notin\mathcal S_N\}
 \leqslant C\sum_{j\geqslant1}(2^{-j}\delta)^3
 =\frac C7\delta^3.
$$
Subtracting from $\delta$ proves the stated measure estimate.

The set $\mathcal S_N$ is Borel. Indeed, for each label, the set
of potentials for which it is collapsed is closed. The condition
of infinitely many collapsed labels in each residue class is
therefore expressed by countable unions and intersections of
closed sets. The condition of fundamental period one is the
complement of $\bigcup_{d\geqslant2}\{v:v(\cdot+1/d)=v\}$.

Let $\mathcal U$ be an open neighborhood of $v_0$.
Choose $\delta>0$ so small that $C\delta^2<1$ and
$v_0+u+\varepsilon W\in\mathcal U$ whenever
$\|u\|_r\leqslant\delta^{12}$ and $\varepsilon\in[\delta,2\delta]$.
If $S=\mathcal S_N\cap\mathcal U$ were shy, then by definition
there would exist a compactly supported probability measure $\mu$
such that $\mu(S-h)=0$ for every $h$.
By restricting $\mu$ to a ball of positive measure, then translating
and normalizing it, we may assume that
$\operatorname{supp}\mu\subset\{u:\|u\|_r\leqslant\delta^{12}\}$
while preserving this property; see \cite[Fact~2]{HSY92}.
Fubini's theorem and \eqref{cl:translated-cubic} now give
$$
\begin{aligned}
 0
 =\int_\delta^{2\delta}
       \mu(S-v_0-\varepsilon W)\,d\varepsilon=\int \operatorname{Leb}\{\varepsilon\in[\delta,2\delta]:
       v_0+u+\varepsilon W\in S\}\,d\mu(u)\geqslant\delta-C\delta^3>0,
\end{aligned}
$$
a contradiction. Thus $\mathcal S_N\cap\mathcal U$ is non-shy.

Now we prove \eqref{uniform}. Fix $m_{0}\in\mathbb{Z}\setminus N\mathbb{Z}$ and set
$\varrho_{*}=\frac{1-\{m_{0}\alpha\}}{2}, E_{*}=2\lambda\cos(2\pi\varrho_{*}).$
For $f=v_{0}+u+\varepsilon W$ and $\varepsilon_{*}>0$, consider the rescaled operator
$$
\big(\varepsilon_{*}H_{\varepsilon_{*}^{-1}f,\alpha,x}z\big)_{m}
=\varepsilon_{*}(z_{m+1}+z_{m-1})+f(x+m\alpha)z_{m}.
$$
At $\varepsilon_{*}=0$, it is the multiplication operator by $f(x+m\alpha)$. Denote its IDS by $\mathcal{N}_{\varepsilon_{*},f}$. Positive rescaling of the operator and energy preserves gap labels. Its dual equation is
$$
\sum_{n\in\mathbb{Z}}\widehat f_{n}z_{m+n}+2\varepsilon_{*}\cos(2\pi(x+m\alpha))z_{m}=Ez_{m}.
$$
At $(\varepsilon_{*},u,\varepsilon,E)=(0,0,0,E_{*})$, this becomes
$\lambda(z_{m+N}+z_{m-N})=E_{*}z_{m}.$
The solutions $z_{m}=e^{2\pi im\theta_{\sigma,j}}$, where
$
\theta_{\sigma,j}=\frac{\sigma\varrho_{*}+j}{N}, \sigma\in\{-,+\}, 0\leqslant j<N,
$
give a basis of the solution space. Their initial vectors diagonalize the one-step cocycle. And  its positive and negative iterates are bounded. The proof of Lemma~\ref{lem:seed} also gives
$$
\|\theta_{\mathbf{a}}-\theta_{\mathbf{b}}-n\alpha\|_{\mathbb{T}}\geqslant c(1+|n|)^{-\tau},\qquad \mathbf{a}\neq\mathbf{b},\quad n\in\mathbb{Z}.
$$

We now repeat the constructions in the proofs of Propositions~\ref{prop:slow1} and~\ref{prop:slow}. The growth estimate needed for the inverse symbol follows from the bounded iterates above. On a fixed complex phase strip, multiplication by $2\varepsilon_{*}\cos(2\pi(x+m\alpha))$ has norm at most $C|\varepsilon_{*}|$. Fix $a<b<2\pi r$. 
Choosing the symbol contours and weights below $b$, the added terms therefore have norm at most $C(|\varepsilon_{*}|+\|u\|_{r})$ in the same Banach algebras. The Neumann series, the contour projection and the normalization of the Green form give rank-$2N$ solution frames and one-step cocycles
$
C_{\varepsilon_{*},u}(\varepsilon,E,x)\in\mathrm{HSp}(\widehat J_{2N})
$
for real variables. The frames and cocycles are analytic on a fixed parameter domain. We use the same constant change of coordinates at the reference point for all $\varepsilon_{*},u$, so that
$
C_{0,0}(0,E_{*})=D_{*}:=\operatorname{diag}(e^{2\pi i\theta_{\mathbf{a}}}).
$
After restricting to a smaller fixed domain, the difference estimates in those constructions and Cauchy's estimates give
\begin{equation}\label{cl:uniform-comparison} 
\max_{\substack{0\leqslant i\leqslant1\\0\leqslant j\leqslant2}}\ \sup_{\substack{|\varepsilon|<s\\|E-E_{*}|<s}}
\|\partial_{\varepsilon}^{i}\partial_{E}^{j}(C_{\varepsilon_{*},u}-C_{0,0})(\varepsilon,E,\cdot)\|_{h}
\leqslant C(|\varepsilon_{*}|+\|u\|_{r}).
\end{equation}
The domain and constants are independent of sufficiently small $\varepsilon_{*},u$.

 At $\varepsilon_{*}=0$, the equation has constant coefficients, so the slow cocycle is independent of $x$. For real $f$ close to $v_{0}$ and $E$ close to $E_{*}$, the equation $f(\theta)=E$ has exactly $2N$ simple real roots near the phases $\theta_{\sigma,j}$. Choose their real lifts continuously. Counting the intervals where $f\leqslant E$ and applying Lemma \ref{lem:unitary-rotation} give
$$
\mathcal{N}_{0,f}(E)=1+\sum_{j=0}^{N-1}\theta_{-,j}(f,E)-\sum_{j=0}^{N-1}\theta_{+,j}(f,E)
=1-\rho\big(\alpha,C_{0,u}(\varepsilon,E)\big).
$$
Here the lift is chosen so that $\rho(\alpha,D_{*})=2\varrho_{*}=1-\{m_{0}\alpha\}$.

For polynomial $f$ and $\varepsilon_{*}>0$, the finite-range comparison in the proof of Lemma~\ref{lem:localids}, applied to $\varepsilon_{*}^{-1}f$ at energy $\varepsilon_{*}^{-1}E$, gives
$$
\mathcal{N}_{\varepsilon_{*},f}(E)+\rho\big(\alpha,C_{\varepsilon_{*},u}(\varepsilon,E)\big)\in\mathbb{Z}+\alpha\mathbb{Z}.
$$
For fixed $f$, the left-hand side is continuous on the connected real $(\varepsilon_{*},E)$ domain with $\varepsilon_{*}>0$. It is therefore constant. As $\varepsilon_{*}\downarrow0$, the operators converge in norm to multiplication by $f(x+m\alpha)$. Since the distribution function of $f$ is continuous near $E_{*}$, their IDS converge there. The cocycles and their rotation numbers also converge. The identity at $\varepsilon_{*}=0$ thus shows that this constant is $1$. Fourier truncation, the uniform frame estimates and continuity of the IDS and rotation number extend the identity to analytic $f$. Consequently,
\begin{equation}\label{cl:scaled-ids}
\mathcal{N}_{\varepsilon_{*},f}(E)-\{m_{0}\alpha\}
=-\rho\big(\alpha,C_{\varepsilon_{*},u}(\varepsilon,E)\big)+\rho(\alpha,D_{*}).
\end{equation}

At the reference point, $u_{\mathbf a,0}$ is constant, so
$|\psi_\sigma(n)|^2=\mathbf1_{\{n=0\}}$. Lemma~\ref{cl:variation} therefore applies
with $v_{\mathbf a}=v_0'(\theta_{\mathbf a})$ and
$w_{\mathbf a}(W)=W(\theta_{\mathbf a})$.
Evaluation at the distinct points $\theta_{\mathbf a}$ maps
real trigonometric polynomials onto $\mathbb R^{2N}$.
The proof of that lemma shows that $\mathcal G_{m_0}^{(0)}$,
defined by \eqref{cl:generic} with these variations, has
complement a finite union of proper closed real linear
hyperplanes.

Fix $W\in\mathcal G_{m_0}^{(0)}$.
Use the same initial $B_0,D_0,E^{(0)}$ as in the proof of
Lemma~\ref{cl:KAM}, obtained from $C_{0,0}$.
The derivative bounds in \eqref{cl:uniform-comparison} also
hold for the change in the initial logarithmic remainder.
Along $E^{(0)}(\varepsilon)$, the old remainder is
$O(\delta^{12})$ and its $\varepsilon$ derivative is
$O(\delta^{11})$. Hence the weighted norm
\eqref{eq:collapse-parameter-norm} gives
$r_{0}\leqslant C\big(\delta^{12}+|\varepsilon_{*}|+\|u\|_{r}\big)
\leqslant3C\delta^{12}.$
whenever $0\leqslant\varepsilon_*\leqslant\delta^{12}$ and
$\|u\|_r\leqslant\delta^{12}$.
The initial phases are unchanged, and the energy derivative
bounds for the accumulated constants remain uniform.
Lemma~\ref{cl:KAM} thus gives the same $C\delta^3$ exclusion
estimate, uniformly in $\varepsilon_*,u$ and
$|\vartheta-\{m_{0}\alpha\}|\leqslant\delta$.

The symplectic pairing is independent of $\varepsilon_*$.
Proposition~\ref{energyvar} and Lemma~\ref{cl:collapse} give
a strict IDS crossing at each retained energy.
By \eqref{cl:scaled-ids}, its label is collapsed.
Taking the limsup over labels in each residue class gives
\begin{equation}\label{cl:scaled-cubic}
 \operatorname{Leb}\left\{\varepsilon\in[\delta,2\delta]:
 \varepsilon_*^{-1}(v_0+u+\varepsilon W)\in\mathcal S_N\right\}
 \geqslant\delta-C\delta^3
\end{equation}
for $0<\varepsilon_*\leqslant\delta^{12}$ and
$\|u\|_r\leqslant\delta^{12}$.
The constant is independent of $\varepsilon_*,u,\delta$.
Fundamental period one follows from
$|\widehat{v_0+u+\varepsilon W}(1)|
\geqslant\delta|w_1|-\delta^{12}>0$.

For the prescribed sequence $\varepsilon_*^j\to0$, set
$$
 S=\mathcal U\cap\limsup_{j\to\infty}
 \left\{v\in C_r^\omega(\mathbb T,\mathbb R):
       (\varepsilon_*^j)^{-1}v\in\mathcal S_N\right\}.
$$
This set is Borel.
Choose $\delta$ so that $C\delta^2<1$ and
$v_0+u+\varepsilon W\in\mathcal U$ whenever
$\|u\|_r\leqslant\delta^{12}$ and
$\varepsilon\in[\delta,2\delta]$.
For each such $u$, \eqref{cl:scaled-cubic} holds for all
sufficiently large $j$. Every tail union has measure at least
$\delta-C\delta^3$. Continuity of measure from above gives
$$
 \operatorname{Leb}\left\{\varepsilon\in[\delta,2\delta]:
       v_0+u+\varepsilon W\in S\right\}\geqslant\delta-C\delta^3.
$$
If $S$ were shy, choose a witnessing probability measure $\mu$
supported on $\{u:\|u\|_r\leqslant\delta^{12}\}$.
Fubini's theorem gives
$$
 \begin{aligned}
 0
 =\int_\delta^{2\delta}\mu(S-v_0-\varepsilon W)\,d\varepsilon
=\int\operatorname{Leb}\left\{
       \varepsilon\in[\delta,2\delta]:
       v_0+u+\varepsilon W\in S\right\}\,d\mu(u)
 \geqslant\delta-C\delta^3>0,
 \end{aligned}
$$
a contradiction.
\qedhere
\popQED

\subsection*{Proof of Theorem~\ref{thm:persistence-collapsed-labels}}
\pushQED{\qed}
For $m\notin N\mathbb Z$, let $\mathcal G_m$ be the set in
\eqref{cl:generic} with reference label $m$, and let
$
  \mathcal G_{\mathrm{all}}
  =\bigcap_{m\notin N\mathbb Z}\mathcal G_m .
$
The assertion about its complement follows from
Lemma~\ref{cl:variation}. Fix $W\in\mathcal G_{\mathrm{all}}$.
Fix for the moment $m\notin N\mathbb Z$.
Apply Lemmas~\ref{cl:KAM} and~\ref{cl:collapse} at
$\vartheta=\{m\alpha\}$. For sufficiently small $a_m>0$, set
$$
 \mathcal F_m=(0,a_m)\cap\bigcup_{j\geqslant0}
       \mathcal E_{2^{-j-1}a_m}(\{m\alpha\}).
$$
This is a Borel set. At every $\varepsilon\in\mathcal F_m$, label
$m$ is collapsed and the IDS has positive derivative at the
corresponding energy. Summing the measure bounds over the intervals
$[2^{-j-1}a_m,2^{-j}a_m]$ gives
$$
 \operatorname{Leb}\bigl((0,\delta)\setminus\mathcal F_m\bigr)
 \leqslant C_m\delta^3,\qquad 0<\delta\leqslant a_m.
$$
Enumerate $\mathbb Z\setminus N\mathbb Z=\{m_1,m_2,\ldots\}$.
Choose $b_j\downarrow0$ with $b_1<1$, $b_j<a_{m_j}$, and
$C_{m_j}b_j\leqslant 2^{-j}$. Set
$$
 \mathcal E_W=(0,b_1)\setminus
 \bigcup_{j\geq1}\bigl((0,b_j)\setminus\mathcal F_{m_j}\bigr).
$$
For $0<\delta<b_1$,
$$
 \operatorname{Leb}\bigl((0,\delta)\setminus\mathcal E_W\bigr)
 \leqslant\sum_{j\geq1} C_{m_j}\min\{\delta,b_j\}^3
 \leqslant\delta^2\sum_{j\geq1}C_{m_j}b_j
 \leqslant\delta^2.
$$
Finally, $\mathcal E_W\cap(0,b_j)\subset\mathcal F_{m_j}$.
Taking $\varepsilon_{m_j}(W)=b_j$ proves the claim. \qedhere
\popQED

\subsection*{Proof of Corollary~\ref{cor:linear-dispersion}}
\pushQED{\qed}
Fix an accepted energy $E_{\ell}=E_{\infty}$ from Lemma~\ref{cl:KAM}, keeping $\varepsilon$ fixed. Let
$$
Q=\big\langle\partial_{E}F_{\infty}(E_{\ell})\big\rangle_{\mathrm{diag}},\quad
V(E,x)=U_{\varepsilon,E-E_{*}}(x)B_{\infty}(x)e^{-(E-E_{\ell})X(x)},
$$
where $X$ is the solution used in \eqref{cl:energy-normal-form}. Each column of $V(E,x)$ solves the dual equation at energy $E$, and \eqref{cl:energy-normal-form} gives
$$
SV(E,x)=V(E,x+\alpha)D_{\infty}e^{(E-E_{\ell})Q}e^{O((E-E_{\ell})^{2})}.
$$
Write $V_{m,j}(E,x)=(V(E,x)\mathbf{e}_{j})_{m}$ and $\phi_{j}(x)=V_{0,j}(E_{\ell},x)$. Choose real analytic phase lifts $k_{j}(E)$, measured in radians, such that
$$
\big(D_{\infty}e^{(E-E_{\ell})Q}\big)_{jj}=e^{ik_{j}(E)},\qquad 1\leqslant j\leqslant2N.
$$
By \eqref{cl:energy-positivity},
$
(\widehat J_{2N})_{jj}Q_{jj}=\frac{1}{\lambda}\|\phi_{j}\|_{L^{2}(\mathbb{T})}^{2}>0.
$
Consequently, $\phi_{j}\not\equiv0$ and $k_{j}'(E_{\ell})=Q_{jj}/i\neq0$.

Iterating the frame relation a fixed positive or negative number of times gives
$$
\partial_{E}^{r}\big(V_{m,j}(E,x)-e^{imk_{j}(E)}V_{0,j}(E,x+m\alpha)\big)\bigg|_{E=E_{\ell}}=0,\qquad r=0,1,\quad m\in\mathbb{Z}.
$$
Indeed, the quadratic remainder has zero value and zero first energy derivative at $E_{\ell}$, so the iterated relation agrees to first order with that of the diagonal constant.

Consider the zeroth-site Bloch residual
$$
\begin{aligned}
\mathcal{R}_{j}(E,x)
&=\lambda e^{iNk_{j}(E)}V_{0,j}(E,x+N\alpha)+\lambda e^{-iNk_{j}(E)}V_{0,j}(E,x-N\alpha)\\
&\quad+\big(2\cos(2\pi x)-E\big)V_{0,j}(E,x)+\varepsilon\sum_{n\in\mathbb{Z}}w_{n}e^{ink_{j}(E)}V_{0,j}(E,x+n\alpha).
\end{aligned}
$$
The series is analytic in $E$ on a sufficiently small complex neighborhood of $E_{\ell}$, since $k_{j}(E_{\ell})$ is real and the coefficients $w_{n}$ decay exponentially. Moreover, the frame is analytic with values in $X_{a}$ and $\sum_{n\in\mathbb{Z}}|w_{n}|e^{a|n|}<\infty$, so its zeroth-site equation may also be differentiated term by term. The preceding coordinate identities and the equation satisfied by the columns of $V$ therefore give
$
\mathcal{R}_{j}(E_{\ell},x)=\partial_{E}\mathcal{R}_{j}(E_{\ell},x)=0.
$
Taylor's formula now yields
$
\|\mathcal{R}_{j}(E,\cdot)\|_{C^{0}}\leqslant C|E-E_{\ell}|^{2}.
$

Set $k_{j}=k_{j}(E_{\ell})$. Since $k_{j}'(E_{\ell})\neq0$, let $E_{j}(k)$ be the local real analytic inverse of $k_{j}(E)$, and define
$$
f_{j}(k,x)=V_{0,j}(E_{j}(k),x),\qquad u_{j,m}(k,x)=e^{imk}f_{j}(k,x+m\alpha).
$$
Then $f_{j}(k_{j},\cdot)=\phi_{j}\not\equiv0$, and direct substitution gives
$$
\big((\widehat H_{\varepsilon,x}-E_{j}(k))u_{j}(k,x)\big)_{m}=e^{imk}\mathcal{R}_{j}(E_{j}(k),x+m\alpha).
$$
 Since $|E_{j}(k)-E_{\ell}|\leqslant C|k-k_{j}|$, the residual is bounded by $C|k-k_{j}|^{2}$ uniformly in $x\in\mathbb{T}$ and $m\in\mathbb{Z}$. Thus $u_{j}$ is a first-order Bloch branch in the stated sense. Finally,
$$
E_{j}(k)=E_{\ell}+\frac{1}{k_{j}'(E_{\ell})}(k-k_{j})+O(|k-k_{j}|^{2}),
$$
which proves \eqref{eq:collapsed-linear-dispersion} with $s_{j}=1/k_{j}'(E_{\ell})\neq0$. The construction applies to each of the infinitely many accepted energies supplied in the proof of Theorem~\ref{thm:collapsed-gaps}. \qedhere
\popQED

\section*{Acknowledgements}
Part of the main results of this work was announced at the CIRM conference
``KAM Theory and PDEs''.
X.L. was partially supported by NSF DMS-2052899, DMS-
2155211, Simons 896624 and an AMS-Simons Travel Grant. Z.W. was supported
by NSFC grant (124B2011). Q.Z.
was supported by NSFC grant (12531006, 12526201). D.W., Z.W., Q.Z. were also supported by Nankai Zhide Foundation.

\section*{AI statement}
This work grew out of a decade of reflection and discussions with
experts, to whom we are deeply grateful. In the introduction and this statement, we have sought to
trace the origins and development of its ideas. We consider this
especially important in the age of AI, a conviction reinforced
after hearing A.~Avila's talk in Marseille.

The key ideas for Avila's problem were already in place after Q.Z. completed \cite{GJYZ26} (the polynomial case) and \cite{CCB26} (the large-coupling analytic case). Through discussions with an LLM, Q.Z. recognized that the method of \cite{CCB26}, rather than its result, could be extended to arbitrary $\lambda>1$. Discussions with the physicists Xiongjun Liu and Xinchi Zhou revealed dispersion as the underlying physical mechanism, while the LLM identified the key reference \cite{LL22}.

In 2015, during his postdoctoral stay in Irvine, Q.Z.\ began
working on the DTMP for nonzero analytic perturbations of the
second harmonic. There, he and S.~Jitomirskaya recognized how
resonances can open spectral gaps. Q.Z.\ subsequently discussed
the problem with L.~Ge and J.~You. In 2016, during Q.Z.'s visit
to the University of Toronto, M.~Goldstein posed the question
on the finiteness of collapsed gaps at large coupling. In 2024, discussions
with D.~Xu led us to recognize the key role of monotonicity,
as first suggested by A.~Avila. As a preparatory step, X.L.\
then began working on  Type~I DTMP.  These efforts led to \cite{LW,LXZ26,WXZ} and a separate manuscript developing a complete KAM scheme for $\mathrm{HSp}(J_{2N})$.\footnote{In preparation, after extensive discussions with R.~Krikorian.} Yet one case resisted all attempts at exclusion. In discussions with an LLM, the model raised objections to each proposed argument and suggested that this case might be unavoidable. The authors eventually recognized that this very case provides a
mechanism for gap collapse. We explained it by extending the
monotonicity argument of \cite{LXZ26} and connecting it with the
theory of positive symplectic paths, which X.L.\ learned through
discussions with Yiming Long, Jian Wang, and Li Wu.

All mathematical strategies originated with the human authors. The LLM also assisted with editing and the preparation of figures.

\end{document}